\documentclass[11pt]{article}
\usepackage{amsthm,amsfonts,latexsym,amssymb,amsopn,amscd,epic,eepic,epsfig}
\usepackage[all,2cell]{xy}
\UseAllTwocells
\usepackage{doc}
\newcommand{\N}{\mathbb{N}}
\newcommand{\Z}{\mathbb{Z}}

\title{Homotopy invariants of higher dimensional categories and concurrency in computer science}
\author{Philippe Gaucher\\
Institut de Recherche Math\'ematique Avanc\'ee  
\\U.L.P. et C.N.R.S.\\
7 rue Ren\'e Descartes\\ 
67084 Strasbourg Cedex\\
France
\\email : gaucher@irma.u-strasbg.fr}

\date{December 1999}

\newtheorem{thm}{Theorem}[section] 
\newtheorem{prop}[thm]{Proposition}

\newtheorem{cor}[thm]{Corollary}
\newtheorem{conj}[thm]{Conjecture}
\newtheorem{defn}{Definition}[section]

\newcommand{\be}{\begin{equation}}
\newcommand{\ee}{\end{equation}}
\newcommand{\bea}{\begin{eqnarray}}
\newcommand{\eea}{\end{eqnarray}}
\newcommand{\beas}{\begin{eqnarray*}}
\newcommand{\eeas}{\end{eqnarray*}}

\newcommand{\C}{\mathcal{C}}
\newcommand{\D}{\mathcal{D}}

\newcommand{\bd}{\begin{defn}}
\newcommand{\ed}{\end{defn}}
\newcommand{\bcd}{\begin{defn}}
\newcommand{\ecd}{\end{defn}}
\newcommand{\bp}{\begin{prop}}
\newcommand{\ep}{\end{prop}}
\newcommand{\bth}{\begin{thm}}
\renewcommand{\eth}{\end{thm}}
\newcommand{\bi}{\begin{enumerate}}
\newcommand{\ei}{\end{enumerate}}
\newcommand{\br}{\begin{rem}}
\newcommand{\er}{\end{rem}}
\newcommand{\bpf}{\begin{proof}}
\newcommand{\epf}{\end{proof}}

\newcommand{\p}\times

\newcommand{\iso}{\cong}

\newcommand{\U}{\Upsilon}
\newcommand{\ot}{\otimes}
\newcommand{\de}{\partial}

\renewcommand{\sec}[1]{\smash{\mathop{ 
{\hbox to 20mm{\leftarrowfill}}}\limits_{\scriptstyle#1}}}
\newcommand{\sect}[2]{\smash{\mathop{\hbox to 20mm{\rightarrowfill}}
\limits^{
\scriptstyle#1}_{\sec{#2}}}}

\newcommand{\fl}[1]{\ar@{->}[l]_{#1}}
\newcommand{\fr}[1]{\ar@{->}[r]^{#1}}
\newcommand{\fd}[1]{\ar@{->}[d]|{#1}}
\newcommand{\fu}[1]{\ar@{->}[u]^{#1}}
\newcommand{\f}[2]{\ar@{->}[#1]|{#2}}
\newcommand{\ff}[2]{\ar@2{->}[#1]_{#2}}

\newcommand{\coin}[2]{\begin{picture}(15,15)(0,0)
\put(2,0){\vector(1,0){10}}\put(13,0){#1}
\put(2,0){\vector(0,1){10}}\put(2,13){#2}
\end{picture}}
\begin{document}

\maketitle
%\abstract{bla}

\begin{abstract}

  The strict globular $\omega$-categories formalize the execution
  paths of a parallel automaton and the homotopies between them.  One
  associates to such (and any) $\omega$-category $\C$ three homology
  theories. The first one is called the globular homology. It contains
  the oriented loops of $\C$. The two other ones are called the
  negative (resp. positive ) corner homology. They contain in a
  certain manner the branching areas of execution paths or negative
  corners (resp. the merging areas of execution paths or positive
  corners) of $\C$.  Two natural linear maps called the negative
  (resp. the positive ) Hurewicz morphism from the globular homology
  to the negative (resp. positive) corner homology are constructed.
  We explain the reason why these constructions allow the
  reinterpretation of some geometric problems coming from computer
  science. 

\end{abstract}

\tableofcontents

\section{Introduction}

The use of geometric notions to describe the behaviour of concurrent
machines is certainly not new since {\em progress graphs} \cite{dij},
HDA \cite{Pratt}\cite{Gla} and simplicial models of \cite{MHNSSimple}
\cite{MHSRSet}. The purpose of this article is to provide a new
setting for the homotopy of execution paths in concurrent automata, in
order to improve the homological approach of \cite{HDA}.  We can point
out that some other approaches of this question already exist.  With
partially ordered topological spaces in \cite{HDA2}.  And with partial
posets in \cite{Sokolowski_1}.  Pratt already noticed that a good
structure to deal with execution paths and homotopies between them is
the structure of globular $\omega$-category \cite{Pratt}.  The aim of
this paper is threefold. First the use of the concept of globular
$\omega$-category to describe concurrent machines is justified on some
well-known examples and in a very informal way.  Secondly we associate
to every globular $\omega$-category three homology theories and two
natural maps between them. We explain why their content is interesting
for some geometric problems coming from computer science.  Thirdly, as
in algebraic topology, a notion of homotopic $\omega$-categories is
proposed and we prove that the preceding homology theories are
invariant with respect to it (only in a particular case for the corner
homologies).

Now here is the organization of the paper. In
Section~\ref{presentation}, we make precise the notion of paths and
homotopies between paths, homotopies between homotopies, etc... The
notion of globular $\omega$-category is recalled. The link between the
usual formalization of concurrent automata using cubical sets and the
new formalization by globular $\omega$-categories is explained. We
give in a very informal way the geometric and computer science meaning
of the homology theories which will be constructed in this paper. In
Section~\ref{glob_homo}, the globular homology of a globular
$\omega$-category is defined. Some examples of globular cycles are
given and the globular complex is related to a derived functor.
Section~\ref{corner} is devoted to the construction of the negative
and positive corner homologies of an $\omega$-category. Afterwards we
introduce in Section~\ref{filling} a technical tool to fill shells in
the cubical singular nerve of an $\omega$-category.  Next in
Section~\ref{corner_nerve}, this notion of filling of shell is used to
construct two families of connections on the cubical singular nerve of
an $\omega$-category. It allows us to prove that the corner homologies
are the homologies associated to two new simplicial nerves.
Afterwards we construct in Section~\ref{Hurewicz_oriente} the two
natural maps from the globular homology to the negative and positive
corner homologies.  In Section~\ref{ditop}, a notion of homotopy
equivalence of $\omega$-categories is proposed.  We will prove the
following property : \textit{Let $f$ and $g$ be two non
$1$-contracting $\omega$-functors from $\C$ to $\D$.  If $f$ and $g$
are homotopic, then for any natural number $n$,
$H_n^{gl}(f)=H_n^{gl}(g)$} (Theorem~\ref{globulaire_homotopie}). In a
very particular case, it is also possible to relate the homotopy of
paths in $\omega$-categories to the corner homology theories $H_n^{-}$
and $H_n^{+}$ (Theorem~\ref{orientee_homotopie}).  In
Section~\ref{perspective}, some conjectures and perspectives both in
mathematics and computer science are exposed. In
Section~\ref{morealgo}, we prove that the cubical singular nerve of
the free $\omega$-category generated by a cubical set $K$ is nothing
else but the free cubical $\omega$-category generated by $K$.  It
allows us to propose a direct construction of the globular homology
and of the corner homologies of a cubical set without using any
globular $\omega$-category.

\section{Presentation of the geometric ideas of this work}\label{presentation}

\subsection{Execution paths and homotopies between them in a very
  informal way}

A sequential machine (i.e. without concurrency) is a set of states,
also called $0$-transitions, and a set of $1$-transitions from a given
state to another one. A concurrent machine, like the previous one,
consists of a set of states and a set of $1$-transitions but has also
the capability of carrying out several $1$-transitions at the same
time.

\begin{figure}
\[
\xymatrix{{\beta} \fr{v'} &
{\delta} \\
{\alpha} \fu{u} \fr{v}& 
{\gamma} \ff{lu}{A}\fu{u'}}\]
\caption{A $2$-transition}
\label{A2transition}
\end{figure}
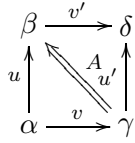

In Figure~\ref{A2transition}, if we work in cartesian coordinates in
such a way that $A=[0,1]\p[0,1]$ with $\alpha=(0,0)$ and
$\delta=(1,1)$, the set of continuous maps $(c_1,c_2)$ from $[0,1]$ to
$A$ such that $c_1(0)=c_2(0)=0$, $c_1(1)=c_2(1)=1$ and $t\leqslant t'$
implies $c_1(t)\leqslant c_1(t')$ and $c_2(t)\leqslant c_2(t')$
represents all the simultaneous possible executions of $u$ and $v$.
Coordinates represent the evolution of $u$ and $v$, that means the
local time taken to execute $u$ or $v$. If $(c_1,c_2)(]0,1[)$ is
entirely included in the interior of $A$, it is a ``true
parallelism''. If $(c_1,c_2)(]0,1[)$ is entirely included in the edge
of the square, that means that $u$ and $v$ are sequentially carried
out by the automaton. More generally, the concurrent execution of $n$
$1$-transitions can be represented by a $n$-cube. This is already
noticed for example in \cite{Pratt} and \cite{HDA}.

\bd\label{def_cubique} A \textit{cubical set} consists of a family of
sets $(K_n)_{n\geqslant 0}$, of a family of face maps
$\xymatrix@1{K_n\fr{\de_i^\alpha} &K_{n-1}}$ for $\alpha\in\{-,+\}$
and of a family of degeneracy maps
$\xymatrix@1{K_{n-1}\fr{\epsilon_i}&K_{n}}$ with $1\leqslant
i\leqslant n$ which satisfy the following relations
\begin{enumerate}
\item $\de_i^\alpha \de_j^\beta = \de_{j-1}^\beta \de_i^\alpha$    
for all $i<j\leqslant n$ and $\alpha,\beta\in\{-,+\}$
\item $\epsilon_i\epsilon_j=\epsilon_{j+1}\epsilon_i$   
for all $i\leqslant j\leqslant n$
\item $\de_i^\alpha \epsilon_j=\epsilon_{j-1}\de_i^\alpha$       for
  $i<j\leqslant n$ 
and $\alpha\in\{-,+\}$
\item $\de_i^\alpha \epsilon_j=\epsilon_{j}\de_{i-1}^\alpha$   for
  $i>j\leqslant n$ and 
$\alpha\in\{-,+\}$
\item $\de_i^\alpha \epsilon_i=Id$
\end{enumerate}
\ed

The corresponding category of cubical sets, with an obvious
definition of its morphisms, is isomorphic to the category of
presheaves $Sets^{\square^{op}}$ over a small category $\square$.
This latter can be described in a nice way as follows
\cite{Crans_Tensor_product}.  The objects of $\square$ are the sets
$\underline{n}=\{1,...,n\}$ where $n$ is a natural number greater or
equal than $1$ and an arrow $f$ from $\underline{n}$ to
$\underline{m}$ is a function $f^*$ from $\underline{m}$ to
$\underline{n}\cup \{-,+\}$ such that $f^*(k)\leqslant f^*(k')\in
\underline{n}$ implies $k\leqslant k'$ and $f^*(k)= f^*(k')\in
\underline{n}$ implies $k=k'$.

Let $X$ be a topological space. Let $[0,1]$ denote the interval
between $0$ and $1$. Set $K_n=C^0([0,1]^n,X)$ the set of continuous
maps from the $n$-box $[0,1]^n$ to $X$. Set \beas
&&\de_i^-(f)(x_1,\dots,x_p)=f(x_1,\dots,[0]_i,\dots,x_p)\\ 
&&\de_i^+(f)(x_1,\dots,x_p)=f(x_1,\dots,[1]_i,\dots,x_p)\\ 
&&\epsilon_i(f)(x_1,\dots,x_p)=f(x_1,\dots,\widehat{x_i},\dots,x_p)\\ 
\eeas Then we obtain a cubical set $K$ which is called the cubical
singular nerve of the topological space $X$.

In Figure~\ref{A2transition}, we call $A$ an homotopy between the two
$1$-paths $uv'$ and $vu'$. We call $2$-path an homotopy between two
$1$-paths and by induction on $n\geqslant 2$, we call $n$-path an
homotopy between two $(n-1)$-paths. This notion of homotopy is
different from the classical one in the sense that only the $1$-paths
$uv'$ and $vu'$ are homotopic and because the $1$-paths are oriented.
For example, still in Figure~\ref{A2transition}, $u$ is homotopic
neither with $u'$ nor with $v$ or $v'$.

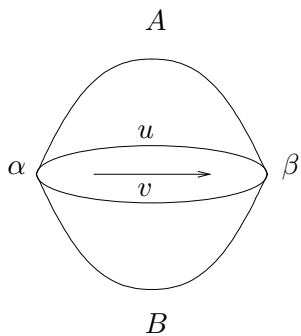
\begin{figure}
\begin{center}
\setlength{\unitlength}{0.0005in}
\begingroup\makeatletter\ifx\SetFigFont\undefined%
\gdef\SetFigFont#1#2#3#4#5{%
  \reset@font\fontsize{#1}{#2pt}%
  \fontfamily{#3}\fontseries{#4}\fontshape{#5}%
  \selectfont}%
\fi\endgroup%
{\renewcommand{\dashlinestretch}{30}
\begin{picture}(2942,3300)(0,-10)
\put(1500,1677){\ellipse{2400}{600}}
\path(900,1677)(2100,1677)
\path(1980.000,1647.000)(2100.000,1677.000)(1980.000,1707.000)
\path(300,1677) (318.909,1718.749)
        (337.671,1759.625)
        (356.293,1799.637)
        (374.782,1838.790)
        (393.144,1877.092)
        (411.387,1914.549)
        (447.542,1986.958)
        (483.302,2056.070)
        (518.720,2121.942)
        (553.853,2184.627)
        (588.756,2244.182)
        (623.482,2300.660)
        (658.088,2354.118)
        (692.627,2404.609)
        (727.155,2452.188)
        (761.728,2496.912)
        (796.399,2538.833)
        (831.224,2578.009)
        (866.258,2614.492)
        (937.170,2679.605)
        (1009.577,2734.610)
        (1083.918,2779.948)
        (1121.950,2799.128)
        (1160.631,2816.057)
        (1200.014,2830.788)
        (1240.156,2843.377)
        (1281.111,2853.879)
        (1322.933,2862.348)
        (1365.678,2868.840)
        (1409.401,2873.409)
        (1454.157,2876.111)
        (1500.000,2877.000)

\path(1500,2877)        (1545.843,2876.113)
        (1590.599,2873.413)
        (1634.322,2868.845)
        (1677.067,2862.355)
        (1718.889,2853.887)
        (1759.844,2843.386)
        (1799.986,2830.798)
        (1839.369,2816.068)
        (1878.050,2799.140)
        (1916.082,2779.961)
        (1990.423,2734.624)
        (2062.830,2679.620)
        (2133.742,2614.508)
        (2168.776,2578.024)
        (2203.601,2538.848)
        (2238.272,2496.926)
        (2272.845,2452.202)
        (2307.373,2404.622)
        (2341.912,2354.130)
        (2376.518,2300.672)
        (2411.244,2244.193)
        (2446.147,2184.638)
        (2481.280,2121.951)
        (2516.698,2056.078)
        (2552.458,1986.964)
        (2588.613,1914.554)
        (2606.856,1877.096)
        (2625.218,1838.794)
        (2643.707,1799.640)
        (2662.329,1759.627)
        (2681.091,1718.750)
        (2700.000,1677.000)

\path(300,1677) (318.909,1635.251)
        (337.671,1594.374)
        (356.293,1554.362)
        (374.782,1515.208)
        (393.144,1476.906)
        (411.387,1439.448)
        (447.542,1367.039)
        (483.302,1297.926)
        (518.720,1232.054)
        (553.853,1169.368)
        (588.756,1109.812)
        (623.482,1053.334)
        (658.088,999.876)
        (692.627,949.385)
        (727.155,901.805)
        (761.728,857.081)
        (796.399,815.159)
        (831.224,775.984)
        (866.258,739.500)
        (937.170,674.388)
        (1009.577,619.383)
        (1083.918,574.046)
        (1121.950,554.866)
        (1160.631,537.938)
        (1200.014,523.207)
        (1240.156,510.618)
        (1281.111,500.117)
        (1322.933,491.648)
        (1365.678,485.157)
        (1409.401,480.589)
        (1454.157,477.888)
        (1500.000,477.000)

\path(1500,477) (1545.843,477.888)
        (1590.599,480.589)
        (1634.322,485.157)
        (1677.067,491.648)
        (1718.889,500.117)
        (1759.844,510.618)
        (1799.986,523.207)
        (1839.369,537.938)
        (1878.050,554.866)
        (1916.082,574.046)
        (1990.423,619.383)
        (2062.830,674.388)
        (2133.742,739.500)
        (2168.776,775.984)
        (2203.601,815.159)
        (2238.272,857.081)
        (2272.845,901.805)
        (2307.373,949.385)
        (2341.912,999.876)
        (2376.518,1053.334)
        (2411.244,1109.812)
        (2446.147,1169.368)
        (2481.280,1232.054)
        (2516.698,1297.926)
        (2552.458,1367.039)
        (2588.613,1439.448)
        (2606.856,1476.906)
        (2625.218,1515.208)
        (2643.707,1554.362)
        (2662.329,1594.374)
        (2681.091,1635.251)
        (2700.000,1677.000)

\put(0,1677){\makebox(0,0)[lb]{\smash{{{\SetFigFont{11}{14.4}{\rmdefault}{\mddefault}{\updefault}$\alpha$}}}}}
\put(2850,1677){\makebox(0,0)[lb]{\smash{{{\SetFigFont{11}{14.4}{\rmdefault}{\mddefault}{\updefault}$\beta$}}}}}
\put(1425,3177){\makebox(0,0)[lb]{\smash{{{\SetFigFont{11}{14.4}{\rmdefault}{\mddefault}{\updefault}$A$}}}}}
\put(1425,27){\makebox(0,0)[lb]{\smash{{{\SetFigFont{11}{14.4}{\rmdefault}{\mddefault}{\updefault}$B$}}}}}
\put(1350,2052){\makebox(0,0)[lb]{\smash{{{\SetFigFont{11}{14.4}{\rmdefault}{\mddefault}{\updefault}$u$}}}}}
\put(1350,1452){\makebox(0,0)[lb]{\smash{{{\SetFigFont{11}{14.4}{\rmdefault}{\mddefault}{\updefault}$v$}}}}}
\end{picture}
}
\end{center}
\caption{A $3$-dimensional hole}
\label{3hole}
\end{figure}

In Figure~\ref{3hole}, there are an initial state $\alpha$, a final
state $\beta$, two $1$-transitions $u$ and $v$ and two $2$-transitions
or homotopies $A$ and $B$ between $u$ and $v$. Choosing an orientation
for $A$ and $B$, for example $s_1 A=s_1 B=u$ and $t_1 A=t_1 B=v$ ($s$
for source and $t$ for target), we see that $s_0 s_1 A = s_0 t_1 A =
\alpha$ and $t_0 s_1 A = t_0 t_1 A = \beta$. These are precisely the
globular equations which appear in the axioms of globular
$\omega$-categories.

\begin{figure}
\[\xymatrix{
{\alpha} \rruppertwocell<10>{A} \rrlowertwocell<-10>{B} \ar[rr]&& {\beta}
}
\]\caption{Composition of two $2$-morphims}
\label{composition-of-2trans}
\end{figure}
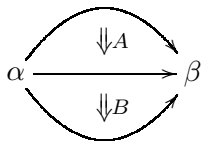

The concatenation yields an associative composition law on the set of
$1$-transitions of an $\omega$-category.  It turns out that there is
also a natural composition law on the set of $2$-transitions. In
Figure~\ref{composition-of-2trans}, with the convention of orientation
$t_1A=s_1B$, we can compose $A$ and $B$. Denote this composition by
$A*_1 B$. We see that $s_1(A*_1 B)=s_1 A$, $t_1(A*_1 B)=t_1 B$.  The
composition of higher dimensional morphisms must be associative
because it corresponds to the concatenation of the real execution
paths contained in $A$ and $B$.

Thus the geometric properties of transitions of concurrent machines
can be encoded in a structure of cubical set. And their associated
set of execution paths and homotopies between them have a natural
structure of globular $\omega$-category.  All these ideas already
appear in \cite{Pratt}. Pratt uses the term of $n$-complex which is in
fact nothing else but a small $n$-category. We use the notations of
\cite{Tensor_product} and \cite{oriental} for the following
definition. The following definition already appears in
\cite{Brown-Higgins0}

\bd\label{omega_categories} An \textit{$\omega$-category} is a set $A$
endowed with two families of maps $(s_n=d_n^-)_{n\geqslant 0}$ and
$(t_n=d_n^+)_{n\geqslant 0}$ from $A$ to $A$ and with a family of partially
defined 2-ary operations $(*_n)_{n\geqslant 0}$ where for any
$n\geqslant 0$, $*_n$ is a map from $\{(a,b)\in A\p A,
t_n(a)=s_n(b)\}$ to $A$ ($(a,b)$ being carried over $a *_n b$) which
satisfies the following axioms for all $\alpha$ and $\beta$ in
$\{-,+\}$ :

\begin{enumerate}
\item $d_m^\beta d_n^\alpha x=
\left\{\begin{CD}d_m^\beta x \hbox{  if $m<n$}\\  d_n^\alpha x \hbox{  if $m\geqslant n$}
\end{CD}\right.$
\item $s_n x *_n x= x *_n t_n x = x$
\item if  $x *_n y$ is well-defined, then  $s_n(x *_n y)=s_n x$, $t_n(x *_n y)=t_n y$
and for  $m\neq n$, $d_m^\alpha(x *_n y)=d_m^\alpha x *_n d_m^\alpha y$
\item as soon as the two members of the following equality exist, then 
$(x *_n y) *_n z= x *_n (y *_n z)$
\item if $m\neq n$ and if the two members of the equality make sense, then  
$(x *_n y)*_m (z*_n w)=(x *_m z) *_n (y *_m w)$
\item for any  $x$ in $A$, there exists a natural number $n$ such that $s_n x=t_n x=x$
(the smallest of these numbers is called the dimension of $x$ and is denoted
by $dim(x)$).
\end{enumerate}
\ed

We will sometimes use the notations $d^-_n:=s_n$ and $d^+_n=t_n$. If
$x$ is a morphism of an $\omega$-category $\C$, we call $s_n(x)$ the
$n$-source of $x$ and $t_n(x)$ the $n$-target of $x$.  The category of
all $\omega$-categories (with the obvious morphisms) is denoted by
$\omega Cat$. The corresponding morphisms are called
$\omega$-functors.

If $S$ is a set, the free abelian group generated by $S$ is denoted by
$\Z S$. By definition, an element of $\Z S$ is a formal linear
combination of elements of $S$.

\bd Let $\C$ be an $\omega$-category. Let $\C_n$ be the set of
$n$-dimensional morphisms of $\C$. Two $n$-morphisms $x$ and $y$ are
\textit{homotopic} if there exists $z\in \Z\C_{n+1}$ such that $s_n
z-t_n z= x-y$. This property is denoted by $x\sim y$.  \ed

If $x\sim y$, then the pair $(x,y)$ belongs to the reflexive,
symmetric and transitive closure of the binary relation generated by all
pairs $(s_n(u),t_n(u))$ where $u$ runs over $\C_{n+1}$.

\begin{figure}
\begin{center}
\setlength{\unitlength}{0.0005in}
\begingroup\makeatletter\ifx\SetFigFont\undefined%
\gdef\SetFigFont#1#2#3#4#5{%
  \reset@font\fontsize{#1}{#2pt}%
  \fontfamily{#3}\fontseries{#4}\fontshape{#5}%
  \selectfont}%
\fi\endgroup%
{\renewcommand{\dashlinestretch}{30}
\begin{picture}(2267,1974)(0,-10)
\path(225,12)(825,12)
\path(705.000,-18.000)(825.000,12.000)(705.000,42.000)
\path(825,12)(1425,12)
\path(1305.000,-18.000)(1425.000,12.000)(1305.000,42.000)
\path(1425,12)(2025,12)
\path(1905.000,-18.000)(2025.000,12.000)(1905.000,42.000)
\path(225,12)(225,612)
\path(255.000,492.000)(225.000,612.000)(195.000,492.000)
\path(825,12)(825,612)
\path(855.000,492.000)(825.000,612.000)(795.000,492.000)
\path(225,612)(225,1212)
\path(255.000,1092.000)(225.000,1212.000)(195.000,1092.000)
\path(825,612)(825,1212)
\path(855.000,1092.000)(825.000,1212.000)(795.000,1092.000)
\path(225,612)(825,612)
\path(705.000,582.000)(825.000,612.000)(705.000,642.000)
\path(825,612)(1425,612)
\path(1305.000,582.000)(1425.000,612.000)(1305.000,642.000)
\path(1425,12)(1425,612)
\path(1455.000,492.000)(1425.000,612.000)(1395.000,492.000)
\path(1425,612)(2025,612)
\path(1905.000,582.000)(2025.000,612.000)(1905.000,642.000)
\path(2025,12)(2025,612)
\path(2055.000,492.000)(2025.000,612.000)(1995.000,492.000)
\path(1425,612)(1425,1212)
\path(1455.000,1092.000)(1425.000,1212.000)(1395.000,1092.000)
\path(225,1212)(825,1212)
\path(705.000,1182.000)(825.000,1212.000)(705.000,1242.000)
\path(825,1212)(1425,1212)
\path(1305.000,1182.000)(1425.000,1212.000)(1305.000,1242.000)
\path(1425,1212)(2025,1212)
\path(1905.000,1182.000)(2025.000,1212.000)(1905.000,1242.000)
\path(2025,612)(2025,1212)
\path(2055.000,1092.000)(2025.000,1212.000)(1995.000,1092.000)
\path(225,1212)(225,1812)
\path(255.000,1692.000)(225.000,1812.000)(195.000,1692.000)
\path(225,1812)(825,1812)
\path(705.000,1782.000)(825.000,1812.000)(705.000,1842.000)
\path(825,1212)(825,1812)
\path(855.000,1692.000)(825.000,1812.000)(795.000,1692.000)
\path(825,1812)(1425,1812)
\path(1305.000,1782.000)(1425.000,1812.000)(1305.000,1842.000)
\path(1425,1212)(1425,1812)
\path(1455.000,1692.000)(1425.000,1812.000)(1395.000,1692.000)
\path(1425,1812)(2025,1812)
\path(1905.000,1782.000)(2025.000,1812.000)(1905.000,1842.000)
\path(2025,1212)(2025,1812)
\path(2055.000,1692.000)(2025.000,1812.000)(1995.000,1692.000)
\path(975,612)(825,837)
\path(1125,612)(825,1062)
\path(1275,612)(900,1212)
\path(1425,687)(1050,1212)
\path(1425,987)(1275,1212)
\put(0,237){\makebox(0,0)[lb]{\smash{{{\SetFigFont{11}{14.4}{\rmdefault}{\mddefault}{\updefault}$u$}}}}}
\put(525,687){\makebox(0,0)[lb]{\smash{{{\SetFigFont{11}{14.4}{\rmdefault}{\mddefault}{\updefault}$v$}}}}}
\put(750,912){\makebox(0,0)[lb]{\smash{{{\SetFigFont{11}{14.4}{\rmdefault}{\mddefault}{\updefault}$w$}}}}}
\put(1125,1287){\makebox(0,0)[lb]{\smash{{{\SetFigFont{11}{14.4}{\rmdefault}{\mddefault}{\updefault}$x$}}}}}
\put(1500,1437){\makebox(0,0)[lb]{\smash{{{\SetFigFont{11}{14.4}{\rmdefault}{\mddefault}{\updefault}$y$}}}}}
\put(1650,1887){\makebox(0,0)[lb]{\smash{{{\SetFigFont{11}{14.4}{\rmdefault}{\mddefault}{\updefault}$z$}}}}}
\put(0,837){\makebox(0,0)[lb]{\smash{{{\SetFigFont{11}{14.4}{\rmdefault}{\mddefault}{\updefault}$a$}}}}}
\put(525,1362){\makebox(0,0)[lb]{\smash{{{\SetFigFont{11}{14.4}{\rmdefault}{\mddefault}{\updefault}$b$}}}}}
\put(1725,1287){\makebox(0,0)[lb]{\smash{{{\SetFigFont{11}{14.4}{\rmdefault}{\mddefault}{\updefault}$c$}}}}}
\put(2175,1437){\makebox(0,0)[lb]{\smash{{{\SetFigFont{11}{14.4}{\rmdefault}{\mddefault}{\updefault}$d$}}}}}
\put(450,87){\makebox(0,0)[lb]{\smash{{{\SetFigFont{11}{14.4}{\rmdefault}{\mddefault}{\updefault}$e$}}}}}
\put(1125,87){\makebox(0,0)[lb]{\smash{{{\SetFigFont{11}{14.4}{\rmdefault}{\mddefault}{\updefault}$f$}}}}}
\put(1725,87){\makebox(0,0)[lb]{\smash{{{\SetFigFont{11}{14.4}{\rmdefault}{\mddefault}{\updefault}$g$}}}}}
\put(2175,237){\makebox(0,0)[lb]{\smash{{{\SetFigFont{11}{14.4}{\rmdefault}{\mddefault}{\updefault}$h$}}}}}
\put(2175,837){\makebox(0,0)[lb]{\smash{{{\SetFigFont{11}{14.4}{\rmdefault}{\mddefault}{\updefault}$i$}}}}}
\end{picture}
}
\beas
&&\gamma_1=u *_0 v *_0 w *_0 x
*_0 y *_0 z\\
&&\gamma_2=u *_0 a *_0 b *_0 x *_0 c *_0 d\\
&&\gamma_3=u *_0 v *_0 w
*_0 x *_0 c *_0 d\\
&&\gamma_4= e *_0 f *_0 g *_0 h *_0 i *_0 d
\eeas
\end{center}

\caption{Example of distributed database}
\label{trou}
\end{figure}

Figure~\ref{trou} is a very simple example of a distributed database.
The hole in the middle corresponds to a mutual exclusion.  See
\cite{FGR} for a complete treatment.  The two $1$-paths $\gamma_1$ and
$\gamma_2$ are homotopic because there exists a $2$-morphism between
$\gamma_1$ and $\gamma_3$ and another one between $\gamma_2$ and
$\gamma_3$. On the other hand none of the previous three $1$-paths is
homotopic to $\gamma_4$ because of the oriented hole in the middle.

\subsection{The free $\omega$-category generated by a cubical set}

How may we mathematically associate to every cubical set $K$ its
corresponding set of execution paths and higher dimensional homotopies
?  The link between the two formalizations is as follows.  We need to
describe precisely the free $\omega$-category associated to each
$n$-cube of $K$. In a very informal way, it consists of seeing the
faces of a $n$-cube as words of length $n$ in the alphabet
$\{-,0,+\}$. The term $0_n$ means $0$ $n$ times, i.e. the interior of
the $n$-cube. And we say that the face $k_1\dots k_n$ is at the source
of $x$ if $k_i\neq 0$ implies $k_i=(-1)^i$ and we say that $k_1\dots
k_n$ is at the target of $x$ if $k_i\neq 0$ implies $k_i=(-1)^{i+1}$.
We will make precise the construction of $I^n$ in
Section~\ref{construction-of-In}.  Once this is done, it suffices to
paste the $\omega$-categories associated to every $n$-cube of $K$ in
the same way that they are pasted in $K$.  More concretely, every
cubical set $K$ is in a canonical way the direct limit of the
elementary $n$-cubes included in it. This is due to the fact that any
functor from a small category to the category $Sets$ of sets is a
canonical direct limit of representable functors, the set of those
functors being dense in the set of all set-valued functors. More
precisely, we have $K=\int^{\underline{n}\in \square} K_n
.\square(-,\underline{n})$ where the integral sign is the coend
construction \cite{cat} and $K_n .\square(-,\underline{n})$ means the
sum of ``cardinal of $K_n$'' copies of $\square(-,\underline{n})$.  So
$F(K)=\int^{\underline{n}\in \square} K_n . I^n$ is a
$\omega$-category containing as $1$-morphisms all arrows of $K$ and
all possible compositions of these arrows, as $2$-morphisms all
homotopies between the execution paths, etc...

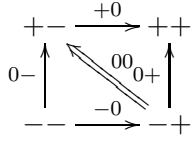
\begin{figure}
\[\xymatrix{{+-} \fr{+0} &
{++} \\
{--} \fu{0-} \fr{-0}& 
{-+} \ff{lu}{00}\fu{0+}}\]
\caption{The $\omega$-category $I^2$}
\label{I2}
\end{figure}

Consider for example the $2$-cube of Figure~\ref{I2}. The $2$-face
$00$ is oriented from the side $\{-0,0+\}$ to the side $\{0-,+0\}$.
The corresponding $\omega$-category $I^2$ contains all possible
compositions of faces of the $2$-cube. Therefore, as set, we have
$$I^2=\{--,-+,+-,++,-0,0-,+0,0+,-0 *_0 0+, 0- *_0 +0, 00\}.$$

In Figure~\ref{composition-3-squares}, the three $2$-morphisms $A$,
$B$ and $C$ are not drawn and are supposed to be oriented to the north
west. The composition of the three squares $A$, $B$ and $C$ is equal
to $$(A*_0 w *_0 x)*_1 (v' *_0 ((u'*_0 C)*_1(B*_0 w')))$$
in the free $\omega$-category generated by this cubical set.

\begin{figure}
\begin{center}
\setlength{\unitlength}{0.0005in}
\begingroup\makeatletter\ifx\SetFigFont\undefined%
\gdef\SetFigFont#1#2#3#4#5{%
  \reset@font\fontsize{#1}{#2pt}%
  \fontfamily{#3}\fontseries{#4}\fontshape{#5}%
  \selectfont}%
\fi\endgroup%
{\renewcommand{\dashlinestretch}{30}
\begin{picture}(1667,1716)(0,-10)
\path(225,243)(825,243)
\path(705.000,213.000)(825.000,243.000)(705.000,273.000)
\path(825,243)(825,843)
\path(855.000,723.000)(825.000,843.000)(795.000,723.000)
\path(825,843)(1425,843)
\path(1305.000,813.000)(1425.000,843.000)(1305.000,873.000)
\path(1425,843)(1425,1443)
\path(1455.000,1323.000)(1425.000,1443.000)(1395.000,1323.000)
\path(825,1443)(1425,1443)
\path(1305.000,1413.000)(1425.000,1443.000)(1305.000,1473.000)
\path(825,843)(825,1443)
\path(855.000,1323.000)(825.000,1443.000)(795.000,1323.000)
\path(225,843)(825,843)
\path(705.000,813.000)(825.000,843.000)(705.000,873.000)
\path(225,243)(225,843)
\path(255.000,723.000)(225.000,843.000)(195.000,723.000)
\path(225,843)(225,1443)
\path(255.000,1323.000)(225.000,1443.000)(195.000,1323.000)
\path(225,1443)(825,1443)
\path(705.000,1413.000)(825.000,1443.000)(705.000,1473.000)
\path(225,243)  (284.975,294.605)
        (327.622,334.462)
        (375.000,393.000)

\path(375,393)  (398.004,458.632)
        (406.619,499.652)
        (414.281,543.015)
        (421.721,586.376)
        (429.668,627.391)
        (450.000,693.000)

\path(450,693)  (473.560,736.003)
        (504.629,785.580)
        (541.443,838.997)
        (582.240,893.520)
        (625.257,946.416)
        (668.731,994.950)
        (710.900,1036.389)
        (750.000,1068.000)

\path(750,1068) (793.287,1091.465)
        (847.642,1111.586)
        (909.383,1129.362)
        (974.827,1145.790)
        (1040.294,1161.870)
        (1102.100,1178.599)
        (1156.562,1196.976)
        (1200.000,1218.000)

\path(1200,1218)        (1239.118,1246.543)
        (1286.141,1290.041)
        (1346.344,1353.769)
        (1383.036,1394.867)
        (1425.000,1443.000)

\path(1369.242,1332.587)(1425.000,1443.000)(1323.843,1371.816)
\put(450,468){\makebox(0,0)[lb]{\smash{{{\SetFigFont{11}{14.4}{\rmdefault}{\mddefault}{\updefault}A}}}}}
\put(450,1068){\makebox(0,0)[lb]{\smash{{{\SetFigFont{11}{14.4}{\rmdefault}{\mddefault}{\updefault}B}}}}}
\put(1050,1068){\makebox(0,0)[lb]{\smash{{{\SetFigFont{11}{14.4}{\rmdefault}{\mddefault}{\updefault}C}}}}}
\put(450,18){\makebox(0,0)[lb]{\smash{{{\SetFigFont{11}{14.4}{\rmdefault}{\mddefault}{\updefault}u}}}}}
\put(900,468){\makebox(0,0)[lb]{\smash{{{\SetFigFont{11}{14.4}{\rmdefault}{\mddefault}{\updefault}v}}}}}
\put(1125,618){\makebox(0,0)[lb]{\smash{{{\SetFigFont{11}{14.4}{\rmdefault}{\mddefault}{\updefault}w}}}}}
\put(1575,1068){\makebox(0,0)[lb]{\smash{{{\SetFigFont{11}{14.4}{\rmdefault}{\mddefault}{\updefault}x}}}}}
\put(0,468){\makebox(0,0)[lb]{\smash{{{\SetFigFont{11}{14.4}{\rmdefault}{\mddefault}{\updefault}v'}}}}}
\put(0,1068){\makebox(0,0)[lb]{\smash{{{\SetFigFont{11}{14.4}{\rmdefault}{\mddefault}{\updefault}x'}}}}}
\put(450,1593){\makebox(0,0)[lb]{\smash{{{\SetFigFont{11}{14.4}{\rmdefault}{\mddefault}{\updefault}u''}}}}}
\put(1050,1593){\makebox(0,0)[lb]{\smash{{{\SetFigFont{11}{14.4}{\rmdefault}{\mddefault}{\updefault}w'}}}}}
\put(450,843){\makebox(0,0)[lb]{\smash{{{\SetFigFont{11}{14.4}{\rmdefault}{\mddefault}{\updefault}u'}}}}}
\end{picture}
}
\end{center}
$$(A*_0 w *_0 x)*_1 (v' *_0 ((u'*_0 C)*_1(B*_0 w')))$$
\caption{Composition of three squares}
\label{composition-3-squares}
\end{figure}

The map $F$ induces a functor from the category of cubical sets to the
category of $\omega$-categories.  Indeed this is the left Kan
extension of the functor $Q$ from $\square$ to $\omega Cat$ defined as
follows \cite{cat}.  It maps $\underline{n}$ to $I^n$. Let
$\epsilon_i$ be the surjective morphism from $\underline{n}$ to
$\underline{n-1}$ for $1\leqslant i\leqslant n$ defined by
$(\epsilon_i)^*(l)=l$ if $l<i$ and $(\epsilon_i)^*(l)=l+1$.  Let
$\de_i^\alpha$ be the injective morphism from $\underline{n-1}$ to
$\underline{n}$ for $1\leqslant i\leqslant n$ and for $\alpha=\pm$
defined by $(\de_i^\alpha)^*(l)=l$ if $l<i$,
$(\de_i^\alpha)^*(l)=\alpha$ for $l=i$ and $(\de_i^\alpha)^*(l)=l-1$
for $l>i$. Then any morphism of $\square$ is a composition of
$\epsilon_i$ and of $\de_i^\alpha$. And $Q$ is the unique functor
which maps $\epsilon_i$ to $\epsilon_i$ and $\de_i^\alpha$ to
$\de_i^\alpha$ \cite{Crans_Tensor_product}.  This way, the notion of
$\omega$-category can be understood as a generalization of the notion
of cubical set. Every cubical set can be seen as an $\omega$-category.
The converse is false. We will see in Section~\ref{ditop} why this
categoric setting is very well adapted for the development of an
analogue of algebraic topology in the computer science framework.

\subsection{$\omega$-categories up to homotopy}

Now we want to give, in a very informal way, examples of
$\omega$-categories which have the same set of execution paths up to
homotopy and to explain the potential interest of this notion. We will
propose a definition of homotopic $\omega$-categories in
Definition~\ref{def_morphisme_homotope} and \ref{def_cat_homotope}.

\begin{figure}
\[\xymatrix{
\ar@/^20pt/[rr]^{A}\ar@/_20pt/[rr]_{B}&&
}
\]
\caption{The oriented globe $G_1$}
\label{G_1}
\end{figure}
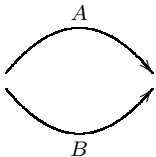

Up to path homotopy, the $\omega$-category of Figure~\ref{trou} must
be the same as the $\omega$-category of Figure~\ref{G_1} because there
are only two execution paths up to homotopy in each case.

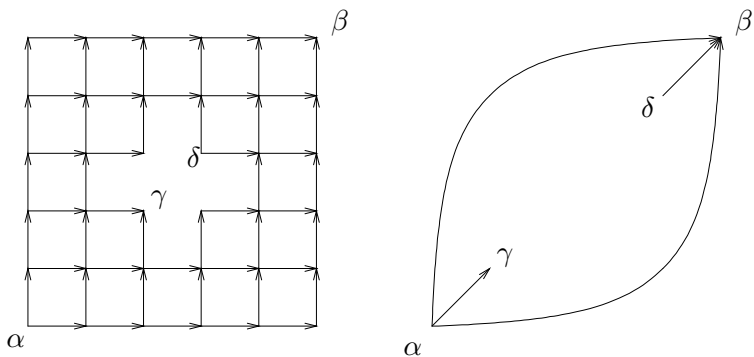
\begin{figure}
\begin{center}
\setlength{\unitlength}{0.0005in}
\begingroup\makeatletter\ifx\SetFigFont\undefined%
\gdef\SetFigFont#1#2#3#4#5{%
  \reset@font\fontsize{#1}{#2pt}%
  \fontfamily{#3}\fontseries{#4}\fontshape{#5}%
  \selectfont}%
\fi\endgroup%
{\renewcommand{\dashlinestretch}{30}
\begin{picture}(7667,3516)(0,-10)
\path(825,1518)(1425,1518)
\path(1305.000,1488.000)(1425.000,1518.000)(1305.000,1548.000)
\path(825,1518)(825,2118)
\path(855.000,1998.000)(825.000,2118.000)(795.000,1998.000)
\path(825,2118)(1425,2118)
\path(1305.000,2088.000)(1425.000,2118.000)(1305.000,2148.000)
\path(1425,2118)(1425,2718)
\path(1455.000,2598.000)(1425.000,2718.000)(1395.000,2598.000)
\path(1425,2718)(2025,2718)
\path(1905.000,2688.000)(2025.000,2718.000)(1905.000,2748.000)
\path(2025,2118)(2025,2718)
\path(2055.000,2598.000)(2025.000,2718.000)(1995.000,2598.000)
\path(2025,2118)(2625,2118)
\path(2505.000,2088.000)(2625.000,2118.000)(2505.000,2148.000)
\path(2625,1518)(2625,2118)
\path(2655.000,1998.000)(2625.000,2118.000)(2595.000,1998.000)
\path(2025,1518)(2625,1518)
\path(2505.000,1488.000)(2625.000,1518.000)(2505.000,1548.000)
\path(2025,918)(2025,1518)
\path(2055.000,1398.000)(2025.000,1518.000)(1995.000,1398.000)
\path(1425,918)(1425,1518)
\path(1455.000,1398.000)(1425.000,1518.000)(1395.000,1398.000)
\path(1425,918)(2025,918)
\path(1905.000,888.000)(2025.000,918.000)(1905.000,948.000)
\path(225,318)(825,318)
\path(705.000,288.000)(825.000,318.000)(705.000,348.000)
\path(225,318)(225,918)
\path(255.000,798.000)(225.000,918.000)(195.000,798.000)
\path(225,918)(825,918)
\path(705.000,888.000)(825.000,918.000)(705.000,948.000)
\path(225,918)(225,1518)
\path(255.000,1398.000)(225.000,1518.000)(195.000,1398.000)
\path(225,1518)(825,1518)
\path(705.000,1488.000)(825.000,1518.000)(705.000,1548.000)
\path(225,1518)(225,2118)
\path(255.000,1998.000)(225.000,2118.000)(195.000,1998.000)
\path(225,2118)(825,2118)
\path(705.000,2088.000)(825.000,2118.000)(705.000,2148.000)
\path(225,2118)(225,2718)
\path(255.000,2598.000)(225.000,2718.000)(195.000,2598.000)
\path(825,2118)(825,2718)
\path(855.000,2598.000)(825.000,2718.000)(795.000,2598.000)
\path(225,2718)(825,2718)
\path(705.000,2688.000)(825.000,2718.000)(705.000,2748.000)
\path(225,2718)(225,3318)
\path(255.000,3198.000)(225.000,3318.000)(195.000,3198.000)
\path(225,3318)(825,3318)
\path(705.000,3288.000)(825.000,3318.000)(705.000,3348.000)
\path(825,2718)(1425,2718)
\path(1305.000,2688.000)(1425.000,2718.000)(1305.000,2748.000)
\path(825,2718)(825,3318)
\path(855.000,3198.000)(825.000,3318.000)(795.000,3198.000)
\path(825,3318)(1425,3318)
\path(1305.000,3288.000)(1425.000,3318.000)(1305.000,3348.000)
\path(1425,2718)(1425,3318)
\path(1455.000,3198.000)(1425.000,3318.000)(1395.000,3198.000)
\path(1425,3318)(2025,3318)
\path(1905.000,3288.000)(2025.000,3318.000)(1905.000,3348.000)
\path(2025,2718)(2025,3318)
\path(2055.000,3198.000)(2025.000,3318.000)(1995.000,3198.000)
\path(2025,2718)(2625,2718)
\path(2505.000,2688.000)(2625.000,2718.000)(2505.000,2748.000)
\path(2625,2118)(2625,2718)
\path(2655.000,2598.000)(2625.000,2718.000)(2595.000,2598.000)
\path(2025,3318)(2625,3318)
\path(2505.000,3288.000)(2625.000,3318.000)(2505.000,3348.000)
\path(2625,2718)(2625,3318)
\path(2655.000,3198.000)(2625.000,3318.000)(2595.000,3198.000)
\path(2625,2718)(3225,2718)
\path(3105.000,2688.000)(3225.000,2718.000)(3105.000,2748.000)
\path(3225,2718)(3225,3318)
\path(3255.000,3198.000)(3225.000,3318.000)(3195.000,3198.000)
\path(2625,3318)(3225,3318)
\path(3105.000,3288.000)(3225.000,3318.000)(3105.000,3348.000)
\path(2625,2118)(3225,2118)
\path(3105.000,2088.000)(3225.000,2118.000)(3105.000,2148.000)
\path(3225,2118)(3225,2718)
\path(3255.000,2598.000)(3225.000,2718.000)(3195.000,2598.000)
\path(2625,1518)(3225,1518)
\path(3105.000,1488.000)(3225.000,1518.000)(3105.000,1548.000)
\path(3225,1518)(3225,2118)
\path(3255.000,1998.000)(3225.000,2118.000)(3195.000,1998.000)
\path(2025,918)(2625,918)
\path(2505.000,888.000)(2625.000,918.000)(2505.000,948.000)
\path(2625,918)(2625,1518)
\path(2655.000,1398.000)(2625.000,1518.000)(2595.000,1398.000)
\path(2625,918)(3225,918)
\path(3105.000,888.000)(3225.000,918.000)(3105.000,948.000)
\drawline(3225,918)(3225,918)
\path(3225,918)(3225,1518)
\path(3255.000,1398.000)(3225.000,1518.000)(3195.000,1398.000)
\path(825,918)(825,1518)
\path(855.000,1398.000)(825.000,1518.000)(795.000,1398.000)
\path(825,918)(1425,918)
\path(1305.000,888.000)(1425.000,918.000)(1305.000,948.000)
\path(825,318)(825,918)
\path(855.000,798.000)(825.000,918.000)(795.000,798.000)
\path(825,318)(1425,318)
\path(1305.000,288.000)(1425.000,318.000)(1305.000,348.000)
\path(1425,318)(1425,918)
\path(1455.000,798.000)(1425.000,918.000)(1395.000,798.000)
\path(1425,318)(2025,318)
\path(1905.000,288.000)(2025.000,318.000)(1905.000,348.000)
\path(2025,318)(2025,918)
\path(2055.000,798.000)(2025.000,918.000)(1995.000,798.000)
\path(2025,318)(2625,318)
\path(2505.000,288.000)(2625.000,318.000)(2505.000,348.000)
\path(2625,318)(2625,918)
\path(2655.000,798.000)(2625.000,918.000)(2595.000,798.000)
\path(2625,318)(3225,318)
\path(3105.000,288.000)(3225.000,318.000)(3105.000,348.000)
\path(3225,318)(3225,918)
\path(3255.000,798.000)(3225.000,918.000)(3195.000,798.000)
\path(4425,318)(5025,918)
\path(4961.360,811.934)(5025.000,918.000)(4918.934,854.360)
\path(6825,2718)(7425,3318)
\path(7361.360,3211.934)(7425.000,3318.000)(7318.934,3254.360)
\path(4425,318) (4427.339,382.962)
        (4429.830,446.769)
        (4432.477,509.434)
        (4435.283,570.971)
        (4438.252,631.393)
        (4441.387,690.715)
        (4444.692,748.949)
        (4448.169,806.111)
        (4451.824,862.213)
        (4455.658,917.269)
        (4459.675,971.293)
        (4463.879,1024.299)
        (4468.273,1076.300)
        (4472.861,1127.310)
        (4477.646,1177.343)
        (4482.632,1226.413)
        (4487.821,1274.533)
        (4493.218,1321.718)
        (4498.826,1367.980)
        (4504.647,1413.333)
        (4510.687,1457.792)
        (4516.947,1501.370)
        (4523.432,1544.081)
        (4530.145,1585.938)
        (4537.090,1626.956)
        (4544.269,1667.147)
        (4551.686,1706.526)
        (4559.345,1745.107)
        (4575.403,1819.927)
        (4592.468,1891.717)
        (4610.568,1960.588)
        (4629.731,2026.650)
        (4649.985,2090.011)
        (4671.356,2150.782)
        (4693.873,2209.073)
        (4717.562,2264.994)
        (4742.452,2318.654)
        (4768.569,2370.163)
        (4795.942,2419.632)
        (4824.597,2467.170)
        (4854.562,2512.886)
        (4885.865,2556.892)
        (4918.532,2599.296)
        (4952.592,2640.209)
        (4988.073,2679.740)
        (5025.000,2718.000)

\path(5025,2718)        (5063.260,2754.927)
        (5102.791,2790.408)
        (5143.704,2824.468)
        (5186.108,2857.135)
        (5230.114,2888.438)
        (5275.830,2918.403)
        (5323.368,2947.058)
        (5372.837,2974.431)
        (5424.346,3000.548)
        (5478.006,3025.438)
        (5533.927,3049.127)
        (5592.218,3071.644)
        (5652.989,3093.015)
        (5716.350,3113.269)
        (5782.412,3132.432)
        (5851.282,3150.533)
        (5923.073,3167.597)
        (5997.893,3183.655)
        (6036.474,3191.314)
        (6075.853,3198.731)
        (6116.044,3205.910)
        (6157.062,3212.855)
        (6198.919,3219.568)
        (6241.630,3226.053)
        (6285.208,3232.313)
        (6329.667,3238.353)
        (6375.020,3244.174)
        (6421.282,3249.782)
        (6468.467,3255.179)
        (6516.587,3260.368)
        (6565.657,3265.354)
        (6615.690,3270.139)
        (6666.700,3274.727)
        (6718.701,3279.121)
        (6771.707,3283.325)
        (6825.731,3287.342)
        (6880.787,3291.176)
        (6936.889,3294.831)
        (6994.051,3298.308)
        (7052.285,3301.613)
        (7111.607,3304.748)
        (7172.029,3307.717)
        (7233.566,3310.523)
        (7296.231,3313.170)
        (7360.038,3315.661)
        (7425.000,3318.000)

\path(7306.107,3283.876)(7425.000,3318.000)(7304.036,3343.840)
\path(4425,318) (4489.963,320.339)
        (4553.770,322.830)
        (4616.436,325.477)
        (4677.973,328.283)
        (4738.395,331.252)
        (4797.717,334.387)
        (4855.952,337.692)
        (4913.114,341.169)
        (4969.216,344.824)
        (5024.273,348.658)
        (5078.297,352.675)
        (5131.303,356.879)
        (5183.305,361.273)
        (5234.315,365.861)
        (5284.349,370.646)
        (5333.419,375.632)
        (5381.539,380.821)
        (5428.724,386.218)
        (5474.986,391.826)
        (5520.340,397.647)
        (5564.799,403.687)
        (5608.377,409.947)
        (5651.088,416.432)
        (5692.945,423.145)
        (5733.963,430.090)
        (5774.154,437.269)
        (5813.534,444.686)
        (5852.114,452.345)
        (5926.934,468.403)
        (5998.725,485.467)
        (6067.596,503.568)
        (6133.657,522.731)
        (6197.018,542.985)
        (6257.789,564.356)
        (6316.080,586.873)
        (6372.000,610.562)
        (6425.660,635.452)
        (6477.169,661.569)
        (6526.637,688.942)
        (6574.174,717.597)
        (6619.890,747.562)
        (6663.895,778.865)
        (6706.299,811.532)
        (6747.211,845.592)
        (6786.741,881.073)
        (6825.000,918.000)

\path(6825,918) (6861.927,956.260)
        (6897.406,995.791)
        (6931.465,1036.704)
        (6964.132,1079.108)
        (6995.434,1123.114)
        (7025.399,1168.830)
        (7054.053,1216.368)
        (7081.425,1265.837)
        (7107.542,1317.346)
        (7132.431,1371.006)
        (7156.120,1426.927)
        (7178.637,1485.218)
        (7200.008,1545.989)
        (7220.261,1609.350)
        (7239.425,1675.412)
        (7257.525,1744.283)
        (7274.590,1816.073)
        (7290.647,1890.893)
        (7298.306,1929.474)
        (7305.724,1968.853)
        (7312.903,2009.044)
        (7319.848,2050.062)
        (7326.561,2091.919)
        (7333.046,2134.630)
        (7339.307,2178.208)
        (7345.346,2222.667)
        (7351.168,2268.020)
        (7356.776,2314.282)
        (7362.173,2361.467)
        (7367.363,2409.587)
        (7372.348,2458.657)
        (7377.133,2508.690)
        (7381.722,2559.700)
        (7386.116,2611.701)
        (7390.321,2664.707)
        (7394.338,2718.731)
        (7398.173,2773.787)
        (7401.827,2829.889)
        (7405.305,2887.051)
        (7408.610,2945.285)
        (7411.746,3004.607)
        (7414.715,3065.029)
        (7417.522,3126.566)
        (7420.169,3189.231)
        (7422.661,3253.038)
        (7425.000,3318.000)

\path(7450.840,3197.036)(7425.000,3318.000)(7390.875,3199.107)
\put(0,93){\makebox(0,0)[lb]{\smash{{{\SetFigFont{11}{14.4}{\rmdefault}{\mddefault}{\updefault}$\alpha$}}}}}
\put(3375,3393){\makebox(0,0)[lb]{\smash{{{\SetFigFont{11}{14.4}{\rmdefault}{\mddefault}{\updefault}$\beta$}}}}}
\put(1500,1593){\makebox(0,0)[lb]{\smash{{{\SetFigFont{11}{14.4}{\rmdefault}{\mddefault}{\updefault}$\gamma$}}}}}
\put(1875,1968){\makebox(0,0)[lb]{\smash{{{\SetFigFont{11}{14.4}{\rmdefault}{\mddefault}{\updefault}$\delta$}}}}}
\put(4125,18){\makebox(0,0)[lb]{\smash{{{\SetFigFont{11}{14.4}{\rmdefault}{\mddefault}{\updefault}$\alpha$}}}}}
\put(7575,3393){\makebox(0,0)[lb]{\smash{{{\SetFigFont{11}{14.4}{\rmdefault}{\mddefault}{\updefault}$\beta$}}}}}
\put(5100,993){\makebox(0,0)[lb]{\smash{{{\SetFigFont{11}{14.4}{\rmdefault}{\mddefault}{\updefault}$\gamma$}}}}}
\put(6600,2493){\makebox(0,0)[lb]{\smash{{{\SetFigFont{11}{14.4}{\rmdefault}{\mddefault}{\updefault}$\delta$}}}}}
\end{picture}
}
\end{center}
\caption{The Swiss Flag}
\label{SwissFlag}
\end{figure}

In Figure~\ref{SwissFlag}, the left side is the Swiss Flag example. It
is again an example of cubical set appearing in the theory of
distributed databases as explained in \cite{FGR}.  The globular
$\omega$-category on the right side should be homotopic to the left
one.

Our claim is that the most interesting computer-scientific properties
of two concurrent machines are the same if the corresponding globular
$\omega$-categories are homotopy equivalent. In
Figure~\ref{SwissFlag}, the state $\gamma$ corresponds to a deadlock
of the corresponding concurrent machine. The deadlock appears also on
the right. Compare the number of possible execution paths on the left
and the only four execution paths on the right which are essentially
the same !  This means that an algorithm which could be able to take
in account this notion of $\omega$-category up to homotopy would be
more efficient that any other algorithm.

Instead of dealing directly with $\omega$-category up to homotopy, a
more fruitful approach consists of building some functors from
$\omega$-categories to, for example, abelian groups, invariant up to
homotopy. These functors contain, at least theoretically, a relevant
geometric information because of their invariance up to homotopy in
the above sense. A usual way to construct such invariants consists of
constructing functors from the category of $\omega$-categories to the
category $Comp(Ab)$ of chain complexes of groups and to consider the
associated homology groups.

\bd\label{chain} A \textit{chain complex of groups} is a family of
abelian groups $(C_n)_{n\geqslant 0}$ together with a family of linear
maps $\de_n:C_{n+1}\longrightarrow C_n$ such that $\de_n\circ
\de_{n+1}=0$ for any $n\geqslant 0$. \ed

Since the image of $\de_{n+1}$ is included in the kernel of $\de_n$,
the quotient group $$H_n(C_*,\de_*)=Ker(\de_n)/Im(\de_{n+1})$$ is well 
defined. It is called the $n$-th homology group of the group complex
$(C_n)_{n\geqslant 0}$. The map $H_n$ yields a functor from $Comp(Ab)$
to the category $Ab$ of abelian groups.  See for example \cite{Rotman}
or \cite{Weibel} for an introduction to the theory of these
mathematical objects.

The first homology theory will be the globular homology $H_*^{gl}(\C)$
(Definition~\ref{def_globulaire}). An example of globular cycle of
dimension $1$ is $\gamma_1-\gamma_4$ of Figure~\ref{trou}. We call it
an oriented $1$-dimensional loop. An example of globular cycle of
dimension $2$ is $A-B$ of Figure~\ref{3hole}.

The two other homology theories will be called the negative and
positive corner homologies $H_*^\pm (\C)$
(Definition~\ref{def_orientee}). The cycles of the negative one
correspond to the branching areas of execution paths (or negative
corner) and the positive one to the merging areas of execution paths
(or positive corner). In the case of $\gamma_1-\gamma_4$ of
Figure~\ref{trou}, there is one branching area on the left and one
merging area on the right. Idem for Figure~\ref{3hole}.

\begin{figure}
\begin{center}
\setlength{\unitlength}{0.0005in}
\begingroup\makeatletter\ifx\SetFigFont\undefined%
\gdef\SetFigFont#1#2#3#4#5{%
  \reset@font\fontsize{#1}{#2pt}%
  \fontfamily{#3}\fontseries{#4}\fontshape{#5}%
  \selectfont}%
\fi\endgroup%
{\renewcommand{\dashlinestretch}{30}
\begin{picture}(7532,3061)(0,-10)
\path(2100,3023)(5100,3023)(5100,23)
        (2100,23)(2100,3023)
\path(2700,1223)(3300,1223)(3300,623)
        (3900,623)(3900,1223)(4500,1223)
        (4500,1823)(3900,1823)(3900,2423)
        (3300,2423)(3300,1823)(2700,1823)
        (2700,1223)(2700,1223)
\path(2700,1223)(2700,1823)
\path(2760.000,1583.000)(2700.000,1823.000)(2640.000,1583.000)
\path(2700,1223)(3300,1223)
\path(3060.000,1163.000)(3300.000,1223.000)(3060.000,1283.000)
\path(4500,1223)(4500,1823)
\path(4560.000,1583.000)(4500.000,1823.000)(4440.000,1583.000)
\path(3900,1823)(4500,1823)
\path(4260.000,1763.000)(4500.000,1823.000)(4260.000,1883.000)
\path(2100,23)  (2103.675,87.925)
        (2107.461,151.696)
        (2111.360,214.326)
        (2115.375,275.829)
        (2119.511,336.218)
        (2123.771,395.508)
        (2128.158,453.713)
        (2132.675,510.845)
        (2137.327,566.919)
        (2142.115,621.948)
        (2147.045,675.947)
        (2152.119,728.928)
        (2157.341,780.906)
        (2162.714,831.895)
        (2168.242,881.907)
        (2173.928,930.957)
        (2179.775,979.059)
        (2185.787,1026.226)
        (2191.967,1072.472)
        (2198.319,1117.811)
        (2204.847,1162.257)
        (2211.553,1205.822)
        (2218.441,1248.522)
        (2225.515,1290.369)
        (2232.777,1331.377)
        (2240.233,1371.561)
        (2247.883,1410.934)
        (2255.733,1449.509)
        (2272.045,1524.322)
        (2289.195,1596.110)
        (2307.211,1664.983)
        (2326.119,1731.052)
        (2345.948,1794.425)
        (2366.726,1855.212)
        (2388.479,1913.525)
        (2411.234,1969.471)
        (2435.021,2023.162)
        (2459.865,2074.707)
        (2485.795,2124.217)
        (2512.837,2171.799)
        (2541.020,2217.566)
        (2570.370,2261.626)
        (2600.916,2304.090)
        (2632.685,2345.067)
        (2665.704,2384.667)
        (2700.000,2423.000)

\path(2700,2423)        (2739.210,2459.673)
        (2784.785,2492.695)
        (2835.847,2522.335)
        (2891.520,2548.864)
        (2950.926,2572.553)
        (3013.189,2593.671)
        (3077.432,2612.490)
        (3142.778,2629.280)
        (3208.349,2644.311)
        (3273.268,2657.854)
        (3336.660,2670.179)
        (3397.646,2681.556)
        (3455.350,2692.256)
        (3508.895,2702.550)
        (3557.404,2712.708)
        (3600.000,2723.000)

\path(3600,2723)        (3649.096,2735.745)
        (3701.582,2748.778)
        (3758.008,2762.208)
        (3818.923,2776.147)
        (3884.877,2790.703)
        (3956.418,2805.986)
        (3994.455,2813.935)
        (4034.096,2822.107)
        (4075.408,2830.516)
        (4118.460,2839.175)
        (4163.321,2848.098)
        (4210.060,2857.300)
        (4258.745,2866.793)
        (4309.445,2876.592)
        (4362.228,2886.710)
        (4417.164,2897.160)
        (4474.321,2907.958)
        (4533.767,2919.116)
        (4595.571,2930.648)
        (4659.803,2942.567)
        (4726.529,2954.889)
        (4795.820,2967.625)
        (4867.744,2980.791)
        (4942.370,2994.400)
        (4980.717,3001.374)
        (5019.765,3008.465)
        (5059.524,3015.673)
        (5100.000,3023.000)

\path(4761.770,2870.355)(5100.000,3023.000)(4729.729,3047.480)
\path(2100,23)  (2164.926,26.675)
        (2228.697,30.460)
        (2291.327,34.358)
        (2352.830,38.374)
        (2413.220,42.509)
        (2472.511,46.768)
        (2530.716,51.155)
        (2587.848,55.672)
        (2643.923,60.323)
        (2698.952,65.112)
        (2752.951,70.041)
        (2805.933,75.115)
        (2857.911,80.336)
        (2908.900,85.709)
        (2958.913,91.236)
        (3007.963,96.922)
        (3056.065,102.769)
        (3103.232,108.781)
        (3149.479,114.961)
        (3194.818,121.313)
        (3239.263,127.840)
        (3282.829,134.546)
        (3325.529,141.434)
        (3367.376,148.508)
        (3408.384,155.770)
        (3448.568,163.225)
        (3487.941,170.876)
        (3526.516,178.726)
        (3601.329,195.038)
        (3673.117,212.188)
        (3741.991,230.203)
        (3808.059,249.112)
        (3871.432,268.941)
        (3932.220,289.719)
        (3990.532,311.472)
        (4046.478,334.228)
        (4100.168,358.015)
        (4151.713,382.859)
        (4201.222,408.789)
        (4248.804,435.833)
        (4294.570,464.016)
        (4338.630,493.367)
        (4381.092,523.914)
        (4422.069,555.683)
        (4461.668,588.703)
        (4500.000,623.000)

\path(4500,623) (4536.673,662.210)
        (4569.695,707.785)
        (4599.335,758.847)
        (4625.864,814.520)
        (4649.553,873.926)
        (4670.671,936.189)
        (4689.490,1000.432)
        (4706.280,1065.777)
        (4721.311,1131.349)
        (4734.854,1196.268)
        (4747.179,1259.660)
        (4758.556,1320.646)
        (4769.256,1378.350)
        (4779.550,1431.895)
        (4789.708,1480.404)
        (4800.000,1523.000)

\path(4800,1523)        (4812.745,1572.096)
        (4825.778,1624.582)
        (4839.208,1681.008)
        (4853.147,1741.923)
        (4867.703,1807.877)
        (4882.986,1879.418)
        (4890.935,1917.455)
        (4899.107,1957.096)
        (4907.516,1998.408)
        (4916.175,2041.460)
        (4925.098,2086.321)
        (4934.300,2133.060)
        (4943.793,2181.745)
        (4953.592,2232.445)
        (4963.710,2285.228)
        (4974.160,2340.164)
        (4984.958,2397.321)
        (4996.116,2456.767)
        (5007.648,2518.571)
        (5019.567,2582.803)
        (5031.889,2649.529)
        (5044.625,2718.820)
        (5057.791,2790.744)
        (5071.400,2865.370)
        (5078.374,2903.717)
        (5085.465,2942.765)
        (5092.673,2982.524)
        (5100.000,3023.000)

\path(5124.480,2652.729)(5100.000,3023.000)(4947.355,2684.770)
\path(6000,2723)        (5959.548,2718.843)
        (5919.814,2714.704)
        (5880.788,2710.582)
        (5842.463,2706.476)
        (5767.879,2698.302)
        (5695.994,2690.168)
        (5626.739,2682.060)
        (5560.045,2673.965)
        (5495.843,2665.870)
        (5434.065,2657.759)
        (5374.642,2649.621)
        (5317.506,2641.440)
        (5262.587,2633.203)
        (5209.818,2624.896)
        (5159.129,2616.507)
        (5110.451,2608.020)
        (5063.717,2599.422)
        (5018.858,2590.700)
        (4975.804,2581.840)
        (4934.487,2572.827)
        (4894.839,2563.649)
        (4856.790,2554.291)
        (4785.218,2534.982)
        (4719.221,2514.791)
        (4658.250,2493.606)
        (4601.756,2471.320)
        (4549.189,2447.821)
        (4500.000,2423.000)

\path(4500,2423)        (4445.093,2386.531)
        (4381.478,2329.062)
        (4344.207,2290.256)
        (4302.124,2243.562)
        (4254.348,2188.104)
        (4200.000,2123.000)

\path(4253.221,2234.658)(4200.000,2123.000)(4299.505,2196.477)
\path(1500,398) (1560.238,405.791)
        (1618.343,413.430)
        (1674.366,420.927)
        (1728.359,428.293)
        (1780.374,435.537)
        (1830.461,442.671)
        (1878.673,449.705)
        (1925.060,456.648)
        (1969.675,463.511)
        (2012.568,470.305)
        (2053.792,477.039)
        (2093.398,483.724)
        (2167.961,496.989)
        (2236.669,510.181)
        (2299.934,523.384)
        (2358.168,536.678)
        (2411.784,550.148)
        (2461.193,563.874)
        (2506.807,577.940)
        (2549.038,592.428)
        (2588.299,607.421)
        (2625.000,623.000)

\path(2625,623) (2675.104,650.018)
        (2737.117,693.009)
        (2774.788,722.142)
        (2818.072,757.245)
        (2867.850,798.978)
        (2925.000,848.000)

\path(2853.878,746.799)(2925.000,848.000)(2814.621,792.174)
\put(6225,2723){\makebox(0,0)[lb]{\smash{{{\SetFigFont{11}{14.4}{\rmdefault}{\mddefault}{\updefault}unreachable area}}}}}
\put(0,98){\makebox(0,0)[lb]{\smash{{{\SetFigFont{11}{14.4}{\rmdefault}{\mddefault}{\updefault}unsafe area}}}}}
\end{picture}
}
\end{center}
\caption{Unsafe area and unreachable area in a concurrent machine with semaphores}
\label{forbiddenarea}
\end{figure}
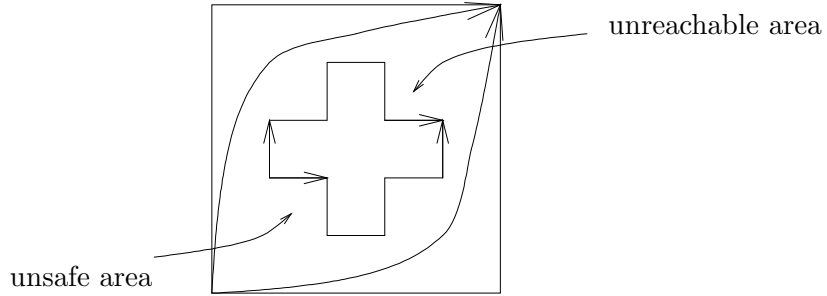

The idea afterwards is to associate to any oriented loop of any
dimension its corresponding negative or positive corners. This is the
underlying geometric meaning of the morphisms $h_*^\pm $ from
$H_*^{gl}(\C)$ to $H_*^\pm (\C)$ (Proposition~\ref{hmoins} and
\ref{hplus}). We can immediately see an application of these maps.
Looking back to the Swiss Flag example of Figure~\ref{SwissFlag}, it
is clear that the cokernel of $h_1^-$ is not empty, because of the
deadlock and the forbidden area. A negative corner which yields a non
trivial element in this cokernel is drawn in
Figure~\ref{forbiddenarea}. In the same way, the cokernel of $h_1^+$
in the Swiss Flag example is still not empty, because of the
unreachable state and the unreachable area.  A positive corner which
yields a non zero element of this cokernel is represented in
Figure~\ref{forbiddenarea}.

These geometric remarks ensure that the homology groups that we are
going to construct contain relevant information about the geometry of
concurrency.

\section{Globular homology of $\omega$-category}\label{glob_homo}

\subsection{Definition}

The starting point is an $\omega$-category $\C$.

\bd\label{def_globulaire} Let $(C_*^{gl}(\C),\de^{gl})$ be the chain
complex defined as follows : $C_0^{gl}(\C)=\Z\C_0\oplus \Z\C_0$, and
for $n\geqslant 1$, $C_n^{gl}(\C)=\Z\C_n$, $\de^{gl}(x)=(s_0x,t_0x)$
if $x\in\Z\C_1$ and for $n\geqslant 1$, $x\in\Z\C_{n+1}$ implies
$\de^{gl}(x)=s_nx-t_nx$. This complex is called the \textit{globular
  complex} of $\C$ and its corresponding homology the globular
homology.  \ed

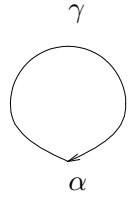
\begin{figure}
\begin{center}
\setlength{\unitlength}{0.0005in}
\begingroup\makeatletter\ifx\SetFigFont\undefined%
\gdef\SetFigFont#1#2#3#4#5{%
  \reset@font\fontsize{#1}{#2pt}%
  \fontfamily{#3}\fontseries{#4}\fontshape{#5}%
  \selectfont}%
\fi\endgroup%
{\renewcommand{\dashlinestretch}{30}
\begin{picture}(1216,1941)(0,-10)
\path(608,318)  (566.687,336.835)
        (527.104,355.390)
        (453.020,391.770)
        (385.527,427.358)
        (324.406,462.375)
        (269.438,497.041)
        (220.402,531.574)
        (177.080,566.196)
        (139.250,601.125)
        (79.191,672.785)
        (56.523,709.956)
        (38.469,748.312)
        (24.809,788.076)
        (15.324,829.465)
        (9.794,872.700)
        (8.000,918.000)

\path(8,918)    (11.113,979.121)
        (20.247,1038.527)
        (35.093,1095.913)
        (55.345,1150.969)
        (80.695,1203.388)
        (110.835,1252.863)
        (145.457,1299.086)
        (184.254,1341.750)
        (226.918,1380.546)
        (273.141,1415.168)
        (322.616,1445.307)
        (375.035,1470.656)
        (430.091,1490.908)
        (487.476,1505.754)
        (546.881,1514.887)
        (608.000,1518.000)

\path(608,1518) (669.119,1514.887)
        (728.524,1505.753)
        (785.909,1490.907)
        (840.965,1470.655)
        (893.384,1445.305)
        (942.859,1415.165)
        (989.082,1380.543)
        (1031.746,1341.746)
        (1070.543,1299.082)
        (1105.165,1252.859)
        (1135.305,1203.384)
        (1160.655,1150.965)
        (1180.907,1095.909)
        (1195.753,1038.524)
        (1204.887,979.119)
        (1208.000,918.000)

\path(1208,918) (1206.206,872.700)
        (1200.676,829.465)
        (1191.191,788.076)
        (1177.531,748.312)
        (1159.477,709.956)
        (1136.809,672.785)
        (1076.750,601.125)
        (1038.920,566.196)
        (995.598,531.574)
        (946.562,497.041)
        (891.594,462.375)
        (830.473,427.358)
        (762.980,391.770)
        (688.896,355.390)
        (649.313,336.835)
        (608.000,318.000)

\path(705.120,394.601)(608.000,318.000)(729.742,339.886)
\put(608,18){\makebox(0,0)[lb]{\smash{{{\SetFigFont{11}{14.4}{\rmdefault}{\mddefault}{\updefault}$\alpha$}}}}}
\put(608,1818){\makebox(0,0)[lb]{\smash{{{\SetFigFont{11}{14.4}{\rmdefault}{\mddefault}{\updefault}$\gamma$}}}}}
\end{picture}
}
\end{center}
\caption{A loop which does not give rise to a globular cycle}
\label{dim1}
\end{figure}

\begin{figure}
\begin{center}
\setlength{\unitlength}{0.0005in}
\begingroup\makeatletter\ifx\SetFigFont\undefined%
\gdef\SetFigFont#1#2#3#4#5{%
  \reset@font\fontsize{#1}{#2pt}%
  \fontfamily{#3}\fontseries{#4}\fontshape{#5}%
  \selectfont}%
\fi\endgroup%
{\renewcommand{\dashlinestretch}{30}
\begin{picture}(2782,2366)(0,-10)
\path(927,27)(2727,627)
\path(2622.645,560.592)(2727.000,627.000)(2603.671,617.513)
\path(927,27)   (867.257,57.638)
        (809.966,87.701)
        (755.086,117.232)
        (702.575,146.271)
        (652.392,174.859)
        (604.497,203.039)
        (558.847,230.850)
        (515.402,258.334)
        (474.120,285.532)
        (434.961,312.486)
        (362.844,365.825)
        (298.722,418.681)
        (242.265,471.383)
        (193.143,524.260)
        (151.028,577.642)
        (115.588,631.859)
        (86.496,687.240)
        (63.420,744.115)
        (46.032,802.814)
        (34.002,863.666)
        (27.000,927.000)

\path(27,927)   (25.103,977.593)
        (26.544,1027.705)
        (31.183,1077.299)
        (38.878,1126.339)
        (49.492,1174.788)
        (62.882,1222.611)
        (78.910,1269.771)
        (97.435,1316.231)
        (118.317,1361.955)
        (141.416,1406.908)
        (166.592,1451.051)
        (193.705,1494.351)
        (222.614,1536.769)
        (253.181,1578.270)
        (285.264,1618.817)
        (318.724,1658.374)
        (353.420,1696.904)
        (389.213,1734.372)
        (463.527,1805.975)
        (501.769,1840.037)
        (540.547,1872.891)
        (579.721,1904.501)
        (619.151,1934.830)
        (658.697,1963.842)
        (698.219,1991.501)
        (737.576,2017.770)
        (776.630,2042.613)
        (815.239,2065.994)
        (853.264,2087.876)
        (927.000,2127.000)

\path(927,2127) (963.887,2144.824)
        (1002.559,2162.285)
        (1042.902,2179.306)
        (1084.801,2195.807)
        (1128.143,2211.710)
        (1172.813,2226.935)
        (1218.698,2241.404)
        (1265.684,2255.039)
        (1313.655,2267.760)
        (1362.500,2279.489)
        (1412.103,2290.146)
        (1462.350,2299.654)
        (1513.127,2307.933)
        (1564.321,2314.905)
        (1615.817,2320.490)
        (1667.501,2324.610)
        (1719.260,2327.186)
        (1770.978,2328.140)
        (1822.543,2327.392)
        (1873.840,2324.865)
        (1924.755,2320.478)
        (1975.174,2314.153)
        (2024.983,2305.812)
        (2074.068,2295.376)
        (2122.315,2282.765)
        (2169.610,2267.902)
        (2215.839,2250.707)
        (2260.888,2231.102)
        (2304.642,2209.007)
        (2346.988,2184.344)
        (2387.812,2157.035)
        (2427.000,2127.000)

\path(2427,2127)        (2482.869,2077.449)
        (2533.277,2024.568)
        (2578.332,1967.808)
        (2618.145,1906.620)
        (2652.826,1840.454)
        (2682.484,1768.761)
        (2695.464,1730.671)
        (2707.230,1690.992)
        (2717.795,1649.657)
        (2727.172,1606.598)
        (2735.377,1561.744)
        (2742.422,1515.028)
        (2748.322,1466.381)
        (2753.089,1415.734)
        (2756.739,1363.019)
        (2759.283,1308.167)
        (2760.737,1251.109)
        (2761.114,1191.776)
        (2760.427,1130.101)
        (2758.691,1066.014)
        (2755.919,999.446)
        (2752.124,930.330)
        (2747.322,858.596)
        (2741.524,784.175)
        (2738.257,745.936)
        (2734.746,706.999)
        (2730.993,667.357)
        (2727.000,627.000)

\path(1377,1827)        (1435.107,1818.765)
        (1488.944,1810.137)
        (1538.703,1801.035)
        (1584.577,1791.375)
        (1626.758,1781.076)
        (1665.439,1770.055)
        (1733.066,1745.516)
        (1788.998,1717.101)
        (1834.772,1684.149)
        (1871.927,1646.002)
        (1902.000,1602.000)

\path(1902,1602)        (1926.137,1551.002)
        (1941.138,1496.662)
        (1946.785,1437.223)
        (1942.856,1370.925)
        (1929.133,1296.012)
        (1918.530,1254.775)
        (1905.396,1210.725)
        (1889.703,1163.642)
        (1871.425,1113.307)
        (1850.533,1059.500)
        (1827.000,1002.000)

\path(702,1227) (711.420,1186.627)
        (720.579,1149.152)
        (738.327,1082.351)
        (755.685,1025.498)
        (773.092,977.494)
        (809.812,903.640)
        (852.000,852.000)

\path(852,852)  (915.453,806.099)
        (956.145,786.347)
        (1004.572,768.142)
        (1062.053,751.045)
        (1129.906,734.615)
        (1168.134,726.513)
        (1209.449,718.413)
        (1254.016,710.260)
        (1302.000,702.000)

\put(1227,1227){\makebox(0,0)[lb]{\smash{{{\SetFigFont{11}{14.4}{\rmdefault}{\mddefault}{\updefault}$A$}}}}}
\put(1827,27){\makebox(0,0)[lb]{\smash{{{\SetFigFont{11}{14.4}{\rmdefault}{\mddefault}{\updefault}$\gamma$}}}}}
\end{picture}
}
\end{center}
\caption{Example of globular cycle in higher dimension}
\label{dimp}
\end{figure}
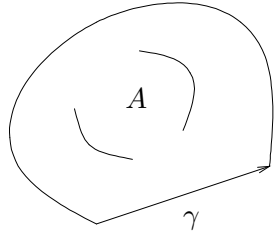

There is a difference between the $1$-dimensional case and the other
cases. A loop as in Figure~\ref{dim1} where $\gamma$ is a $1$-morphism
such that $s_0\gamma=t_0\gamma=\alpha$ does not yield a globular
$1$-cycle. However a loop as in Figure~\ref{dimp} where $A$ is a
$n$-morphism with $n\geqslant 2$ and such that
$s_{n-1}A=t_{n-1}A=\gamma$ yields a globular $n$-cycle.

\subsection{Functorial property of the globular homology}

Now an important technical definition for the sequel.

\bd\label{noncontractant} Let $f$ be an $\omega$-functor from $\C$ to
$\D$.  The morphism $f$ is \textit{non $1$-contracting} if for any
$1$-dimensional $x\in \C$, the morphism $f(x)$ is a $1$-dimensional
morphism of $\D$. \ed

The category of $\omega$-categories with the non $1$-contracting
$\omega$-functors is denoted by $\omega Cat_1$. The category of
cubical sets equipped with the non $1$-contracting morphisms is
denoted by $Sets_1^{\square^{op}}$.

If $f$ is a non $1$-contracting $\omega$-functor from $\C$ to $\D$, then
for any morphism $x\in \C$ of dimension greater than $1$, $f(x)$ is of
dimension greater than one as well. This is due to the equality $f(s_1
x)=s_1 f(x)$.

Let $f$ be an $\omega$-functor from $\C$ to $\D$. Then $f$ induces for
all $n\geq 0$ a linear morphism $f_*$ from $\Z\C_n$ to $\Z\D_n$ by
setting $f_*(x)= f(x) \hbox{ mod }\D_{\leqslant n-1}$ : this notation
meaning that $f_*(x)= f(x)$ if $f(x)$ is $n$-dimensional and
$f_*(x)=0$ otherwise. For $n\geq 2$,
$f_*(s_{n-1}-t_{n-1})(x)=(s_{n-1}-t_{n-1})f_*(x)$, therefore
$f_*\de^{gl}(x)=\de^{gl}f_*(x)$. The latter equality is not anymore
true if $x$ is $1$-dimensional because an $\omega$-functor can
contract $1$-morphisms and because $\de^{gl}(x)=(s_0(x),t_0(x))$. So
the globular homology does not yield a functor from $\omega Cat$ to
$Ab$ but only a functor from $\omega Cat_1$ to $Ab$.

\subsection{Homological property}

Now we give a homological property of the globular complex to justify
this construction. The starting point is the small category $Glob$
defined as follows : the objects are all natural numbers and the
arrows are generated by $s$ and $t$ in $Glob(m,m-1)$ for any $m>0$ and
quotiented by the relations $ss=st$, $ts=tt$. We can depict $Glob$
like this :

\[\xymatrix@1{\ar@<-1ex>@{->}[r]_{s}\ar@<1ex>@{->}[r]^{t}& 3 
\ar@<-1ex>@{->}[r]_{s}\ar@<1ex>@{->}[r]^{t}& 2 
\ar@<-1ex>@{->}[r]_{s}\ar@<1ex>@{->}[r]^{t}& 1 
\ar@<-1ex>@{->}[r]_{s}\ar@<1ex>@{->}[r]^{t}& 0}\]

\bd \label{globset} A \textit{globular set} is a covariant functor
from $Glob$ to $Sets$.  The corresponding category is denoted by
$[Glob,Sets]$. A \textit{globular group} is a covariant functor from
$Glob$ to the category $Ab$ of abelian groups. The corresponding
category is denoted by $[Glob,Ab]$.  \ed

The notion of globular set already appears in many works and is
certainly not new \cite{petittoposglob}\cite{Penon_weakcat}\cite{Batanin_1}.

If $\C$ is an $\omega$-category, we denote by $\C_{n}$ the set of
$n$-dimensional morphisms of $\C$ with $n\geqslant 0$.

\bd\label{transglob} Let $Gl$ be the map from $\omega$-categories to
globular groups defined as follows. If $n\geqslant 0$, set
$Gl(\C)_n=\Z\C_{n}$. For any $s,t\in Glob(n+1,n)$, set
$Gl(s)(x)=s_n(x)$ and $Gl(t)(x)=t_n(x)$.  \ed

Unfortunately, $Gl(-)$ is not a functor because an $\omega$-functor
might be $n$-contracting for $n\geqslant 2$. That is, suppose that $f$
is an $\omega$-functor from $\C$ to $\D$ such that for a $2$-morphism
$x$ of $\C$, $f(x)$ is $1$-dimensional. Then $Gl(f)(x)=0\in Gl(\D)_2$
and $s_1 Gl(f)(x)=0$ but $Gl(f)(s_1(x))=f(s_1(x))\in Gl(\D)_1$ and
$f(s_1(x))\neq 0$.

If $M$ is a globular group, let $H(M)$ be the cokernel of the additive
map from $M_1$ to $M_0\oplus M_0$ which maps $x$ to $(s(x),t(x))$. The
map $H$ induces a right exact additive functor from $[Glob,Ab]$ to
$Ab$. Since $[Glob,Ab]$ has enough projectives, we can deal with the
left derived functors $L_n(H)$ of $H$ (see \cite{Weibel} or
\cite{Rotman} for the definition of projective object and right exact
functor).

\bp For any $\omega$-category $\C$, we have 
$H_*^{gl}(\C)\iso L_*(H)(Gl(\C))$.
\ep

Before proving this theorem, we need to recall standard facts about
category of diagrams. We are going to solve exercice (2.3.13) of
\cite{Weibel} because in the sequel we need a precise description of a
family of projective globular groups which allows to resolve any
globular group.

Let $ev_k$ be the functor from $[Glob,Ab]$ to $Ab$ such that
$ev_k(M)=M(k)$, $k$ being a natural number. This functor is exact and
by the special adjoint functor theorem has a left adjoint denoted
by $k_!$. We need to explicit $k_!$ for the sequel.

\bp If $M$ is an abelian group, set 
$$k_!(M)(l)=\bigoplus_{h\in Glob(k,l)}M_h$$
where $M_h$ is a copy of
$M$. If $x\in M$, let $h.x$ be the corresponding element of
$k_!(M)(l)$. If $f:l\longrightarrow l'$ is an arrow of $I$, then we
set $k_!(M)(f)(h.x)=(fh).x$.  Then $k_!$ is a globular group and this
is the left adjoint of $ev_k$.  \ep

\bpf Let ${N}$ be a globular group. We introduce the map
$$F:\xymatrix@1{{[Glob,Ab]\left(k_!(M),{N}\right)}\fr{}&
  {Ab\left(M,{N}(k)\right)}}$$ defined by $F(u)(x)= u(Id_k.x)$ and the
map $$\xymatrix@1{ G:Ab(M,{N}(k))\fr{} &
  [Glob,Ab]\left(k_!(M),{N}\right) }$$ defined by
$G(v)(f.x)={N}(f)(v(x))$. The arrow $G(v)$ is certainly a morphism of
globular groups from $k_!(M)$ to ${N}$. In fact, if
$\xymatrix@1{l\fr{f}& l'}$ is an arrow of $Glob$ and if $h.x$ is an
element of $k_!(M)(l)$, then
$G(v)(fh.x)={N}(fh)(v(x))={N}(f)(G(v)(h.x))$. Therefore the diagram
$$\xymatrix{ {k_!(M)(l)} \fr{G(v)}\fd{k_!(M)(f)} &
{{N}(l)}\fd{{N}(f)}\\ {k_!(M)(l')}\fr{G(v)}& {{N}(l')}} $$ commutes.
Now we have to verify that $F$ and $G$ are inverse of each other.
Indeed, $$F(G(v))(x)=G(v)(Id_k.x)={N}(Id)(v(x))$$ therefore
$F(G(v))=v$. Conversely,
$$G(F(u))(h.x)={N}(h)(F(u)(x))={N}(h)(u(Id_k.x))=u(k_!(M)(h)(Id_k.x))=u(h.x)$$
therefore $G(F(u))=u$.
\epf

Now we set $$\mathcal{F}=\left\{\bigoplus_{k\in\N}k_!(L)/L \hbox{
free module}\right\}$$ And we can state the proposition

\bp\label{proj_canonique} All elements of $\mathcal{F}$ are projective
globular groups. Any globular group can be resolved by elements of
$\mathcal{F}$.  \ep

\bpf Let ${X}$ be a globular group. For any $k\in\N$, let $L_k$ be a
free abelian group and $L_k \twoheadrightarrow {X}(k)$ an epimorphism
of abelian groups. Then the epimorphisms $\xymatrix@1{L_k\fr{}&
  {X}(k)}$ for all $k$ induce a natural transformation
$\xymatrix@1{\bigoplus_{k\in\N} k_!(L_k)\fr{} & X}$ which is certainly
itself an epimorphism. Left adjoint functors and coproduct preserve
projective objects \cite{borceux1}.  Hence the conclusion.  \epf

We are in position to prove the proposition :

\bp For any $\omega$-category $\C$, the equality 
$H_*^{gl}(\C)\iso L_*(H)(Gl(\C))$ holds.
\ep

\bpf If $M$ is a globular group, we introduce the complex of abelian
groups $(C^{gl}_*(M),\de^{gl})$ defined as follows :
$C^{gl}_0(M)=M_0\oplus M_0$ and for $n\geqslant 1$, $C^{gl}_n(M)=M_n$,
with the differential map $\de^{gl}(x)=(s(x),t(x))$ if $x\in M_1$ and
$\de^{gl}(x)=s(x)-t(x)$ if $x\in M_n$ with $n\geqslant 2$. We have
$H_0(C^{gl}_*(M),\de)=H(M)$. Let $k$ be a natural number and let $L$
be a free abelian group. If $p>0$, let us prove that
$H_p(C^{gl}_*(k_!(L)))=0$. Let $X=x_p $ be a cycle of $C_p(k_!(L))$.
By construction, for all $p>k$, one has $k_!(L)(p)=0$.  Therefore if
$p> k$, then $X=0$ hence $H_p(C^{gl}_*(k_!(L)))=0$ whenever $p>k$. Now let
us see the case $p\leqslant k$. We have
$0=\de^{gl}(X)=s_{p-1}(x_p)-t_{p-1}(x_p)$.  Then there exists
$x_{p}^s$ and $x_{p}^t$ such that
$x_p=s^{k-p}.x^s_{p}+t^{k-p}.x^t_{p}$. The equality
$s_{p-1}(x_p)=t_{p-1}(x_p)$ implies
$s^{k-p+1}.x^s_{p}+s^{k-p+1}.x^t_{p}=t^{k-p+1}.x^s_{p}+t^{k-p+1}.x^t_{p}$.
Therefore $x_p=0$.

Now we deduce from Proposition~\ref{proj_canonique} that for any
projective globular group $P$, $$H_p(C^{gl}_*(P),\de)=0$$ for all
$p>0$. It is thus easy to check that for any natural number $p$, the
equality $H_p(C^{gl}_*(M),\de)=L_p(H)(M)$ holds. Indeed, the case
$p=0$ is trivial and the commutative diagram
$$\xymatrix{H_{p+1}(C^{gl}_*(P),\de)\fr{}\fd{}&H_{p+1}(C^{gl}_*(M),\de)\fr{}\fd{}&
H_{p}(C^{gl}_*(K),\de)\fd{}\fr{} & H_{p}(C^{gl}_*(P),\de)\fd{}\\
L_{p+1}(H)(P) \fr{} &L_{p+1}(H)(M)
\fr{}&L_{p}(H)(K) \fr{}& L_{p}(H)(P)\\}$$
with $P$ projective allows to make the induction on $p$.

\epf

\section{Positive and negative corner homology  of an $\omega$-category}\label{corner}

\subsection{The pasting scheme $\Lambda^n$ and the $\omega$-category
  $I^n$}\label{construction-of-In}

We need to describe precisely the $\omega$-category associated to the
$n$-cube.

\bd\label{recollement} A \textit{pasting scheme}  is a triple $(A,E,B)$ where 
$A$ is a $\N$-graded set, and $E^i_j$ and $B^i_j$ two binary relations
over $A_i\p A_j$ with $j\leqslant i$ satisfying
\begin{enumerate}
\item[(1)] the set $E^i_i$ is the diagonal of  $A_i$
\item[(2)] for  $k>0$, and for any  $x\in A_k$, there exists  $y\in A_{k-1}$
  with $xE^k_{k-1} y$
\item[(3)] for  $k<n$, $wE^n_k x$ if and only if there exists  $u$ and $v$
 such that $wE^n_{n-1}uE^{n-1}_k x$ and $wE^n_{n-1}vB^{n-1}_k x$
\item[(4)] if $wE^n_{
n-1}zE^{n-1}_k x$, then either  $wE^n_k x$ or there exists
  $v$ such that  $wB^n_{n-1}vE^{n-1}_k x$
\end{enumerate}
and such that $(A,B,E)$ satisfies the same properties. If $x\in A_i$,
we set $dim(x)=i$. 
\ed

If $x$ is an element of the pasting scheme $(A,E,B)$, we denote by
$R(x)$ the smallest pasting scheme of $(A,E,B)$ containing $x$.

Intuitively, a pasting scheme is a pasting of faces of several
dimensions together \cite{CPS}. Kapranov and Voedvosky have their own
formalization using some particular chain complexes of abelian groups
\cite{KapVod}. We can see Figure~\ref{composition-of-2trans} as a
pasting scheme $(S,E,B)$. It suffices to set
$S=\{\alpha,\beta,s_1A,t_1A,t_1B,A,B\}$ endowed with the binary
relations $B^2_1=\{(A,s_1 A),(B,s_1B)\}$, $E^2_1=\{(A,t_1
A),(B,t_1B)\}$, $B^2_0=E^2_0=\{\}$,
$B^1_0=\{(s_1A,\alpha),(t_1A,\alpha), (t_1B,\alpha)\}$,
$E^1_0=\{(s_1A,\beta),(t_1A,\beta), (t_1B,\beta)\}$.
Figure~\ref{morecomplicated} shows another more complicated example of
pasting scheme.

\begin{figure}
\[
\xymatrixrowsep{1.5pc}
\xymatrixcolsep{3pc}
\xymatrix{
&&\relax\rtwocell<0>^{w^{}\;\;}{\omit}
&\relax\ddtwocell<0>{\omit}
\drtwocell<0>^{\;\;x^{}}{<3>}
\ddrrtwocell<\omit>{<8>}\\
&&&&\relax\drtwocell<0>^{\;\;y^{}}{\omit}\\
\alpha\relax
\uurrcompositemap<2>_{u^{}}^{v^{}}{<.5>}
\drtwocell<0>_{f^{}\;}{\omit}
&&&\relax\urtwocell<0>{\omit}
&&\beta\\
&\relax \urrtwocell<0>{\omit}
\xcompositemap[-1,4]{}%
<-4.5>_{g^{}}^{h^{}}{\omit}\\
}\]
\caption{A  pasting scheme}\label{morecomplicated}\end{figure}

We only want here to recall the construction of the free
$\omega$-category $I^n$ generated by the faces of the $n$-cube.  For
more details see \cite{Crans_Tensor_product}, for an analogous
construction for simplices see \cite{oriental}, and for some explicit
calculations on $I^n$ see \cite{explicite}.

Set $\underline{n}=\{1,...,n\}$ and let $\Lambda^n$ be the set of maps
from $\underline{n}$ to $\{-,0,+\}$. We say that an element $x$ of
$\Lambda^n$ is of dimension $p$ if $x^{-1}(0)$ is a set of $p$
elements. We can identify the elements of $\Lambda^n$ with the words
of length $n$ in the alphabet $\{-,0,+\}$. The set $\Lambda^n$ is
supposed to be graded by the dimension of its elements. The set
$\Lambda^0$ is the set of maps from the empty set to $\{-,0,+\}$ and
therefore it is a singleton.

Let $y\in \Lambda^i$. Let $r_y$ be the map from $(\Lambda^n)_i$ to
$(\Lambda^n)_{dim(y)}$ defined as follows (with $x\in (\Lambda^n)_i$)
: for $k\in\underline{n}$, $x(k)\neq 0$ implies $r_y(x)(k)=x(k)$ and
if $x(k)$ is the $l$-th zero of the sequence $x(1),...,x(n)$, then
$r_y(x)(k)=y(\ell)$. If for any $\ell$ between $1$ and $i$,
$y(\ell)\neq 0$ implies $y(\ell)=(-)^\ell$, then we set
$b_y(x):=r_y(x)$. If for any $\ell$ between $1$ and $i$, $y(\ell)\neq
0$ implies $y(\ell)=(-)^{\ell+1}$, then we set $e_y(x):=r_y(x)$. We
thus introduce the following binary relations : the set $B^i_j$ of
pairs $(x,z)$ in $(\Lambda^n)_i\p (\Lambda^n)_j$ such that there
exists $y$ such that $z=b_y(x)$ and the set $E^i_j$ of pairs $(x,z)$
in $(\Lambda^n)_i\p (\Lambda^n)_j$ such that there exists $y$ such
that $z=e_y(x)$. Then $\Lambda^n$ is a pasting scheme. We have

\bth 
If $X\subset\Lambda^n$, let $R(X)$ be the sub-pasting scheme of
$(\Lambda^n,B^i_j,E^i_j)$ generated by $X$. There is one and only one
$\omega$-category $I^n$ such that
\begin{enumerate}
\item the underlying set of $I^n$ is included in the set of
%well-formed 
sub-pasting schemes of $(\Lambda^n,B^i_j,E^i_j)$ and it contains 
all pasting schemes like $R(\{x\})$ where $x$ runs over $\Lambda^n$
\item all elements of $I^n$ are a composition of $R(\{x\})$ where $x$ runs over $\Lambda^n$
\item for $x$ $p$-dimensional with $p\geqslant 1$, one has
\beas
&&s_{p-1}(R(\{x\}))=R\left(\{b_y(x),dim(y)=p-1\}\right)\\
&&t_{p-1}(R(\{x\}))=R\left(\{e_y(x),dim(y)=p-1\}\right)
\eeas
\item if $R(X)$ and $R(Y)$ are two elements of $I^n$ such that 
$t_p(R(X))=s_p(R(Y))$ for some $p$, then $R(X\cup Y)\in I^n$ and 
$R(X\cup Y)=R(X) *_p R(Y)$.
\end{enumerate}
\eth

The oriented $2$-cube is drawn in Figure~\ref{I2}.  With the rules
exposed in the above theorem, we can calculate $s_2R(00)$.  We have
actually $s_2R(00)=R(\{-0,0+\})$. But $t_0R(-0)=R(-+)=s_0R(0+)$.  Then
$s_2R(00)=R(\{-0\}\cup \{0+\})=R(-0)*_0 R(0+) $.

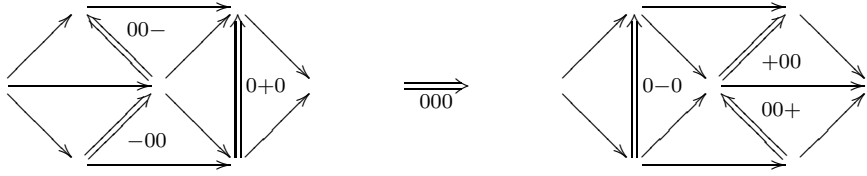
\begin{figure}
\[\xymatrix{
&\ar@{->}[rr]& &\ar@{->}[rd]&&&&
&\ar@{->}[dr]\ar@{->}[rr]&&\ar@{->}[dr]&\\
\ar@{->}[ru]\ar@{->}[rr]\ar@{->}[rd]&&\ar@{->}[ru]\ar@{->}[dr]\ff{lu}{00-}&&&\ff{r}{000}&&
\ar@{->}[ur]\ar@{->}[dr]&&\ff{ur}{+00}\ar@{->}[rr]&&\\
&\ar@{->}[rr]\ff{ru}{-00}&&\ar@{->}[ru]\ff{uu}{0+0}&&&&
&\ar@{->}[ru]\ar@{->}[rr]\ff{uu}{0-0}&&\ff{ul}{00+}\ar@{->}[ru]\\
}\]
\caption{The $\omega$-category $I^3$}
\label{I3}
\end{figure}

The $\omega$-category generated by a $3$-cube is drawn in
Figure~\ref{I3}. Let us give the example of the calculation of $s_2
R(000)$. We have
$$s_2R(000)=R(\{-00,0+0,00-\})=R(\{-00,0++\}\cup
\{-0-,0+0\}\cup \{00-,0++\})$$ 
since  $0++,-0-,0++ \in
R(\{-00,0+0,00-\})$. We verify easily that
$t_1R(\{-00,0++\})=s_1R(\{-0-,0+0\})$ and
$t_1R(\{-0-,0+0\})=s_1R(\{00-,0++\})$.  Therefore
$$s_2R(000)=R(\{-00,0++\}) *_1 R(\{-0-,0+0\}) *_1 R(\{00-,0++\}).$$ It
is then easy to verify that $R(\{-00,0++\})=R(-00)*_0 R(0++)$,
$R(\{-0-,0+0\})=R(-0-)*_0 R(0+0)$ and $R(\{00-,0++\})=R(00-)*_0
R(0++)$.

\begin{figure}
\begin{center}
\[\epsfig{file=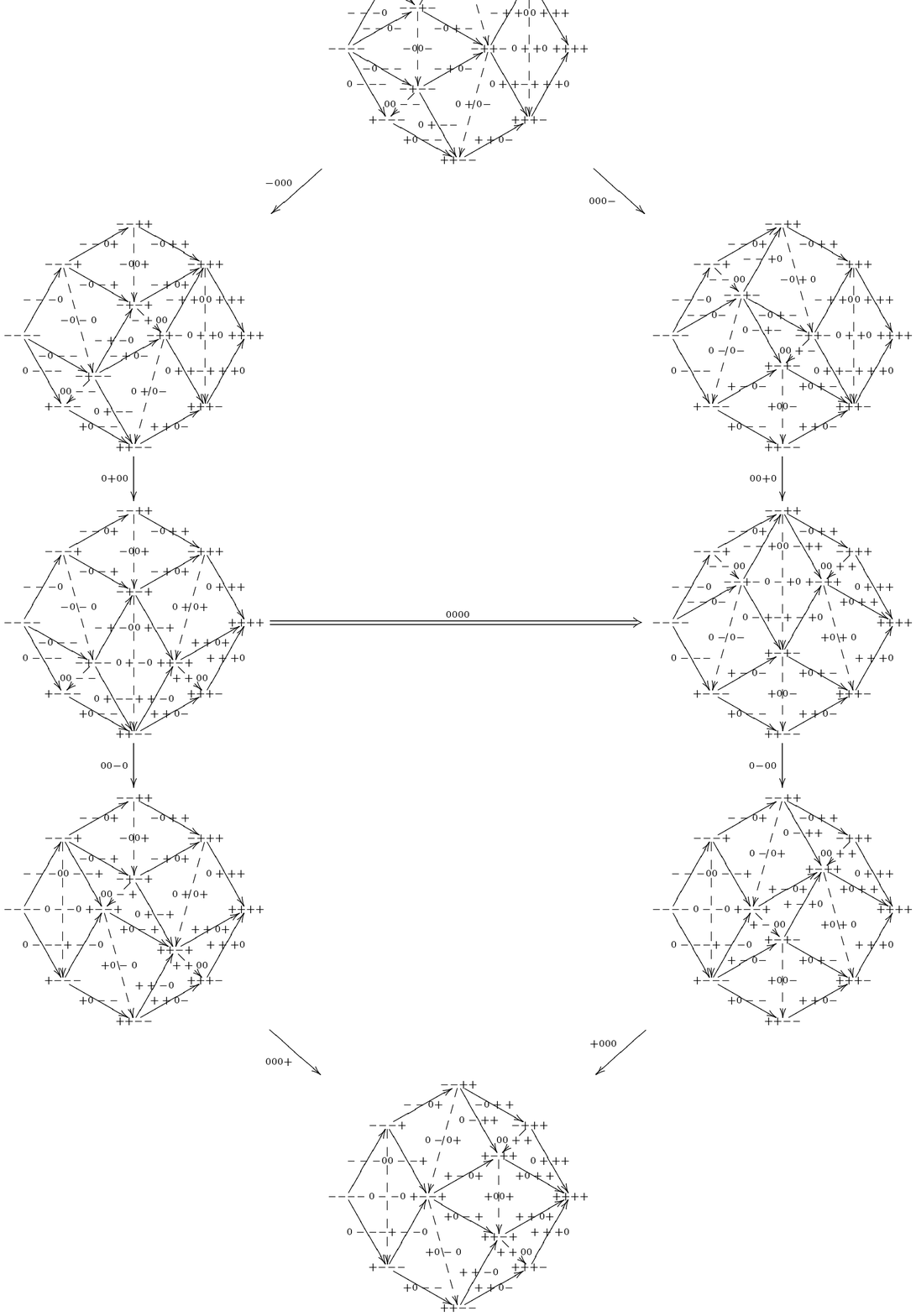,width=14.5cm,height=19cm}\]
\end{center}
\caption{The $\omega$-category $I^4$}
\label{I4}
\end{figure}

The oriented $4$-cube is represented in Figure~\ref{I4}.

%\[\epsfig{file=4cube.eps,width=14cm,height=19cm}\]

\subsection{The cubical singular nerve}

The map which sends every $\omega$-category $\C$ to
$\mathcal{N}^\square(\C)_*=\omega Cat(I^*,\C)$ induces a functor from
$\omega Cat$ to the category of cubical sets. If $x$ is an
element of $\omega Cat(I^n,\C)$, $\epsilon_i(x)$ is the
$\omega$-functor from $I^{n+1}$ to $\C$ defined by
$\epsilon_i(x)(k_1...k_{n+1})=x(k_1...\widehat{k_i}...k_{n+1})$ for
all $i$ between $1$ and $n+1$ and $\de_i^\alpha(x)$ is the
$\omega$-functor from $I^{n-1}$ to $\C$ defined by
$\de_i^\alpha(x)(k_1...k_{n-1})=x(k_1...k_{i-1}\alpha k_i...k_{n-1})$
for all $i$ between $1$ and $n$.

The arrow $\de_i^\alpha$ for a given $i$ such that $1\leqslant
i\leqslant n$ induces a natural transformation from $\omega
Cat(I^n,-)$ to $\omega Cat(I^{n-1},-)$ and therefore, by Yoneda,
corresponds to an $\omega$-functor $\delta^\alpha_i$ from $I^{n-1}$ to
$I^n$.  This functor is defined on the faces of $I^{n-1}$ by
$\delta^\alpha_i(k_1...k_{n-1})=R(k_1...[\alpha]_i...k_{n-1})$. The
notation $[...]_i$ means that the term inside the brackets are in the
$i$-th place.

\bd\label{cubical_singular} The cubical set $\omega Cat(I^*,\C)$
is called the \textit{cubical singular nerve} of the $\omega$-category
$\C$. \ed

This functor is a right adjoint. Its left adjoint is the functor $F$.

\subsection{Construction of the corner homologies}

The starting point is the cubical singular nerve $\omega Cat(I^*,\C)$
of $\C$ which contains all $n$-cubes included in $\C$. The main idea
to build the positive and negative corner homology of an
$\omega$-category $\C$ is to separate the two differential maps
$\de^-=\sum_i (-1)^{i+1}\de_i^-$ and $\de^+=\sum_i (-1)^{i+1}\de_i^+$
and to separately consider the chain complexes of groups $(\Z\omega
Cat(I^*,\C),\de^\pm )$ (a bit as in \cite{HDA} where the author
separates the horizontal and vertical differential maps of a
bicomplex). However the following proposition holds :

\bp Both chain complexes of groups $(\Z\omega Cat(I^*,\C),\de^-)$ and
$(\Z\omega Cat(I^*,\C),\linebreak[4]\de^+)$ are acyclic. \ep

\bpf It turns out that $\epsilon_1$ is a chain retraction of $(\Z\omega
Cat(I^*,\C),\de^-)$. If $x\in \Z\omega Cat(I^0,\C)$, then $\de^\alpha
\epsilon_1 x = \de^\alpha_1 \epsilon_1 x=x$. And for 
$x\in \Z\omega Cat(I^n,\C)$, with
$n\geqslant 1$, we have actually :

\beas
\de^\alpha \circ
\epsilon_1 (x)+\epsilon_1\circ\de^\alpha(x)
&=&\sum_{i=1}^{i=n+1}(-1)^{i+1}\de_i^\alpha\epsilon_1 (x)+ 
\sum_{i=1}^{i=n}(-1)^{i+1}\epsilon_1\de_i^\alpha (x)\\
&=&x + \sum_{i=2}^{i=n+1}(-1)^{i+1}\epsilon_1\de_{i-1}^\alpha (x)+ 
\sum_{i=1}^{i=n}(-1)^{i+1}\epsilon_1\de_i^\alpha (x)\\
&=&x+  \sum_{i=1}^{i=n}(-1)^{i}\epsilon_1\de_{i}^\alpha (x)+ 
\sum_{i=1}^{i=n}(-1)^{i+1}\epsilon_1\de_i^\alpha (x)\\
&=&x
\eeas

\epf

The previous proposition entails the following definition :

\bd\label{def_orientee} Let $\C$ be an $\omega$-category and
$\alpha\in\{-,+\}$.  We denote by $\omega Cat(I^n,\C)^\alpha$ the
subset of elements $x$ of $\omega Cat(I^n,\C)$ satisfying the
following conditions : 
\begin{itemize}
\item the element $x$ is a non degenerate element
of the cubical nerve
\item any element of the form
$\de^\alpha_{i_1}...\de^\alpha_{i_p}(x)$ is non degenerate in the
cubical nerve. 
\end{itemize}
Then $\de^\alpha(\Z\omega Cat(I^{*+1},\C)^\alpha)\subset
\Z\omega Cat(I^{*},\C)^\alpha$ by construction. We thus set
$$H_*^\alpha(\C,\Z)=H_*(\Z \omega Cat(I^{*},\C)^\alpha,\de^\alpha)$$
and we call these homology theories the \textit{negative (or positive
according to $\alpha$) corner homology} of $\C$. The cycles are called
the \textit{negative (or positive) corners} of $\C$. The maps $H_*^\pm $ induce
functors from $\omega Cat_1$ to $Ab$.  \ed

The second part of the definition is essential. Indeed if $a$ is a
$1$-dimensional morphism of $\C$, then the following element of
$\omega Cat(I^2,\C)$ is non degenerate although its image by $\de_1^-$
and $\de_2^-$ are degenerate elements of $Cat(I^1,\C)$ :

$$\xymatrix{s_0(a) \fr{a} \ff{rd}{a}&
t_0(a) \\
s_0(a) \fu{s_0(a)} \fr{s_0(a)}& 
s_0(a) \fu{a}}$$

The following proposition characterizes the elements of $\omega
Cat(I^{*},\C)^\alpha$

\bp\label{caracterisation} Assume that $x$ is an element of $\omega
Cat(I^{*},\C)$. Then $x$ is in $\omega Cat(I^{*},\C)^\alpha$ if and
only if all $x(\alpha ...0 ...\alpha)$ (the notation
$\alpha...0...\alpha$ meaning that $0$ appears only once) are
$1$-dimensional morphisms of $\C$.  \ep

\bpf If $x$ is in $\omega Cat(I^{n},\C)^\alpha$, then all
$\de^\alpha_{i_1}...\de^\alpha_{i_p}(x)$ are non degenerate in the
cubical nerve. But $y\in \omega Cat(I^1,\C)$ is non degenerate if and
only if $y(R(0))$ is $1$-dimensional. Hence the necessity of the
condition. Conversely assume that $x\notin \omega
Cat(I^{n},\C)^\alpha$. Then there exists $i$ between $1$ and $n$ such
that $\de^\alpha_{i_1}...\de^\alpha_{i_p} x=\epsilon_i(z)$ with $p<n$
and some $i_1$, ..., $i_p$ and with $z\in \omega Cat(I^{n-p-1},\C)$.
Then $$\{x(\alpha ...[-]_i ...\alpha),x(\alpha
...[0]_i...\alpha),x(\alpha ...[+]_i...\alpha)\}$$
is a singleton for
some $i$ therefore $x(\alpha ...[0]_i...\alpha)$ is $0$-dimensional.
Hence the sufficiency of the condition.  \epf

As the globular homology, the corner homologies do not yield functors
from $\omega Cat$ to $Ab$.

\subsection{Examples of corners}

\bp Let $\C$ be an $\omega$-category. The group $H_0^-(\C)$ is the
free abelian group generated by the final states of $\C$. The group
$H_0^+(\C)$ is the free abelian group generated by the initial states
of $\C$.  \ep

\bpf Obvious. \epf

\begin{figure}
\begin{center}
\setlength{\unitlength}{0.0005in}
\begingroup\makeatletter\ifx\SetFigFont\undefined%
\gdef\SetFigFont#1#2#3#4#5{%
  \reset@font\fontsize{#1}{#2pt}%
  \fontfamily{#3}\fontseries{#4}\fontshape{#5}%
  \selectfont}%
\fi\endgroup%
{\renewcommand{\dashlinestretch}{30}
\begin{picture}(624,1014)(0,-10)
\path(12,327)(312,927)
\path(285.167,806.252)(312.000,927.000)(231.502,833.085)
\path(12,327)(612,327)
\path(492.000,297.000)(612.000,327.000)(492.000,357.000)
\put(312,27){\makebox(0,0)[lb]{\smash{{{\SetFigFont{11}{14.4}{\rmdefault}{\mddefault}{\updefault}$y$}}}}}
\put(12,927){\makebox(0,0)[lb]{\smash{{{\SetFigFont{11}{14.4}{\rmdefault}{\mddefault}{\updefault}$x$}}}}}
\end{picture}
}
\end{center}
\caption{A $1$-dimensional negative corner}
\label{coin1}
\end{figure}
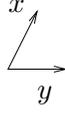

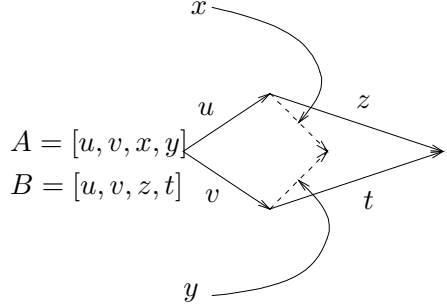
\begin{figure}
\begin{center}
\setlength{\unitlength}{0.0005in}
\begingroup\makeatletter\ifx\SetFigFont\undefined%
\gdef\SetFigFont#1#2#3#4#5{%
  \reset@font\fontsize{#1}{#2pt}%
  \fontfamily{#3}\fontseries{#4}\fontshape{#5}%
  \selectfont}%
\fi\endgroup%
{\renewcommand{\dashlinestretch}{30}
\begin{picture}(4512,3050)(0,-10)
\path(1800,1527)(2700,2127)
\path(2616.795,2035.474)(2700.000,2127.000)(2583.513,2085.397)
\path(1800,1527)(2700,927)
\path(2583.513,968.603)(2700.000,927.000)(2616.795,1018.526)
\path(2700,2127)(4500,1527)
\path(4376.671,1536.487)(4500.000,1527.000)(4395.645,1593.408)
\path(2700,927)(4500,1527)
\path(4395.645,1460.592)(4500.000,1527.000)(4376.671,1517.513)
\dashline{60.000}(2700,2127)(3300,1527)
\path(3193.934,1590.640)(3300.000,1527.000)(3236.360,1633.066)
\dashline{60.000}(2700,927)(3300,1527)
\path(3236.360,1420.934)(3300.000,1527.000)(3193.934,1463.360)
\path(2100,3027)        (2167.468,3011.427)
        (2232.337,2995.807)
        (2294.658,2980.109)
        (2354.482,2964.301)
        (2411.860,2948.354)
        (2466.845,2932.236)
        (2519.488,2915.917)
        (2569.841,2899.365)
        (2617.954,2882.549)
        (2663.880,2865.439)
        (2707.669,2848.004)
        (2749.374,2830.212)
        (2789.045,2812.033)
        (2826.735,2793.437)
        (2896.376,2754.866)
        (2958.709,2714.253)
        (3014.144,2671.350)
        (3063.095,2625.909)
        (3105.974,2577.685)
        (3143.191,2526.429)
        (3175.160,2471.894)
        (3202.292,2413.834)
        (3225.000,2352.000)

\path(3225,2352)        (3237.974,2296.980)
        (3240.524,2242.005)
        (3231.990,2185.540)
        (3211.714,2126.044)
        (3179.036,2061.980)
        (3133.297,1991.810)
        (3105.324,1953.955)
        (3073.838,1913.996)
        (3038.758,1871.742)
        (3000.000,1827.000)

\path(3056.888,1936.835)(3000.000,1827.000)(3101.882,1897.143)
\path(2100,27)  (2171.388,34.597)
        (2240.041,42.639)
        (2306.013,51.153)
        (2369.359,60.167)
        (2430.134,69.709)
        (2488.393,79.805)
        (2544.192,90.483)
        (2597.584,101.771)
        (2648.625,113.696)
        (2697.370,126.286)
        (2743.873,139.568)
        (2788.191,153.570)
        (2830.376,168.318)
        (2870.486,183.841)
        (2944.695,217.320)
        (3011.257,254.227)
        (3070.613,294.780)
        (3123.200,339.201)
        (3169.459,387.709)
        (3209.829,440.523)
        (3244.749,497.863)
        (3274.660,559.949)
        (3300.000,627.000)

\path(3300,627) (3316.639,695.231)
        (3319.568,761.880)
        (3307.907,828.705)
        (3280.778,897.465)
        (3237.301,969.916)
        (3209.158,1008.076)
        (3176.599,1047.817)
        (3139.513,1089.361)
        (3097.791,1132.926)
        (3051.324,1178.733)
        (3000.000,1227.000)

\path(3108.572,1167.737)(3000.000,1227.000)(3067.910,1123.616)
\put(1950,1902){\makebox(0,0)[lb]{\smash{{{\SetFigFont{11}{14.4}{\rmdefault}{\mddefault}{\updefault}$u$}}}}}
\put(2025,1002){\makebox(0,0)[lb]{\smash{{{\SetFigFont{11}{14.4}{\rmdefault}{\mddefault}{\updefault}$v$}}}}}
\put(1875,2952){\makebox(0,0)[lb]{\smash{{{\SetFigFont{11}{14.4}{\rmdefault}{\mddefault}{\updefault}$x$}}}}}
\put(1800,27){\makebox(0,0)[lb]{\smash{{{\SetFigFont{11}{14.4}{\rmdefault}{\mddefault}{\updefault}$y$}}}}}
\put(0,1527){\makebox(0,0)[lb]{\smash{{{\SetFigFont{11}{14.4}{\rmdefault}{\mddefault}{\updefault}$A=[u,v,x,y]$}}}}}
\put(3600,1977){\makebox(0,0)[lb]{\smash{{{\SetFigFont{11}{14.4}{\rmdefault}{\mddefault}{\updefault}$z$}}}}}
\put(3675,927){\makebox(0,0)[lb]{\smash{{{\SetFigFont{11}{14.4}{\rmdefault}{\mddefault}{\updefault}$t$}}}}}
\put(0,1077){\makebox(0,0)[lb]{\smash{{{\SetFigFont{11}{14.4}{\rmdefault}{\mddefault}{\updefault}$B=[u,v,z,t]$}}}}}
\end{picture}
}
\end{center}
\caption{A $2$-dimensional negative corner}
\label{coin2}
\end{figure}

There is in Figure~\ref{coin1} a very simple example of
$1$-dimensional negative corner. It consists of two $\omega$-functors
$x$ and $y$ from $I^1$ to $\C$ such that $x(R(-))=y(R(-))$ and such
that $x(R(0))$ and $y(R(0))$ are $1$-dimensional. Figure~\ref{coin2}
shows an example of $2$-dimensional negative corner. If we suppose $A$
and $B$ to be oriented such that $s_1A=u*_0x$, $t_1A=v*_0y$, $s_1 B=u
*_0z$ and $t_1 B=v*_0 t$, then $A-B$ is a negative corner.

\section{Filling of shells in the cubical singular nerve}\label{filling}

Now here is a technical tool which will enable us to construct some
operations on the cubical singular nerve of a globular
$\omega$-category (Section~\ref{corner_nerve} and
Section~\ref{morealgo}) and to construct the two Hurewicz morphisms
(Section~\ref{Hurewicz_oriente}). The notion of shell and of filling
of (thin or not) shells already appears in \cite{Brown_cube} in the
framework of $\omega$-groupoids, in \cite{phd-Al-Agl} in the framework
of cubical $\omega$-categories (see Definition~\ref{cubcat}).

\subsection{Recall about the freeness of $I^n$}

A key property of $I^n$ is its freeness. It means the following fact.
Let $p$ be some natural number. Let us call a realization
$(\Lambda^p,f_i)$ of $(\Lambda^p,E,B)$ in an $\omega$-category $\C$ a
family of maps $f_i$ from $(\Lambda^p)_i$ to $\tau_i\C$, where
$\tau_i\C$ is the $\omega$-category obtained by only keeping the cells
of $\C$ of dimension lower or equal than $i$.  The realization
$(\Lambda^p,f_i)$ is called $n$-extendable if there exists only one
functor $f$ from $\tau_nI^p$ to $\C$ such that for any $k\leqslant n$,
the commutative diagram holds

$$\xymatrix{
\tau_kI^p\ar@{->}[r]& \tau_nI^p \ar@{->}[d]^{f}\\
(\Lambda^p)_k \ar@{->}[u]^{R()} \ar@{->}[r]^{f_k}& {\C} \\
}$$

By convention, any realization $(\Lambda^p,f_i)$ is $0$-appropriate
and is already $0$-extendable since $f_0$ is already a functor from
$(I^p)_0$ to $\C$.  Suppose set up a notion of $n$-appropriate
realization. Then we have

$$\xymatrix{
(\Lambda^p)_{n+1} \ar@{->}[dr]^{f_{n+1}}\ar@{->}[d]_{R()} &\\
\tau_{n+1}I^p\ar@{-->}[r] \ar@<1ex>@{->}[d]^{t_n} \ar@<-1ex>@{->}[d]_{s_n}& {\C}
\ar@<1ex>@{->}[d]^{t_n} \ar@<-1ex>@{->}[d]_{s_n} \\
\tau_{n}I^p \ar@{->}[r]_{f} & {\C}\\
}$$

A realization $(\Lambda^p,f_i)$ is $(n+1)$-appropriate if $s_n
f_{n+1}=f s_n R$ and $t_n f_{n+1}=f t_n R$, and in this case,
$(\Lambda^p,f_i)$ is $(n+1)$-extendable (cf \cite{CPS} page 224).
Thus to construct an $\omega$-functor from $I^n$ for some natural
number $n$ to an $\omega$-category $\C$, it suffices to construct a
realization which is appropriate in all dimensions.

\subsection{Filling of shells using the freeness of $I^n$}

\bd An element $x\in \omega Cat(I^n,\C)$ is \textit{thin} if
$x(R(0_n))$ is of dimension strictly lower than $n$. An element which
is not thin is called \textit{thick}.  \ed

\bd A \textit{$n$-shell} in the cubical singular nerve is a family of
$2(n+1)$ elements $x_i^\pm $ of $\omega Cat(I^n,\C)$ such that
$\de_i^\alpha x_j^\beta= \de_{j-1}^\beta x_i^\alpha$ for $1\leqslant
i< j\leqslant n+1$ and $\alpha,\beta\in\{-,+\}$.  \ed

\bd
The $n$-shell $(x_i^\pm )$ is \textit{fillable}  if
\begin{enumerate}
\item the sets $\{x_i^{(-)^i},1\leqslant i\leqslant n+1\}$ and
$\{x_i^{(-)^{i+1}},1\leqslant i\leqslant n+1\}$ have each one exactly
one thick element and if the other ones are thin.
\item 
if $x_{i_0}^{(-)^{i_0}}$ and $x_{i_1}^{(-)^{i_1+1}}$ are these two
thick elements then there exists $u\in \C$ such that
$s_n(u)=x_{i_0}^{(-)^{i_0}}(0_n)$ and
$t_n(u)=x_{i_1}^{(-)^{i_1+1}}(0_n)$.
\end{enumerate}
\ed

The main proposition of this section is the following one, which is an
analogue of \cite{phd-Al-Agl} Proposition 2.7.3.

\bp\label{remplissage} Let $(x_i^\pm )$ be a fillable $n$-shell with
$u$ as above.  Then there exists one and only one element $x$ of
$\omega Cat(I^{n+1},\C)$ such that $x(0_{n+1})=u$, and for $1\leqslant
i\leqslant n+1$, and $\alpha\in\{-,+\}$ such that $\de_i^\alpha x=
x_i^\alpha$.  \ep

\bpf The underlying idea of the proof is as follows. If one wants to define an
$\omega$-functor from $I^{n+1}$ to an $\omega$-category $\C$, it
suffices to construct $2(n+1)$ $\omega$-functors from the $2(n+1)$
$n$-faces of $I^{n+1}$ to $\C$ which coincide on the intersection of
their definition domains and to fill correctly the interior of
$I^{n+1}$.

We have necessarily $x(k_1\dots[\pm ]_i\dots k_n)=x_i^\pm (k_1...k_n)$
and $x(0_{n+1})=u$. Therefore there is at most one such realization
$(\Lambda^{n+1},x_i)$.  It suffices to show that $x$ is
$(n+1)$-extendable. It is certainly $0$-appropriate therefore
$0$-extendable. Suppose we have proved that this realization is
$p$-appropriate and therefore $p$-extendable for every $p<n+1$.
First suppose that $p<n$. We want to prove that $(\Lambda^{n+1},x_i)$
is $(p+1)$-appropriate, i.e. that $s_p x(k_1\dots k_{n+1})=x s_p
R(k_1\dots k_{n+1})$ and $t_p x(k_1\dots k_{n+1})=x t_p R(k_1\dots
k_{n+1})$ for $k_1\dots k_{n+1}$ of dimension $p+1$.

Let us verify the first equality. We have
\beas
s_p x(k_1\dots k_{n+1})
&=&s_p (\de_i^{k_i} x)(k_1\dots \widehat{k_i}\dots k_{n+1})\hbox{ since $p+1<n+1$}\\
&=& s_p x_i^{k_i} (k_1\dots \widehat{k_i}\dots k_{n+1})\\
&=& x_i^{k_i} s_p R(k_1\dots \widehat{k_i}\dots k_{n+1})\hbox{ since the  $x_i^\alpha$ are $\omega$-functors}
\eeas
But $s_p R(k_1\dots \widehat{k_i}\dots k_{n+1})=\Psi(R(X_1),...,R(X_s))$ where 
$\Psi$ consists only of compositions of $X_h$ of dimension lower than $p$. Then 
\beas 
s_p x(k_1\dots k_{n+1})
&=&\Psi(x_i^{k_i} R(X_1),...,x_i^{k_i}R(X_s))\\
&=&\Psi((\de_i^{k_i}x)R(X_1),...,(\de_i^{k_i}x)R(X_s))\\
&=&x\Psi(\delta_i^{k_i}R(X_1),...,\delta_i^{k_i} R(X_s))\hbox{ since $x$ is $p$-extendable}\\
&=&x s_p R(k_1\dots k_{n+1})\hbox{ since $\delta_i^{k_i}$ is an $\omega$-functor}
\eeas
It remains to prove that $x$ is $(n+1)$-extendable. 
We have to prove that $s_n x(0_{n+1})=x s_n R(0_{n+1})$ and 
$t_n x(0_{n+1})=x t_n R(0_{n+1})$. Let us verify the first equality. 
We have $s_n R(0_{n+1})=\Psi'(\delta_{i_0}^{(-)^{i_0}}(0_n),Y_1,\dots,Y_s)$
where $\Psi'$ is a function containing only composition maps and where 
$Y_1,\dots,Y_s$ are elements of $I^{n+1}$ of the form $R(\{x\})$ where
$x$ is a face of the $(n+1)$-cube. Then 
\beas
x s_n R(0_{n+1})
&=&x\Psi'(\delta_{i_0}^{(-)^{i_0}}(0_n),Y_1,\dots,Y_s)\\
&=&\Psi'(x (\delta_{i_0}^{(-)^{i_0}}(0_n)),x(Y_1),\dots,x(Y_s))\hbox{ since $x$ is $n$-extendable}\\
&=&x(\delta_{i_0}^{(-)^{i_0}}(0_n))\hbox{ since $dim(x(Y_1)),\dots,dim(x(Y_s))<n$}\\
&=&x_{i_0}^{(-)^{i_0}}(0_n)\\
&=&s_n(u) \hbox{ by definition of $u$}\\
&=&s_n x(0_{n+1})\hbox{ by definition of $x$}
\eeas
\epf

\section{Two new simplicial nerves}\label{corner_nerve}

As an immediate application of Section~\ref{filling} , we construct two
families of connections on the cubical singular nerve of any
$\omega$-category. They will be useful in the sequel.

\subsection{Cubical set with connections}

\bd\label{connections}\cite{phd-Al-Agl}
A \textit{cubical set with connections} consists of a cubical set
$$((K_n)_{n\geqslant 0},\de_i^\alpha,\epsilon_i)$$ together with two
additional families of degeneracy maps
$$\Gamma_i^\alpha:\xymatrix{K_{n}\fr{} & K_{n+1}} $$ with
$\alpha\in\{-,+\}$, $n\geqslant 1$ and $1\leqslant i \leqslant n$ and satisfying the
following axioms :
\begin{enumerate}
\item $\de_i^\alpha\Gamma_j^\beta=\Gamma^\beta_{j-1}\de_i^\alpha$ 
  for all $i<j$ and all $\alpha,\beta\in\{-,+\}$
\item $\de_i^\alpha\Gamma_j^\beta=\Gamma^\beta_{j}\de_{i-1}^\alpha$ 
  for all $i>j+1$ and all $\alpha,\beta\in\{-,+\}$
\item $\de_j^\pm \Gamma_j^\pm = \de_{j+1}^\pm \Gamma_j^\pm =Id$
\item $\de_j^\pm \Gamma_j^\mp = \de_{j+1}^\pm \Gamma_j^\mp= \epsilon_{j}\de_{j}^\pm $
\item $\Gamma_i^\pm \Gamma_j^\pm =\Gamma_{j+1}^\pm 
\Gamma_i^\pm $ if $i\leqslant j$\
\item $\Gamma_i^\pm \Gamma_j^\mp= \Gamma_{j+1}^\mp\Gamma_i^\pm $  if $i<j$
\item $\Gamma_i^\pm \Gamma_j^\mp= \Gamma_{j}^\mp\Gamma_{i-1}^\pm $  if $i>j+1$
\item $\Gamma_i^\pm \epsilon_j=\epsilon_{j+1}\Gamma_i^\pm $ if $i<j$
\item $\Gamma_i^\pm \epsilon_j=\epsilon_i\epsilon_i$  if $i=j$
\item $\Gamma_i^\pm \epsilon_j=\epsilon_j\Gamma_{i-1}^\pm $ if $i>j$
\end{enumerate}
There is an obviously defined small category $\Gamma$, such that the
category of cubical sets with connections is exactly the category
of presheaves over $\Gamma$. Hence the category of cubical set
with connections is denoted by $Sets^{\Gamma^{op}}$. \ed

The category of cubical sets with connections equipped with the
non $1$-contracting morphisms is denoted by $Sets_1^{\Gamma^{op}}$.

Looking back to the cubical singular nerve of a topological space, we
can endow it with connections as follows : \beas
&&\Gamma_i^-(f)(x_1,\dots,x_p)=f(x_1,\dots,max(x_i,x_{i+1}),\dots,x_p)\\ 
&&\Gamma_i^+(f)(x_1,\dots,x_p)=f(x_1,\dots,min(x_i,x_{i+1}),\dots,x_p)
\eeas

By keeping in a cubical set with connections only the morphisms
$\Gamma_i^-$, or by exchanging the role of the face maps $\de_i^+$
and $\de_i^-$ and by keeping only the morphisms $\Gamma_i^+$, we
obtain exactly a cubical set with connections in the sense of
Brown-Higgins.

\subsection{Construction of connections on the cubical singular nerve}

\bth\label{newnerve} Let $\C$ be an $\omega$-category and let $n$ be a
natural number greater or equal than $1$. For any $x$ in $\omega
Cat(I^n,\C)$ and for any $i$ between $1$ and $n$, we introduce two
realizations $\Gamma_i^-(x)$ and $\Gamma_i^+(x)$ from $\Lambda^{n+1}$
to $\C$ by setting \beas
\Gamma_i^-(x)(k_1...k_{n+1})=x(k_1...max(k_i,k_{i+1})...k_{n+1})\\ 
\Gamma_i^+(x)(k_1...k_{n+1})=x(k_1...min(k_i,k_{i+1})...k_{n+1}) \eeas
where the set $\{-,0,+\}$ is ordered by $-<0<+$. Then $\Gamma_i^-(x)$
and $\Gamma_i^+(x)$ yield $\omega$-functors from $I^{n+1}$ to $\C$,
meaning two elements of $\omega Cat(I^{n+1},\C)$.  Moreover, in this
way, the cubical nerve of $\C$ is equipped with a structure of cubical
complex with connections.  \eth

\bpf The construction of  $\Gamma_i^\pm $ is exactly the same as the
one of connections on the cubical singular nerve of a topological
space. Thus there is nothing to verify in the axioms of cubical
complex with connections except the relations mixing the two families
of degeneracies $\Gamma_i^+$ and $\Gamma_i^-$ : all other axioms are
already verified in \cite{Brown_cube}. Therefore it remains to verify
that $\Gamma_i^\pm \Gamma_j^\mp= \Gamma_{j+1}^\mp\Gamma_i^\pm $ if
$i<j$, and $\Gamma_i^\pm \Gamma_j^\mp=
\Gamma_{j}^\mp\Gamma_{i-1}^\pm $ if $i>j+1$, that can be done quickly.
The only remaining point to be verified is that all  realizations 
$\Gamma_i^\pm (x) $ yield $\omega$-functors.

If $x$ is an element of $\omega Cat(I^1,\C)$, then we can depict
$\Gamma_1^-(x)$ as in Figure~\ref{gamma1moins} and $\Gamma_1^+(x)$ as
in Figure~\ref{gamma1plus}. We see immediately that $\Gamma_1^-(x)$
and $\Gamma_1^+(x)$ yield two elements of $\omega Cat(I^2,\C)$.

\begin{figure}
$$\xymatrix{x(+)\fr{x(+)}\ff{rd}{x(0)}&x(+)\\x(-)\fu{x(0)}\fr{x(0)}&
x(+)\fu{x(+)}}$$
\caption{$\Gamma_1^-(x)$}
\label{gamma1moins}
\end{figure}

\begin{figure}
$$\xymatrix{x(-)\fr{x(0)}\ff{rd}{x(0)}&x(+)\\x(-)\fu{x(-)}\fr{x(-)}&
x(-)\fu{x(0)}}$$
\caption{$\Gamma_1^+(x)$}
\label{gamma1plus}
\end{figure}

Suppose we have proved that all realizations $\Gamma_i^\pm (x)$ for
$1\leqslant i\leqslant n$ yield $\omega$-functors if $x\in \omega
Cat(I^n,\C)$. Now we want to prove that all realizations
$\Gamma_i^\pm (y)$ for $1\leqslant i\leqslant n+1$ yield $\omega$-functors if
$y\in \omega Cat(I^{n+1},\C)$.

Because of the axioms of cubical set with connections and because of
the induction hypothesis, the $(\de_j^\alpha\Gamma_i^\pm (y))$ are
$\omega$-functor from $I^n$ to $\C$. The family
$(\de_j^\alpha\Gamma_i^\pm (y))$ is also a $n$-shell. We can fill it
in a canonical way because the top dimensional elements are the same.
\epf

We denote by $Sets^{\Delta^{op}}$ the category of simplicial sets
\cite{May}. If $A$ is a simplicial set, the axioms of simplicial
sets imply that $C(A)=(\Z A_*,\de=\sum (-1)^i \de_i)$, where $\Z
A_n$ means the free abelian group generated by the set $A_n$, is a
chain complex. It is called the unnormalized chain complex of $A$.
The normalized chain complex of $A$ is the quotient chain complex
$N(A)=C(A)/D(A)$ where $D(A)$ is the sub-complex of $C(A)$ generated
by the degenerate elements. It turns out that the canonical morphism
of chain complex from $C(A)$ to $N(A)$ is a quasi-isomorphism
\cite{Weibel}.

As a consequence of the previous construction, we obtain two new
simplicial nerves.

\bp Let $\C$ be an $\omega$-category and $\alpha\in\{-,+\}$.
We set $$\mathcal{N}^\alpha_{n}(\C)=\omega Cat(I^{n+1},\C)^\alpha$$
and for all $n\geqslant 0$ and all $0\leqslant i\leqslant n$,
$$\xymatrix{{\de_i:\mathcal{N}^\alpha_{n}(\C)}\fr{}&{\mathcal{N}^\alpha_{n-1}(\C)}}$$
is the arrow $\de^\alpha_{i+1}$, and
$$\xymatrix{{\epsilon_i:\mathcal{N}^\alpha_{n}(\C)}\fr{}&{\mathcal{N}^\alpha_{n+1}(\C)}}$$
is the arrow $\Gamma^\alpha_{i+1}$. We obtain this way a simplicial
set $$(\mathcal{N}^\alpha_*(\C),\de_i,\epsilon_i)$$ called the
\textit{negative (or positive according to $\alpha$) corner simplicial nerve}
of $\C$. The non normalized chain complex associated to it gives exactly the
corner homology of $\C$ (in degree greater than or equal to $1$).
The maps $\mathcal{N}^\alpha$ induces a functor from $\omega Cat_1$ to
$Sets^{\Delta^{op}}$.   \ep

\bpf The axioms of simplicial sets are immediate consequences of
the axioms of cubical set with connections.
\epf

Notice that the indices are shifted by one. Intuitively, these
simplicial nerves consist of cutting an oriented $n$-hypercube by an
hyperplane close to a corner (the negative one or the positive one) :
the intersection we get is the oriented $(n-1)$-simplex in sense of
\cite{oriental}.

\section{The oriented Hurewicz morphisms}\label{Hurewicz_oriente}

In this section, we construct natural morphisms from the globular
homology of a $\omega$-category $\C$ to its two corner homology
theories. We call these maps the negative and positive oriented
Hurewicz morphisms. Intuitively, they map any oriented loop with
corners to its corresponding negative or positive corners (except for
the $0$-dimensional, see below).

\subsection{The $0$-dimensional case}

The projection from $\Z\C_0\p\Z\C_0$ to $\Z\C_0$ on the first
(resp. the second) component yields a natural group morphism from
$H_0^{gl}(\C)$ to $H_0^{-}(\C)$ (resp. $H_0^{+}(\C)$).  Indeed if
$$X=(s_0(x-y),t_0(x-y)),$$ then the first component (resp. the second
one) of $X$ induces $0$ on the corner homology.  We obtain thus a natural 
morphism
$h_0^\pm $ from $H_0^{gl}(\C)$ to $H_0^\pm (\C)$.

\subsection{The $1$-dimensional case}

If $x$ is a $1$-dimensional morphism of $\C$, let $\square_1(x)$ be the
element of $\omega Cat(I^1,\C)$ defined by 
$$\square_1(x)(R(-))=s_0(x),
\square_1(x)(R(0))=x,\square_1(x)(R(+))=t_0(x).$$ We extend
$\square_1$ by linearity.  If $z$ is a $2$-dimensional morphism of
$\C$, let $\square^-_2(z)$ be the element of $\omega Cat(I^2,\C)$
defined by the diagram
$$\xymatrix{t_0(z) \fr{t_0(z)} \ff{rd}{z}&
t_0(z) \\
s_0(z) \fu{s_1(z)} \fr{t_1(z)}& 
t_0(z) \fu{t_0(z)}}$$
and let $\square^+_2(z)$ be the element of $\omega Cat(I^2,\C)$ defined
by
$$\xymatrix{s_0(z) \fr{s_1(z)} \ff{rd}{z}&
t_0(z) \\
s_0(z) \fu{s_0(z)} \fr{s_0(z)}& 
s_0(z) \fu{t_1(z)}}$$
We extend $\square^-_2()$ and $\square^+_2()$ by linearity.
We get thus the following proposition

\bp The natural linear map $h_1^\pm $ from $\Z\C_1$ to 
$\Z\omega Cat(I^1,\C)^\pm $ which maps $x_1$ to $\square_1(x_1)$ induces a
natural map (still denoted by $h_1^\pm $) from $H_1^{gl}(\C)$ to
$H_1^\pm (\C)$. We call it the $1$-dimensional oriented Hurewicz
morphism.
\ep

\bpf The proof is quite simple. If $x_1$ is a $1$-dimensional globular
cycle, then $$\de^\pm \square_1(x_1)=\square_1(x_1)(\pm )=0$$ because of the
definition of $\square_1(x_1)$. And a $1$-dimensional globular boundary
$s_1(x_2)-t_1(x_2)$ is mapped to $\de^\pm (\square_2^-(x_2))$.
\epf

\subsection{The higher dimensional case}

\bp\label{construction_carre} For any natural number $n$ greater 
or equal than $2$, there exists a unique natural map $\square_n^-$ from
$\C_n$ (the $n$-dimensional cells of $\C$) to $\omega Cat(I^n,\C)$
such that
\begin{enumerate}
\item the equality $\square_n^-(x)(0_n)=x$ holds.
\item if $n\geqslant 3$ and $1\leqslant i\leqslant n-2$,  then
$\de_i^\pm \square_n^-= \Gamma_{n-2}^- \de_i^\pm \square^-_{n-1} s_{n-1}$.
\item if $n\geqslant 2$ and $n-1\leqslant i \leqslant n$, then
$\de_i^- \square_n^-=\square_{n-1}^- d_{n-1}^{(-)^i}$ and
$\de_i^+ \square_n^-= \epsilon_{n-1} \de_{n-1}^+ \square_{n-1}^- s_{n-1}$.
\end{enumerate}
Moreover for $1\leqslant i \leqslant n$, we have $\de_i^\pm \square_n^- s_n u=
\de_i^\pm \square_n^- t_nu$ for any $(n+1)$-morphism $u$.
\ep

\bpf Let $k_1\dots k_{n-1}\in \Lambda^{n-1}$. Then the natural map
$ev_{k_1\dots k_{n-1}}\de_i^\pm \square_n^-$ from $\C_n$ to $\C_{n-1}$
which sends $x\in\C_n$ to $(\de_i^\pm \square_n^- x)(k_1\dots
  k_{n-1})$ corresponds by Yoneda to an $\omega$-functor $f_{k_1\dots k_{n-1}}$
from $2_{n-1}$ to $2_n$, where $2_{n-1}$ (resp. $2_n$) is the free $\omega$-category
generated by a $(n-1)$-morphism $A$ (resp. a $n$-morphism $B$). Set 
$$f_{k_1\dots k_{n-1}}(A)=d_{n_{k_1\dots k_{n-1}}}^{\alpha_{k_1\dots k_{n-1}}}(B)$$
with still the convention $d^-=s$ and $d^+=t$. Then 
$f_{k_1\dots k_{n-1}}=d_{n_{k_1\dots k_{n-1}}}^{\alpha_{k_1\dots k_{n-1}}}$ because of the freeness
of $2_{n-1}$. Moreover, the inequality $n_{k_1\dots k_{n-1}}\leqslant n-1$ holds.
Therefore \beas
(\de_i^\pm \square_n^- s_n u)(k_1\dots k_{n-1})&=& d_{n_{k_1\dots k_{n-1}}}^{\alpha_{k_1\dots k_{n-1}}} s_n u\\
&=& d_{n_{k_1\dots k_{n-1}}}^{\alpha_{k_1\dots k_{n-1}}} t_n u\\
&=&(\de_i^\pm \square_n^- t_n u)(k_1\dots k_{n-1})
\eeas
Therefore for $1\leqslant i \leqslant n$, we have $\de_i^\pm 
\square_n^- s_nu = \de_i^\pm \square_n^- t_nu$ for any $(n+1)$-morphism $u$.

Suppose the proposition proved for $p<n$ with $n\geqslant 2$ and take
a $n$-dimensional morphism $x$. Set $h_i^\pm =\Gamma_{n-2}^- \de_i^\pm 
\square^-_{n-1}{s_{n-1}(x)}$ for $1\leqslant i\leqslant n-2$, and set
$h_i^-=\square^-_{n-1}{d_{n-1}^{(-)^i}x}$ and
$h_i^+=\epsilon_{n-1}\de_{n-1}^+\square^-_{n-1}{s_{n-1}(x)}$ for
$i\geqslant n-1$. We are going to verify that $(h_i^\pm )$ is a
fillable $(n-1)$-shell.  It is sufficient to prove that for any $i$
and any $j$ between $1$ and $n$, and any $\alpha,\beta\in\{-,+\}$,
the equality $\de_i^\alpha h_j^\beta=\de_{j-1}^\beta h_i^\alpha$ holds as soon as
$1\leqslant i<j\leqslant n$.

First treat the case $i<j\leqslant n-2$. We have
\beas
\de_i^\alpha h_j^\beta &=& \de_i^\alpha \Gamma_{n-2}^-\de_j^\beta \square^-_{n-1}{s_{n-1}(x)}\hbox{ since $j<n-1$}\\
&=& \Gamma_{n-3}^-\de_i^\alpha\de_j^\beta \square^-_{n-1}{s_{n-1}(x)} \hbox{ since $i<n-2$}\\
&=&\Gamma_{n-3}^-\de_{j-1}^\beta\de_i^\alpha \square^-_{n-1}{s_{n-1}(x)}\\
&=&\de_{j-1}^\beta h_i^\alpha\hbox{ since $i<n-1$}
\eeas

Now treat the case $i<j=n-1$. We have
\beas
\de_i^\pm h_j^-&=&\de_i^\pm \square^-_{n-1}{d_{n-1}^{(-)^{n-1}}x} \\
&=& \de_{j-1}^-\Gamma_{n-2}^- \de_i^\pm \square^-_{n-1}{d_{n-1}^{(-)^{n-1}}x}\\
&=& \de_{j-1}^- h_i^\pm 
\eeas
We have also 
\beas 
\de_i^\pm h_{n-1}^+ &=& \de_i^\pm \epsilon_{n-1} \de_{n-1}^+ \square^-_{n-1}{s_{n-1}(x)}\\
&=& \epsilon_{n-2} \de_i^\pm \de_{n-1}^+ \square^-_{n-1}{s_{n-1}(x)}\\
&=& \epsilon_{n-2} \de_{n-2}^+ \de_i^\pm \square^-_{n-1}{s_{n-1}(x)}\\
&=& \de_{n-2}^+ \Gamma_{n-2}^- \de_i^\pm \square^-_{n-1}{s_{n-1}(x)}\\
&=& \de_{n-2}^+ h_i^\pm 
\eeas

Now treat the case $i<j=n$ and $i<n-1$. We have 
\beas
\de_i^\pm h_{n}^- &=& \de_i^\pm \square^-_{n-1}{d_{n-1}^{(-)^{n-1}}(x)}\\
&=& \de_{n-1}^-\Gamma_{n-2}^-  \de_i^\pm \square^-_{n-1}{d_{n-1}^{(-)^{n-1}}(x)}\\
&=&\de_{n-1}^- h_i^\pm \eeas
and
\beas
\de_i^\pm h_{n}^+ &=& \de_i^\pm \epsilon_{n-1} \de_{n-1}^+ \square^-_{n-1}{s_{n-1}(x)}\\
&=&\epsilon_{n-2}\de_i^\pm \de_{n-1}^+ \square^-_{n-1}{s_{n-1}(x)}\\
&=&\epsilon_{n-2}\de_{n-2}^+\de_i^\pm \square^-_{n-1}{s_{n-1}(x)}\\
&=& \de_{n-1}^+ \Gamma_{n-2}^- \de_i^\pm \square^-_{n-1}{s_{n-1}(x)}\\
&=&\de_{n-1}^+ h_i^\pm \eeas

Finally treat the case $i=n-1$ and $j=n$. We have
\beas
\de_{n-1}^- h_{n}^-&=&\de_{n-1}^-
\square^-_{n-1}{d_{n-1}^{(-)^{n}}(x)}\\
&=& \de_{n-1}^- \square^-_{n-1}{d_{n-1}^{(-)^{n-1}}(x)}\\&=&
\de_{n-1}^- h_{n-1}^-\\
\de_{n-1}^- h_{n}^+&=&\de_{n-1}^- \epsilon_{n-1}\de_{n-1}^+
\square^-_{n-1}{s_{n-1}(x)}
\\&=&\de_{n-1}^+ \square^-_{n-1}{s_{n-1}(x)}\\&=& \de_{n-1}^+ \square^-_{n-1}{d_{n-1}^{(-)^{n-1}}(x)} \\&=&\de_{n-1}^+ h_{n-1}^-\\
\de_{n-1}^+ h_{n}^-&=&\de_{n-1}^+ \square^-_{n-1}{d_{n-1}^{(-)^n}(x)}
\\
&=& \de_{n-1}^- \epsilon_{n-1} \de_{n-1}^+ \square^-_{n-1}{s_{n-1}(x)}
\\
&=& \de_{n-1}^- h_{n-1}^+\\
\de_{n-1}^+ h_{n}^+&=&\de_{n-1}^+ \epsilon_{n-1} \de_{n-1}^+
\square^-_{n-1}{s_{n-1}(x)}\\
&=&\de_{n-1}^+ h_{n-1}^+
\eeas
\epf

\begin{cor} Let $\C$ be an $\omega$-category and let $n\geqslant 1$. Set
$\C_{\leqslant n}=\C_1 \cup \dots \cup \C_n$
There exists one and only one natural map $\square_n^-$ from 
$\C_{\leqslant n}$ to $\omega Cat(I^n,\C)^-$ such that the following
axioms hold :
\begin{enumerate}
\item if $n\geqslant 3$ and $1\leqslant i\leqslant n-2$, then
$\de_i^\pm \square_n^-= \Gamma_{n-2}^- \de_i^\pm \square^-_{n-1} s_{n-1}$.
\item if $n\geqslant 2$ and $n-1\leqslant i \leqslant n$, then
$\de_i^- \square_n^-=\square_{n-1}^- d_{n-1}^{(-)^i}$ and
$\de_i^+ \square_n^-= \epsilon_{n-1} \de_{n-1}^+ \square_{n-1}^- s_{n-1}$.
\end{enumerate}
Moreover for $1\leqslant i \leqslant n$, we have 
$\de_i^\pm \square_n^- s_n= \de_i^\pm \square_n^- t_n$.
\end{cor}

Let $D^-_*(\C)$ be the acyclic group complex generated by the
degenerate elements of the negative simplicial nerve of $\C$ with the
conventions $D^-_n(\C)\subset \Z \omega Cat(I^n,\C)$ and
$D^-_1(\C)=D^-_0(\C)=0$. We get thus the following proposition :

\bp The natural linear map $h_n^-$ from $\Z\C_n$ to $\Z\omega
Cat(I^n,\C)^-$ which sends $x_n$ to $\square^-_n x_n $ for $n\geqslant
1$ and which associates $(x_0,y_0)\in \Z\C_0\p \Z\C_0$ to $x_0\in
\Z\C_0$ induces a natural complex morphism 
$$h^- : \xymatrix{{C^{gl}_*(\C)} \fr{} & {\Z \omega Cat(I^*,\C)^-/D^-_*(\C)}}$$
\ep

\bpf Let $n\geqslant 2$ and let $x_n\in \Z\C_n$.  We have to compare
$\sum_{j=1}^{j=n} (-1)^{j+1} \de_j^- \square_n^- x_n$ and
$\square_{n-1}^-(s_{n-1}(x_n)-t_{n-1}(x_n))$ modulo elements of
$D^-_*(\C)$. We get 
\beas &&\sum_{j=1}^{j=n} (-1)^{j+1} \de_j^-
\square_n^-x_n\\ &&= \sum_{j=1}^{j=n-2} (-1)^{j+1} \Gamma_{n-2}^-
\de_j^- \square_{n-1}^- s_{n-1} x_n+ (-1)^n \square_{n-1}^-
(d^{(-)^{n-1}}_{n-1}(x_n)-d^{(-)^{n}}_{n-1}(x_n))\\ 
&&=\square_{n-1}^-(s_{n-1}(x_n)-t_{n-1}(x_n)) \hbox{ mod } D^-_*(\C)
\eeas 
Now let us treat the case $n=1$. Let $x_1\in \Z\C_1$. We
immediately see that $h_0^-(s_0(x_1),t_0(x_1))$ and $\de^-(\square_1^-(x_1))$ 
are equal.  \epf

\begin{cor}\label{hmoins} The natural 
  linear map $h_n^-$ from $\Z\C_n$ to $\Z\omega Cat(I^n,\C)^-$ which
  sends $x_n$ to $\square^-_n x_n $ for $n\geqslant 1$ and which
  associates $(x_0,y_0)\in \Z\C_0\p \Z\C_0$ to $x_0\in \Z\C_0$ induces
  a natural linear map from the globular homology to the negative
  corner homology of $\C$.
\end{cor}

\bpf It is due to the fact that for $n\geqslant 2$, the $n$-th
homology group of the quotient chain complex $\Z \omega
Cat(I^*,\C)^-/D^-_*(\C)$ is the $(n-1)$-th homology group of the
normalized chain complex associated to the corner simplicial nerve of
$\C$. \epf

Now let us expose the construction of $h_n^+$ (without proof).

\bp Let $\C$ be an $\omega$-category and let $n\geqslant 1$. Set
$\C_{\leqslant n}=\C_1 \cup \dots \cup \C_n$
There exists one and only one natural map $\square_n^+$ from 
$\C_{\leqslant n}$ to $\omega Cat(I^n,\C)^+$ such that the following
axioms hold :
\begin{enumerate}
\item if $n\geqslant 3$ and $1\leqslant i\leqslant n-2$, then
$\de_i^\pm \square_n^+= \Gamma_{n-2}^+ \de_i^\pm \square^+_{n-1} s_{n-1}$.
\item if $n\geqslant 2$ and $n-1\leqslant i \leqslant n$, then
$\de_i^+ \square_n^+=\square_{n-1}^+ d_{n-1}^{(-)^{i+1}}$ and
$\de_i^- \square_n^+= \epsilon_{n-1} \de_{n-1}^+ \square_{n-1}^+ s_{n-1}$.
\end{enumerate}
Moreover for $1\leqslant i \leqslant n$, we have 
$\de_i^\pm \square_n^+ s_n= \de_i^\pm \square_n^+ t_n$.
\ep

\bp\label{hplus} The natural linear map $h_n^+$ from $\Z\C_n$ to
$\Z\omega Cat(I^n,\C)^+$ which sends $x_n$ to $\square_n^+ x_n $ for
$n\geqslant 1$ and which sends $(x_0,y_0)\in \Z\C_0\oplus \Z\C_0$ to
$y_0$ induces a natural complex morphism 
$$h^+ : \xymatrix{{C^{gl}_*(\C)} \fr{} & {\Z \omega Cat(I^*,\C)^+/D^+_*(\C)}}$$
where $D^+_*(\C)$ is the sub-complex of $\Z \omega Cat(I^*,\C)^+$ generated
by the degenerate elements of the positive corner simplicial nerve.
Therefore $h_n^+$ induces a natural linear map from $H_n^{gl}(\C)$ to
$H_n^{+}(\C)$.  \ep

\section{Toward an ``oriented algebraic topology''}\label{ditop}

\subsection{Homotopic $\omega$-categories}

Now we want to speculate about the notion of homotopic $\omega$-categories.
We proceed like in algebraic topology by defining a notion of homotopy
between non $1$-contracting $\omega$-functors, and hence we deduce a
notion of homotopy equivalence of $\omega$-category.  We are obliged
to work with non $1$-contracting $\omega$-functors because of the
globular and corner homologies.

Intuitively, we could say that two non $1$-contracting $\omega$-functors
$f$ and $g$ from $\C$ to $\D$ are homotopic if $f(\C)$ and $g(\C)$
have the same ``oriented topology''. So a first attempt of definition
could be : the $\omega$-functors $f$ and $g$ are homotopic if for any
$x\in\C$, $f(x)\sim g(x)$. Unfortunatly, $f(x)$ and $g(x)$ do not
have necessarily the same dimension. So this definition does not make
sense, except if $x$ is $0$-dimensional : in this case, $f(x)\sim
g(x)$ means $f(x)=g(x)$.  It is plausible to think that if $f$ and $g$
were homotopic, then $C_*^{gl}(f)$ and $C_*^{gl}(g)$ would be two
chain homotopic morphisms from $C_*^{gl}(\C)$ to $C_*^{gl}(\D)$.  So
we propose this definition :

\bd\label{def_morphisme_homotope} Let $f$ and $g$ be two non
$1$-contracting $\omega$-functors from an $\omega$-category $\C$ to an
$\omega$-category $\D$. The morphisms $f$ and $g$ are
\textit{homotopic} if the following conditions hold :

\begin{enumerate}
\item for any $0$-dimensional $x$ of $\C$, one has $f(x)=g(x)$.
\item there exists a linear map $h_1$ from $\Z \C_1$ to 
$\Z\D_{2}$ such that $(s_1-t_1)h_1(x)=f(x)-g(x)$ for any $1$-morphism
$x$ of $\C$.
\item for any $n\geqslant 2$, there exists a linear map $h_n$ from
$\Z\C_n$ to $\Z\D_{n+1}$ such that for any $n$-morphism $x$ of $\C$, 
we have 
$$h_{n-1}(s_{n-1}-t_{n-1})(x)+(s_n-t_n)h_n(x)= f(x)-g(x) \hbox{ mod $\Z\D_{n-1}$}$$
\end{enumerate}
We denote this property by $f\sim_{(h_*)} g$ or more simply $f\sim g$
whenever it is not necessary to precise the homotopy map.
\ed

\bp The binary relation ``is homotopic to'' is an equivalence relation
on the collection of non $1$-contracting $\omega$-functors from a given
$\omega$-category $\C$ to a given $\omega$-category $\D$.  \ep

\bpf One has $f\sim f$ since $f\sim_{(0)} f$. If $f \sim_{(h_*)} g$, then
$g\sim_{(-h_*)} f$. Now suppose that $f\sim_{(h^1_*)} g$ and
$g\sim_{(h^2_*)}k$. Then $f\sim_{(h^1_*+h^2_*)}k$. \epf

\bp The homotopy equivalence of non $1$-contracting  $\omega$-functors is 
compatible with the composition of non $1$-contracting $\omega$-functors
in the following sense. Take a diagram in $\omega Cat_1$
{$$\xymatrix{{\C}\fr{f} &{\D}
\ar@<1ex>@{->}[r]^{g}\ar@{->}[r]_{h} &{\mathcal{E}} \fr{k} &
{\mathcal{F}}}$$} If $g\sim h$, then $g\circ f \sim h\circ f$, $k\circ
g\sim k\circ h$ and $k\circ g\circ f\sim k\circ h\circ f$.
\ep

\bpf Suppose that $g\sim_{H_*} h$. Define $H^f$ like this ($n\geqslant 1$)
\begin{enumerate}
\item if $x\in\C_n$ and if $dim(f(x))<n$, then $H^f_n(x):=0$
\item if $x\in\C_n$ and if $dim(f(x))=n$, then $H^f_n(x):=H_n(f(x))$.
\end{enumerate}
If $x$ is $1$-dimensional, then 
\beas &&g\circ f(x)- h\circ f(x)\\
&&=(s_1-t_1)H_1(f(x))\hbox{ since $f$ is non $1$-contracting}\\
&&=(s_1-t_1)H^f_1(x)
\eeas
If $x$ is of dimension $n$ greater than $2$, then either $dim(f(x))<n$ and
in this case 
\beas &&g\circ f(x)- h\circ f(x)\hbox{ mod $\Z\mathcal{E}_{n-1}$}\\&=&0\\
&=&H_{n-1}(s_{n-1}-t_{n-1})(f(x))+(s_n-t_n)H^f_n(x)\\
&=&H_{n-1}\circ f \circ (s_{n-1}-t_{n-1})(x) +(s_n-t_n)H^f_n(x) \hbox{ since $f$ is an $\omega$-functor}\\
&=&H^f_{n-1} (s_{n-1}-t_{n-1})(x) +(s_n-t_n)H^f_n(x)
\eeas
since $H_{n-1}\circ f \circ (s_{n-1}-t_{n-1})(x)=0$. Or $dim(f(x))=n$
and
in that case
\beas
&&g\circ f(x)- h\circ f(x)\hbox{ mod $\Z\mathcal{E}_{n-1}$}\\
&=&H_{n-1}(s_{n-1}-t_{n-1})(f(x))+(s_n-t_n)H^f_n(x)\\
&=&H_{n-1}\circ f \circ (s_{n-1}-t_{n-1})(x) +(s_n-t_n)H^f_n(x) \hbox{ since $f$ is an $\omega$-functor}
\eeas
Since $f(x)$ is $n$-dimensional, then $f(s_{n-1}(x))=s_{n-1} \circ f(x)$ and 
 $f(t_{n-1}(x))=t_{n-1} \circ f(x)$ are $(n-1)$-dimensional. Therefore 
$$g\circ f(x)- h\circ f(x)\hbox{ mod $\Z\mathcal{E}_{n-1}$ }=H^f_{n-1}(s_{n-1}-t_{n-1})(x) +(s_n-t_n)H^f_n(x)
$$
Therefore $g\circ f\sim_{H^f} h\circ f$.

Now define $^k\!H$ by $^k\!H_n(x):= k(H_n(x))$ for  $x\in\C_n$. Then for $x$
$1$-dimensional, we have 
\beas
&&k\circ g(x)- k\circ h(x)\\
&=& k\circ (s_1-t_1)H_1(x)\\
&=& (s_1-t_1) ^k\!H_1(x)\hbox{ since $k$ is an $\omega$-functor}
\eeas
And for $x$ of dimension greater than $2$, we have
\beas
&&k\circ g(x)- k\circ h(x) \hbox{ mod $\Z\mathcal{F}_{n-1}$}\\
& =& k\left(H_{n-1}(s_{n-1}-t_{n-1})(x)+(s_n-t_n)H_n(x)\right)\\
&=& \ ^k\!H_{n-1}(s_{n-1}-t_{n-1})(x) + (s_n-t_n)\ ^k\!H_n(x)\hbox{ since $k$ is an $\omega$-functor}
\eeas
Therefore $k\circ g\sim_{^k\!H_*} k\circ h$. \epf

\bp\label{cas_particulier} Let $f$ and $g$ be two non $1$-contracting
$\omega$-functors such that for any $x$ of dimension lower than $n$,
$f(x)$ and $g(x)$ are two homotopic $dim(x)$-dimensional morphisms.
Then $f$ and $g$ are homotopic as $n$-functor from $\tau_n \C$ to
$\tau_n \D$ (when we consider only the morphisms of dimension lower
than $n$). In other terms, the ``oriented topology'' is the same in
dimension lower than $n$.  \ep

\bpf For any $x$ of dimension $1\leqslant d\leqslant n$, there exists
$h_d(x)\in \Z\D_{d+1}$ such that $(s_d-t_d)(h_d(x))=f(x)-g(x)$. By
convention, we take $h_d(x)=0$ whenever $f(x)=g(x)$. However
$s_{d-1}(f(x)-g(x))=s_{d-1}(s_d-t_d)(h_d(x))=0$ and in the same way,
we have $t_{d-1}(f(x)-g(x))=t_{d-1}(s_d-t_d)(h_d(x))=0$.  Therefore
$h_{d-1}(s_{d-1}x)=h_{d-1}(t_{d-1}x)=0$ and
$$h_{d-1}(s_{d-1}-t_{d-1})(x)+(s_d-t_d)h_d(x)= f(x)-g(x)$$
\epf

\bd\label{def_cat_homotope} Let $\C$ and $\D$ be two
$\omega$-categories. They are \textit{homotopic} if and only if
there exists a non $1$-contracting $\omega$-functor $f$ from $\C$ to $\D$
and a non $1$-contracting $\omega$-functor $g$ from $\D$ to $\C$ such that
$f\circ g \sim Id_{\D}$ and $g\circ f \sim Id_{\C}$. We say that $f$
and $g$ are \textit{homotopy equivalences} between the two $\omega$-categories
$\C$ and $\D$.  \ed

\bp  The homotopy equivalence is an equivalence relation indeed on the
collection of $\omega$-categories.
\ep

\bpf This relation is obviously reflexible and symmetric. It remains
to prove the transitivity. Let us consider the following diagram in
$\omega Cat_1$ : {$$\xymatrix{{\C} \ar@<1ex>@{->}[r]^{f} &
\ar@{->}[l]^{g} {\D} \ar@<1ex>@{->}[r]^{h} & \ar@{->}[l]^{k}
{\mathcal{E}}}$$} and suppose that $g\circ f\sim Id_{\C}$, $f\circ
g\sim Id_{\D}$, $h\circ k\sim Id_{\mathcal{E}}$ and $k\circ h\sim
Id_{\mathcal{D}}$. Then $g\circ k \circ h \circ f\sim g \circ f\sim
Id_{\C}$ and $h\circ f\circ g\circ k\sim h\circ k\sim
Id_{\mathcal{E}}$.  \epf

Now some examples of homotopic $\omega$-categories.

\bp For any natural number $p\geqslant 1$ and $q\geqslant 1$, $2_p$ (the free
$\omega$-category generated by one $p$-morphism) and $2_q$ are
homotopic.  \ep

\bpf If $p=q$, it is trivial. Suppose that $p>q$.  Let $f$
be the only $\omega$-functor from $<A>=2_p$ to $<B>=2_q$ 
such that $f(A)=B$ and let $g$ be the unique functor from 
$<B>$ to $<A>$ such that $g(B)=s_q(A)$. Then $f\circ g = Id_{2_q}$
so $f\circ g$ and $Id_{2_q}$ are homotopic  as $\omega$-functors.
Now consider $g\circ f$ and $Id_{2_q}$. Set 
\beas
&&h_r=0 \hbox{ if }1\leqslant r<q\\
&&h_r(s_r A)=0\hbox{ and }h_r(t_r A)=s_{r+1} A\hbox{ if }q\leqslant r<p\\
&&h_r(A)=0\hbox{ if }r\geqslant p
\eeas

First suppose that $q=1$. Then we have $(s_1-t_1)h_1(s_1 A)=0=(g\circ f)(s_1 A) -s_1 A$ and
$(s_1-t_1)h_1(t_1 A)=s_1 A - t_1 A =(g\circ f)(t_1 A)-t_1A$, and for
any $1< r<p$, we have
\beas
&&h_{r-1}(s_{r-1}-t_{r-1})(s_r A )+(s_r-t_r)h_r(s_r A)\\
&=&-s_r A \\
&=&g\circ f (s_r A) - s_r A \hbox{ mod }(2_p)_{r-1}
\eeas
and
\beas
&&h_{r-1}(s_{r-1}-t_{r-1})(t_r A )+(s_r-t_r)h_r(t_r A)\\
&=& -s_r A + s_r A - t_r A\\
&=& s_1 A - t_r A \hbox{ mod }(2_p)_{r-1} \hbox{ since $r>1$}\\
&=&g\circ f (t_r A) - t_r A 
\eeas
In order to  complete the case $q=1$, now suppose that $r= p$. Then 
\beas
&&h_{r-1}(s_{r-1}-t_{r-1})(A)+(s_r-t_r)h_r(A)\\
&=&- A\\
&=&s_1 A - A \hbox{ mod }(2_p)_{p-1} \hbox{ since }p\geqslant 2\\
&=& g\circ f(A) - A 
\eeas
Now suppose that $q>1$. The different cases may be treated in the 
same way. First set $r=1$. Then 
$(s_1-t_1)h_1(s_1 A)=0=(g\circ f)(s_1 A) -s_1 A$ and 
$(s_1-t_1)h_1(t_1 A)=s_1 A - t_1 A =(g\circ f)(t_1 A)-t_1A$, and for
any $1< r<q$, we have
\beas
&&h_{r-1}(s_{r-1}-t_{r-1})(s_r A )+(s_r-t_r)h_r(s_r A)\\
&=&0\\
&=&g\circ f (s_r A) -s_r A
\eeas
and
\beas
&&h_{r-1}(s_{r-1}-t_{r-1})(t_r A )+(s_r-t_r)h_r(t_r A)\\
&=&0\\
&=&g\circ f (t_r A) -t_r A
\eeas
For $r=q$, we have 
\beas
&&h_{q-1}(s_{q-1}-t_{q-1})(s_q A )+(s_q-t_q)h_q(s_q A)\\
&=&0\\
&=&g\circ f (s_q A) -s_q A
\eeas
and
\beas
&&h_{q-1}(s_{q-1}-t_{q-1})(t_q A )+(s_q-t_q)h_q(t_q A)\\
&=&s_q A - t_q A\\
&=&g\circ f (t_q A) -t_q A
\eeas
For $q< r<p$, we have 
\beas
&&h_{r-1}(s_{r-1}-t_{r-1})(s_r A )+(s_r-t_r)h_r(s_r A)\\
&=&-s_r A\\
&=&s_q A -s_r A \hbox{ mod }(2_p)_{r-1}\hbox{ since }q\leqslant r-1\\
&=& g\circ f (s_r A) -s_r A
\eeas
and
\beas
&&h_{r-1}(s_{r-1}-t_{r-1})(t_r A )+(s_r-t_r)h_r(t_r A)\\
&=&-s_r A +s_r A -t_r A\\
&=& s_q A -t_r A \hbox{ mod }(2_p)_{r-1} \hbox{ since }q\leqslant r-1\\
&=& g\circ f (t_r A) -t_r A\hbox{ mod }(2_p)_{r-1}
\eeas
It remains the case $r=p$ : 
\beas
&&h_{p-1}(s_{p-1}-t_{p-1})(A )+(s_p-t_p)h_p(A)\\
&=&-A \\
&=& s_q A -A \hbox{ mod }(2_p)_{p-1}\hbox{ since }q\leqslant p-1\\
&=&g\circ f(A) - A
\eeas
\epf

\bd Let $\C$ be an $\omega$-category. Let $I$ and $F$ be some sets of
$0$-morphims of $\C$. The \textit{bilocalization} of $\C$ with respect to
$I$ and $F$ the sub-category of $\C$ consists of the $n$-morphims $f$
such that $s_0 f \in I$ and $t_0 f\in F$ with the induced structure of
$\omega$-category.  This $\omega$-category is denoted by $\C(I,F)$. \ed

In the sequel, the set $I$ will be always the set of initial states
and the set $F$ the set of final states of the considered
$\omega$-category.

Now we prove that the bilocalization of $\omega$-category $\C$ of
Figure~\ref{trou} with respect to $I$ (the intersection of $u$ and
$e$) and $F$ (the intersection of $z$ and $d$) is homotopic to the
$\omega$-category of Figure~\ref{G_1}, denoted by $G_1[A,B]$. Let $f$
be the unique $\omega$-functor from $\C(s_0(u),t_0(z))$ to $G_1[A,B]$
which maps any $1$-path homotopic to $\gamma_1$ to $A$ and any
$1$-path homotopic to $\gamma_4$ to $B$. Let $g$ be the unique
$\omega$-functor from $G_1[A,B]$ to $\C(s_0(u),t_0(z))$ which maps $A$
to $\gamma_1$ and $B$ to $\gamma_4$. Then $f \circ g = Id_{G_1[A,B]}$.
It remains to prove that $g\circ f$ and $Id_{\C(s_0(u),t_0(z))}$ are
homotopic $\omega$-functors. For any $1$-morphism $x$ of
$\C(s_0(u),t_0(z))$, $g\circ f(x)$ and $x$ are homotopic
$1$-morphisms. Let $h_1(x)$ be the element of $\Z\C(s_0(u),t_0(z))_2$
such that $(s_1-t_1)h_1(x)= g\circ f(x)-x$. Take $h_2=0$. We have to
verify that for any $2$-morphism $C$, we have
$$h_1(s_1-t_1)C+
(s_2-t_2)h_2 C = g\circ f (C)- C\hbox{ mod }\C(s_0(u),t_0(z))_1.$$
Suppose for example that $s_1C$ is homotopic to $\gamma_1$.  Then $h_1
s_1 C$ is the unique element of $\Z \C(s_0(u),t_0(z))_2$ such that
$(s_1-t_1)h_1 s_1 C=\gamma_1 - s_1 C$. And $h_1 t_1 C$ is the unique
element of $\Z \C(s_0(u),t_0(z))_2$ such that $(s_1-t_1)h_1 t_1
C=\gamma_1 - t_1 C$. Then $h_1(s_1-t_1)C$ is the unique element of $\Z
\C(s_0(u),t_0(z))_2$ such that $(s_1-t_1) h_1(s_1-t_1)C=t_1 C - s_1
C$.  Therefore $h_1(s_1-t_1)C=-C$. Hence the result.

Now here are some other examples without proof. In
Figure~\ref{SwissFlag}, the bilocalization of the depicted
$\omega$-categories with respect to their set of initial and final
states are homotopic.

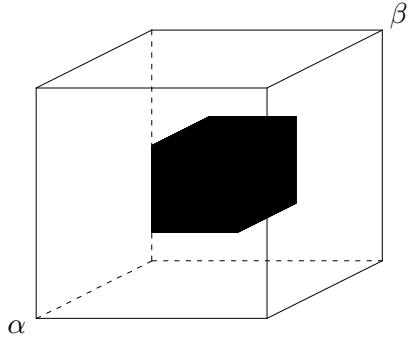
\begin{figure}
\begin{center}
\setlength{\unitlength}{0.0005in}
\begingroup\makeatletter\ifx\SetFigFont\undefined%
\gdef\SetFigFont#1#2#3#4#5{%
  \reset@font\fontsize{#1}{#2pt}%
  \fontfamily{#3}\fontseries{#4}\fontshape{#5}%
  \selectfont}%
\fi\endgroup%
{\renewcommand{\dashlinestretch}{30}
\begin{picture}(4067,3366)(0,-10)
\path(2700,168)(3900,768)
\path(2700,2568)(3900,3168)
\path(300,2568)(1500,3168)
\path(1500,3168)(3900,3168)(3900,768)
\dashline{60.000}(3900,768)(1500,768)(1500,3168)
\dashline{60.000}(300,168)(1500,768)
\whiten\path(1500,1068)(2400,1068)(3000,1368)
        (3000,2268)(2100,2268)(1500,1968)(1500,1068)
\path(1500,1068)(2400,1068)(3000,1368)
        (3000,2268)(2100,2268)(1500,1968)(1500,1068)
\path(300,168)(2700,168)
\path(300,168)(300,2568)
\path(300,2568)(2700,2568)
\path(2700,2568)(2700,168)
\blacken\path(1500,1068)(2400,1068)(3000,1368)
        (3000,2268)(2100,2268)(1500,1968)(1500,1068)
\path(1500,1068)(2400,1068)(3000,1368)
        (3000,2268)(2100,2268)(1500,1968)(1500,1068)
\put(0,18){\makebox(0,0)[lb]{\smash{{{\SetFigFont{11}{14.4}{\rmdefault}{\mddefault}{\updefault}$\alpha$}}}}}
\put(3975,3243){\makebox(0,0)[lb]{\smash{{{\SetFigFont{11}{14.4}{\rmdefault}{\mddefault}{\updefault}$\beta$}}}}}
\end{picture}
}
\end{center}
\caption{A $2$-semaphore}
\label{2semaphore}
\end{figure}

Figure~\ref{2semaphore} represents a cubical set with a
$3$-dimensional cubical hole. Then the bilocalization of this
$3$-category with respect to its sets of initial and final states is
homotopic to $G_2[A,B]$, the $2$-category of Figure~\ref{3hole}
generated by two non homotopic $2$-morphisms $A$ and $B$ having the
same $1$-source and the same $1$-target.

\subsection{Invariance of the globular and corner homologies}

\bth\label{globulaire_homotopie} Let $f$ and  $g$ be two 
non $1$-contracting $\omega$-functors from $\C$ to $\D$. Suppose that $f$
and $g$ are homotopic. Then for all natural number $n$, $f$ and $g$
induce linear maps from $H_n^{gl}(\C)$ to $H_n^{gl}(\D)$ and moreover
$H_{n}^{gl}(f)=H_{n}^{gl}(g)$.
\eth

\bpf Take two homotopic  $\omega$-functors $f$ and $g$. Let 
$x_1$ be a globular $1$-cycle. Then $(f-g)(x_1)=(s_1-t_1)h_1(x)$.
Therefore $(f-g)(x_1)=\de(h_1(x))$. Now take a globular $n$-cycle
$x_n$ with $n\geqslant 2$. Then $s_{n-1} x_n=t_{n-1} x_n$. Therefore
$f(x_n)-g(x_n)=(s_n-t_n)h_n x_n \hbox{ mod }\D_{n-1}$.
\epf

The analogous statement for the corner homologies is still a
conjecture only proved in the following particular case (see
Proposition~\ref{cas_particulier}) :

\bth\label{orientee_homotopie} Let $f$ and $g$ be two non 
$1$-contracting $\omega$-functor from $\C$ to $\D$. Suppose that for any
$x$ of dimension strictly lower than $n$, $f(x)=g(x)$ and such that for
$x$ $n$-dimensional, $f(x)$ and $g(x)$ are two homotopic $n$-morphisms
of $\D$. Then for any $p\leqslant n$, $f_p$ (resp. $g_p$) yield linear maps
from $H^\alpha_p(\C)$ to $H^\alpha_p(\D)$ for $\alpha\in\{-,+\}$ and
moreover, $H^\alpha_p(f)=H^\alpha_p(g)$.
\eth

\bpf We make the proof for $\alpha=-$. The homology of the non
normalized complex associated to a simplicial group is equal to its
homotopy (\cite{Weibel} Theorem 8.3.8). Therefore it suffices to
find an homotopy between $f(x)$ and $g(x)$ in $\mathcal{N}^-(\C)$ for any 
$x\in \omega Cat(I^n,\C)$. We can suppose without loss of generality that there
exists a $(n+1)$-dimensional morphism $u$ of $\C$ such that
$s_n(u)=f(x)(0_n)$ and $t_n(u)=g(x)(0_n)$. Then consider the following
realizations $h_i^\pm $ of $I^{n}$ for $1\leqslant i\leqslant n+1$ : for $i$
between $1$ and $n-1$, $h_i^\pm =\Gamma^-_{n-1}\de_i^\pm x =
\Gamma^-_{n-1}\de_i^\pm y$, $h^-_{n}=d_n^\epsilon(u)$,
$h^-_{n+1}=d_n^{\epsilon+1}(u)$ with $\epsilon$ equal to $0$ or $1$
depending on the parity of $n$, and finally
$h^+_{n}=h^+_{n+1}=\epsilon_{n}\de^+_{n}x=\epsilon_{n}\de^+_{n}y$. We
obtain the fillable $n$-shell which already appears in
Proposition~\ref{construction_carre}. The corresponding element of
$\omega Cat(I^{n+1},\C)$ yields an homotopy in $\mathcal{N}^-(\C)$
between $f(x)$ and $g(x)$.  \epf

\begin{conj}\label{invariance_corner} Let $f$ and  $g$ be two 
non $1$-contracting $\omega$-functors from $\C$ to $\D$. Suppose that $f$
and $g$ are homotopic. Then for all natural number $n$, $f$ and $g$
induce linear maps from $H_n^\pm (\C)$ to $H_n^\pm (\D)$ and moreover
$H_{n}^\pm (f)=H_{n}^\pm (g)$.
\end{conj}

\section{Some open questions and perspectives}\label{perspective}

\subsection{Some interesting problems in mathematics}

\begin{figure}
\[\xymatrix{
{\alpha} \rruppertwocell<10>{A} \rrlowertwocell<-10>{B} \ar[rr]&& {\beta}
}
\]\caption{Composition of two $2$-morphims}
\label{composition-of-2trans2}
\end{figure}
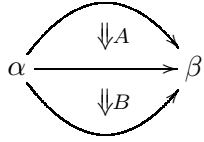

First of all, we come back to the globular homology of an
$\omega$-category. We propose here a small modification of its
definition.  Consider the $\omega$-category $\C$ of
Figure~\ref{composition-of-2trans2} where $A$ and $B$ are two
$2$-morphisms which are supposed to be composable.  Then
$H_2^{gl}(\C)\neq 0$ since $C_3^{gl}(\C)=0$ and since
$(s_1-t_1)(A*_1 B-A-B)=0$. However, this globular $2$-cycle
corresponds to nothing real in $\C$. As consequence of this small
calculation, we obtain that $H_2^{gl}(I^3)\neq 0$ (cf Figure~\ref{I3}). It suffices to
consider for example the globular $2$-cycle $C*_1 D-C-D$ with
$C=R(-00)*_0 R(0++)$ and $D=R(-0-)*_0 R(0+0)$.

\begin{figure}
\[
\xymatrix{
{\alpha}\rrtwocell<10>^u_v{_A} && {\beta} \rrtwocell<10>^x_y{^B}&& {\gamma}\\
}\]
\caption{Two $0$-composable $2$-morphisms}
\label{Two-0-composable-2-morphisms}
\end{figure}
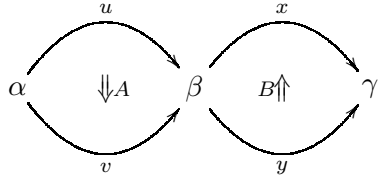

Consider the $\omega$-category of
Figure~\ref{Two-0-composable-2-morphisms} where $\alpha$, $\beta$ and
$\gamma$ are $0$-morphisms, $u$, $v$, $x$ and $y$ are $1$-morphisms
and $A$ and $B$ two $2$-morphims. There are at least two elements of
$\Z\C_2$ between $u *_0 x$ and $v*_0 y$ : $A *_0 y - u*_0 B$, $A*_0x
-v*_0B$. Therefore $A *_0 y - u*_0 B -A*_0x +v*_0B$ is a globular
$2$-cycle which means nothing geometrically.

So we propose to modify the definition of the globular homology as
follows. Let $\widehat{\omega Cat_1}$ be the category whose objects
are globular $\omega$-categories such that any $(n+1)$-morphism $X$ is
invertible with respect to $*_n$ as soon as $n\geqslant 1$ (i.e. there
exists a $(n+1)$-morphism $X^{-1}$ such that $s_n X^{-1}=t_n X$, $t_n
X^{-1}=s_n X$, $X *_n X^{-1}=s_n X$, $X^{-1} *_n X=t_n X$) and whose
morphisms are non $1$-contracting $\omega$-functors. Let us denote by
$\C\mapsto \widehat{\C}$ the left adjoint functor to the forgetful
functor from $\widehat{\omega Cat_1}$ to $\omega Cat_1$.  We set
$C^{gl}_0(\C)=\Z\C_0\oplus \Z\C_0$, $C^{gl}_1(\C)_1=\Z\C_1$, and
$C^{gl}_n(\C)$ for $n\geqslant 2$ is the free abelian group generated by
$\widehat{\C}_n$ quotiented by the relations 
$A + B=A *_{n-1} B \hbox{ mod $\Z\C_{n-1}$}$ 
if $A$ and $B$ are  two $n$-morphisms such that
$t_{n-1}A=s_{n-1}B$. With the same differential map, we obtain a new
globular homology theory. Let us denote it by $H_*^{new-gl}$.

With this new homology theory, the above problems disappear. It is
obvious for the first one and concerning the second one, here is the
reason. In $C^{new-gl}_2(\C)$ one has 
\beas 
A*_0x  -v*_0B&=&A *_0 x + (v*_0 B)^{-1}
\\&=&A *_0 x +v*_0 B^{-1}
\\&=& (A *_0 x) *_1 (v *_0 B^{-1})\\&=& (A *_1 v) *_0 (x *_1 B^{-1})\\&=& A *_0 B^{-1}\\
&=& (u *_1 A) *_0 (B^{-1} *_1 y)\\&=& (u*_0 B^{-1}) *_1 (A *_0 y)\\&=& u*_0 B^{-1} + A *_0 y\\&=& (u*_0 B)^{-1}+ A *_0 y\\
&=& A *_0 y - u*_0 B
\eeas

\begin{figure}
\[\xymatrix{
&\ar@{->}[rr]& &\ar@{->}[rd]&&&&
&\ar@{->}[dr]\ar@{->}[rr]&&\ar@{->}[dr]&\\
\ar@{->}[ru]\ar@{->}[rr]\ar@{->}[rd]&&\ar@{->}[ru]\ar@{->}[dr]\ff{lu}{\square_2^-(B)}&&&\ff{r}{A*_1B}&&
\ar@{->}[ur]\ar@{->}[dr]&&\ff{ur}{t_0A}\ar@{->}[rr]&&\\
&\ar@{->}[rr]\ff{ru}{\square_2^-(A)}&&\ar@{->}[ru]\ff{uu}{t_0A}&&&&
&\ar@{->}[ru]\ar@{->}[rr]\ff{uu}{\square_2^-(A*_1 B)}&&\ff{ul}{t_0A}\ar@{->}[ru]\\
}\]
\caption{$h_2^-(A*_1 B-A-B)$ is a boundary}
\label{bord}
\end{figure}
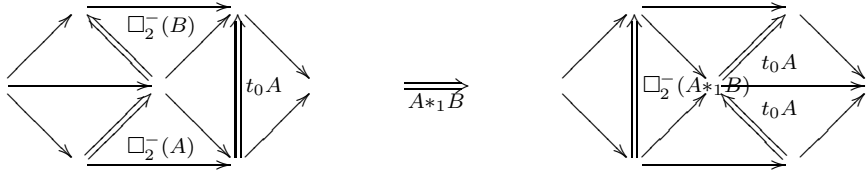

It suffices to consider the labelled $3$-cube of Figure~\ref{bord} to
see that $h_2^-(A*_1 B-A-B)$ is a boundary in the negative corner
homology. In the same way, we can prove that $h_2^+(A*_1 B-A-B)$ is
also a boundary, this time in the posivite corner homology. There is a
canonical linear map
$\xymatrix@1{H_*^{gl}(\C)\fr{}&H_*^{new-gl}(\C)}$. Therefore there
exists at least for $n=0,1,2$ a natural linear map $\widehat{h}_n^\pm $
such that the following diagram commutes :

\[
\xymatrix{H_n^{gl}(\C)\ar@{->}[rr]\ar@{->}[dd]^{h_n^\pm }&&H_*^{new-gl}(\C)\ar@{->}[dd]^{\widehat{h}_n^\pm }\\
&&\\
H_n^\pm (\C)\ar@{->}[rr]&&H_n^\pm (\widehat{\C})\\}
\]

We are led to the following conjectures :

\begin{conj}
\begin{enumerate}
\item The group $H_p^{new-gl}(I^n)$ vanishes for every natural number $n$ and any $p>0$.
\item For every $n\geqslant 2$, if $A$ and $B$ are two $n$-morphisms such that 
$t_{n-1}A=s_{n-1}B$, then $\square_n^\pm (A*_{n-1}B -A -B)$ is a boundary 
in the corresponding corner homology of $\widehat{\C}$. 
\item The canonical map $H_n^\pm (\C)\longrightarrow H_n^\pm (\widehat{\C})$
is an isomorphism for every $n\geqslant 0$
\item As consequence of the above two conjectures, $h_n^\pm $ factorizes
through the new globular homology theory.
\end{enumerate}
\end{conj}

This new definition of the globular homology gives rise to a new
definition of homotopic non $1$-contracting $\omega$-functors and gives
rise to the conjecture that the corner homology theories are still
invariant with respect to this new equivalence relation.

Elements like $A*_1 B-A-B$ could be called thin globular cycles. On the
corner homologies, the analogous elements are the linear combinations
of $x\in\omega Cat(I^n,\C)$ for some given natural number $n$ such
that $x(0_n)$ is of dimension strictly lower than $n$. We have
therefore the following conjecture :

\begin{conj} (About the thin elements of the corner complexes of a free globular 
$\omega$-category $\C$) Let $x_i$ be elements of $\omega
Cat(I^n,\C)^\pm $ and let $\lambda_i$ be natural numbers, where $i$
runs over some set $I$. Suppose that for any $i$, $x_i(0_n)$ is of
dimension strictly lower than $n$. Then $\sum_i \lambda_i x_i$ is a
boundary if and only if it is a cycle.
\end{conj}

With these new facts about the thin elements, Conjecture~\ref{invariance_corner} becomes

\begin{conj} Let $f$ and $g$ be two homotopic non $1$-contracting $\omega$-functors
from $\C$ to $\D$ where $\C$ and $\D$ are two free
$\omega$-categories. Let $x$ be an element of $\omega Cat(I^n,\C)^\pm
$. Then $f(x)-g(x)= B + T$ where $B$ is a boundary of $\Z \omega
Cat(I^n,\D)^\pm $ and $T$ a linear combination of thin elements of $\Z
\omega Cat(I^n,\D)^\pm $.
\end{conj}

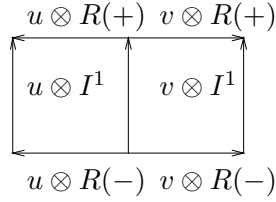
\begin{figure}
\begin{center}
\setlength{\unitlength}{0.0005in}
\begingroup\makeatletter\ifx\SetFigFont\undefined%
\gdef\SetFigFont#1#2#3#4#5{%
  \reset@font\fontsize{#1}{#2pt}%
  \fontfamily{#3}\fontseries{#4}\fontshape{#5}%
  \selectfont}%
\fi\endgroup%
{\renewcommand{\dashlinestretch}{30}
\begin{picture}(2424,1920)(0,-10)
\path(1212,411)(12,411)
\path(132.000,441.000)(12.000,411.000)(132.000,381.000)
\path(1212,411)(2412,411)
\path(2292.000,381.000)(2412.000,411.000)(2292.000,441.000)
\path(1212,411)(1212,1611)
\path(1242.000,1491.000)(1212.000,1611.000)(1182.000,1491.000)
\path(12,411)(12,1611)
\path(42.000,1491.000)(12.000,1611.000)(-18.000,1491.000)
\path(1212,1611)(12,1611)
\path(132.000,1641.000)(12.000,1611.000)(132.000,1581.000)
\path(1212,1611)(2412,1611)
\path(2292.000,1581.000)(2412.000,1611.000)(2292.000,1641.000)
\path(2412,411)(2412,1611)
\path(2442.000,1491.000)(2412.000,1611.000)(2382.000,1491.000)
\put(157,1011){\makebox(0,0)[lb]{\smash{{{\SetFigFont{11}{14.4}{\rmdefault}{\mddefault}{\updefault}$u\ot I^1$}}}}}
\put(157,36){\makebox(0,0)[lb]{\smash{{{\SetFigFont{11}{14.4}{\rmdefault}{\mddefault}{\updefault}$u\ot R(-)$}}}}}
\put(157,1761){\makebox(0,0)[lb]{\smash{{{\SetFigFont{11}{14.4}{\rmdefault}{\mddefault}{\updefault}$u\ot R(+)$}}}}}
\put(1537,1011){\makebox(0,0)[lb]{\smash{{{\SetFigFont{11}{14.4}{\rmdefault}{\mddefault}{\updefault}$v\ot I^1$}}}}}
\put(1537,36){\makebox(0,0)[lb]{\smash{{{\SetFigFont{11}{14.4}{\rmdefault}{\mddefault}{\updefault}$v\ot R(-)$}}}}}
\put(1537,1761){\makebox(0,0)[lb]{\smash{{{\SetFigFont{11}{14.4}{\rmdefault}{\mddefault}{\updefault}$v\ot R(+)$}}}}}
\end{picture}
}
\end{center}
\caption{Filling of corners}
\label{fillcorner}
\end{figure}

There exists a unique biclosed monoidal structure $\ot$ on $\omega
Cat$ such that $I^m \ot I^n=I^{m+n}$. See for example
\cite{Crans_Tensor_product} for an explicit construction using
globular pasting scheme theory. The effect of the functor $-\ot I^1$
is to fill corners as it is showed in Figure~\ref{fillcorner}.  We are
led to the following conjectures.

\begin{conj} 
  If $\C$ is an $\omega$-category, then the $\omega$-functor from $\C$
  to $I^1\ot \C$ which maps $u$ to $R(\mp)\ot u$ yield an isomorphism
  from $H_*^\pm (\C)$ to $H_*^\pm (I^1\ot\C)$ and the zero map from
  $H_*^\mp (\C)$ to $H_*^\mp (I^1\ot\C)$ for $*>0$.
\end{conj}

We can easily check that $u\mapsto R(\mp)\ot u$ induces $0$ from
$H_1^\mp (\C)$ to $H_1^\mp (I^1\ot\C)$. We end up this section with
three other problems and one remark about the corner homologies. The
conjectures are easy to verify in lower dimension.

\begin{conj}
\begin{enumerate}
\item The corner homology groups of $I^n$ vanish in dimension strictly
greater than $0$. In other terms, if $p>0$, then $H_p^\pm (I^n)=0$
\item Let $2_n$ be the free $\omega$-category generated by a $n$-morphism.
Then $p>0$ implies $H_p^\pm (2_n)=0$
\item Let $G_n$ be the oriented $n$-globe with $n>0$, i.e. the free $\omega$-category 
generated by two non homotopic $n$-morphisms having the same $(n-1)$-source
and the same $(n-1)$-target. Then $H_p^\pm (G_n)=0$ if $p\neq n$ and $p>0$.
Moreover  the equality $H_n^\pm (G_n)=\Z$ holds.
\end{enumerate}
\end{conj}

All the conjectures of this section will be the subject of future papers.

\subsection{Perspectives in computer science}

The study of the cokernel of the negative Hurewicz morphism would
allow us to detect the deadlock in concurrent machines. In an
analogous way the cokernel of the positive Hurewicz morphism would
allow us to detect the unreachable states in a concurrent
machine. This is useful for detecting the dead code in a concurrent
machines and for analyzing the safety properties of a machine
\cite{PGPWUsing} : proving that a property is false is equivalent to
proving that some states are unreachable.

We exhibited in Figure~\ref{SwissFlag} a $1$-category homotopic to a
$\omega$-category. The $1$-category we obtain suggests some relations
with the graph of oriented connected components introduced in
\cite{HDA2}.

We think also that some problems of confidentiality in computer
science involve the construction of a relative corner homology. The
problem stands as follows : take a concurrent machine with a flow of
inputs and a flow of outputs, every input and output having a
confidentiality level ; such a machine is confidential if the flow of
inputs of confidentiality level lower than $l$ determines the flow of
outputs of confidentiality level $l$ (otherwise an observer could
deduce from observations of outputs of confidentiality levels $l$ some
information about inputs of confidentiality level greater than $l$).
The geometric problem which arises from this situation stands as
follows : if some $n$-transitions are the inputs and some other ones
are the outputs, the problem is to know whether inputs determine
outputs over the set of all possible execution paths of the
machine. In the $1$-dimensional case, using bicomplexes, we already
found out a relation between this problem and the vertical and
horizontal $H_1$ and we suspect that in higher dimension this problem
is related in some way with the relative oriented Hurewicz morphisms
\cite{Gau_dim1_2} \cite{Gau_dim1_0} \cite{Gau_dim1_1}.

\section{Direct construction of the globular and corner homologies of a cubical set}\label{cubique_libre}\label{morealgo}

In this last section we explain how to obtain the globular and corner
homologies of a cubical set by using the free cubical
$\omega$-category generated by it, instead of considering the free
globular one. This approach could be useful in an algorithmic
viewpoint.

\subsection{Cubical $\omega$-category}

The notion of cubical $\omega$-category appears in the (already cited)
works of Brown, Higgins, Al-Agl.

\bd\label{cubcat} A \textit{cubical $\omega$-category} consists of a cubical set with
connections $$((K_n)_{n\geqslant
0},\de_i^\alpha,\epsilon_i,\Gamma_i^\alpha)$$ together with a family
of associative operations $+_j$ defined on $\{(x,y)\in K_n\p K_n,
\de^+_i x=\de_i^- y\}$ for $1\leqslant j\leqslant n$ such that
\begin{enumerate}
\item $(x+_j y)+_j z=x+_j(y+_j z)$
\item $\de^-_j(x +_j y)=\de^-_j(x)$
\item $\de^+_j(x +_j y)=\de^+_j(y)$
\item $\de^\alpha_i(x +_j y)=\left\{
\begin{CD}
&&\de^\alpha_i(x) +_{j-1} \de^\alpha_i(y)\hbox{ if $i<j$}\\
&&\de^\alpha_i(x) +_{j} \de^\alpha_i(y)\hbox{ if $i>j$}\\ \end{CD}
\right.$
\item $(x +_i y) +_j (z +_i t)=(x +_j z) +_i (y +_j t)$. We will denoted
the two members of this equality by $$
\left[\begin{array}{cc} x & z\\
y & t\\
\end{array}
\right] \coin{j}{i}
$$
\item $\epsilon_i(x +_j y)=\left\{
\begin{CD}&&\epsilon_i(x) +_{j+1} \epsilon_i(y) \hbox{ if $i\leqslant j$}\\
&&\epsilon_i(x) +_{j} \epsilon_i(y) \hbox{ if $i>j$}\\
\end{CD}\right.$
\item $\Gamma_i^\pm (x +_j y)=\left\{
\begin{CD}&&\Gamma_i^\pm (x) +_{j+1} \Gamma_i^\pm (y)\hbox{ if $i<j$}\\
&&\Gamma_i^\pm (x) +_{j} \Gamma_i^\pm (y)\hbox{ if $i> j$}
\end{CD}\right.$
\item If $i=j$, $\Gamma_i^- (x +_j y)=\left[\begin{array}{cc}
\epsilon_{j+1}(y) & \Gamma_j^- (y)\\
\Gamma_j^- (x) & \epsilon_j(y)\\
\end{array}\right] \coin{j+1}{j}$
\item If $i=j$, $\Gamma_i^+ (x +_j y)=\left[\begin{array}{cc}
\epsilon_{j}(x) & \Gamma_j^+ (y)\\
\Gamma_j^+ (x) & \epsilon_{j+1}(x)\\
\end{array}\right] \coin{j+1}{j}$
\item $\Gamma_j^+ x +_{j+1} \Gamma_j^- x = \epsilon_j x$ and $\Gamma_j^+ x +_{j} \Gamma_j^- x = \epsilon_{j+1} x$
\item $\epsilon_i \de_i^- x +_i x = x +_i \epsilon_i \de_i^+ x = x$
\end{enumerate}
The corresponding category with the obvious morphisms is denoted by
$\infty Cat$.
\ed

Look back again to the cubical singular nerve of a topological space
$X$.  We can equip it with operations $+_j$ as follows :

$$(f +_j g)(x_1,\dots,x_p) =
\left\{\begin{array}{c} f(x_1,\dots,2x_i,\dots,x_p)\hbox{ if $x_i\leqslant
1/2$}\\ g(x_1,\dots,2x_i-1,\dots,x_p)\hbox{ if $x_i\geqslant
1/2$}\end{array}\right.$$

All axioms of cubical $\omega$-categories are satisfied except the
associativity axiom.

It turns out that $\omega Cat$ and $\infty Cat$ are equivalent. The
$2$-dimensional case is solved in \cite{dim2_1} (which is followed by
\cite{dim2_2}) and the $3$-dimensional case is solved in
\cite{phd-Al-Agl}. Recently Richard Steiner developed the methods of 
Al-Agl to prove the result in all dimensions, as conjectured in
Al-Agl. The corresponding result for groupoids was already known from
earlier results of Brown-Higgins \cite{Brown-Higgins0}
\cite{Brown_cube}. The category equivalence is
realized by the functor $\gamma:\infty Cat \longrightarrow \omega Cat$
defined as follows ($G\in \infty Cat$) :

$$(\gamma G)_n=\{x\in G_n,\de_j^\alpha x\in
\epsilon_1^{j-1}G_{n-j}\hbox{ for $1\leqslant j\leqslant n$},
\alpha=0,1\}$$

Using general category theory arguments, one can prove that the
forgetful functor $U$ from $\infty Cat$ to $Sets^{\square^{op}}$ has a
left adjoint functor $\rho$ which defines therefore the free cubical
$\omega$-category generated by a cubical set. We will see in
Section~\ref{morealgo} that it can be constructed explicitely by
considering the cubical singular complex of the free globular
$\omega$-category generated by $K$.

\subsection{The globular and corner homologies of a cubical set}

First of all, take a look  at the cubical singular nerve :

\bp Let $\C$ be a globular $\omega$-category. 
For any strictly positive natural number $n$ and any $j$ between $1$
and $n$, there exists one and only one natural map $+_j$ from the set
of pairs $(x,y)$ of $\mathcal{N}^\square(\C)_n\p
\mathcal{N}^\square(\C)_n$ such that  $\de_j^+(x)=\de_j^-(x)$ to the
set $\mathcal{N}^\square(\C)_n$ which satisfies the following properties :
\bea
&&\de^-_j(x +_j y)=\de^-_j(x)\\
&&\de^+_j(x +_j y)=\de^+_j(x)\\
&&\de^\alpha_i(x +_j y)=\left\{
\begin{CD}
&&\de^\alpha_i(x) +_{j-1} \de^\alpha_i(y)\hbox{ if $i<j$}\\
&&\de^\alpha_i(x) +_{j} \de^\alpha_i(y)\hbox{ if $i>j$}\\ \end{CD}
\right.
\eea
Moreover, these operations induce a structure of cubical
$\omega$-category on $\mathcal{N}^\square(\C)$.
\ep

\bpf We give only a sketch of proof.

\textit{Step 1}. First of all, we observe that the
functor from $\omega Cat$ to the category $Sets$ of sets $$\C\mapsto
\omega Cat(I^p,\C)\p_j \omega Cat(I^p,\C)=\{(x,y)\in \omega
Cat(I^p,\C)\p \omega Cat(I^p,\C),\de_i^+x=\de_i^-y\}$$ is
representable. We denote by $I^p +_j I^p$ the representing
$\omega$-category.  It is equal to the direct limit of the diagram
$$\xymatrix{I^p && I^p\\ &
  I^{p-1}\ar@{->}[ul]_{\delta_i^+}\ar@{->}[ur]^{\delta_i^-}}$$
We denote by $\phi_-$ and $\phi_+$ the two canonical embeddings of $I^p$ 
in $I^p +_j I^p$ respectively in the first and the second term.

\textit{Step 2}. Using Yoneda, constructing a natural map
$$+_j:\xymatrix{\omega Cat(I^n,\C)\p_j \omega Cat(I^n,\C)\fr{}&\omega
  Cat(I^n,\C)}$$
is equivalent to construct an $\omega$-functor
$\eta_{n,j}$ from $I^n$ to $I^n +_j I^n$ satisfying the dual
properties. If $i<j$, then the natural transformation of functors
$(\de_i^\alpha,\de_i^\alpha)$ yields an $\omega$-functor
$(\delta_i^\alpha,\delta_i^\alpha)$ from $I^{n-1} +_{j-1} I^{n-1}$ to
$I^{n} +_{j} I^{n}$. It is easy to see that this $\omega$-functor
comes from the morphism of pasting scheme which associates to
$(\beta,k_1...k_{n-1})\in I^{n-1} +_{j-1} I^{n-1}$
$(\beta,k_1...[\alpha]_i...k_{n-1})\in I^{n} +_{j} I^{n}$ with
$\beta\in\{-,+\}$.  If $i>j$, then the natural transformation of
functors $(\de_i^\alpha,\de_i^\alpha)$ yields an $\omega$-functor
$(\delta_i^\alpha,\delta_i^\alpha)$ from $I^{n-1} +_{j} I^{n-1}$ to
$I^{n} +_{j} I^{n}$. It is easy to see that this $\omega$-functor
comes from the morphism of pasting scheme which associates to
$(\beta,k_1...k_{n-1})\in I^{n-1} +_{j} I^{n-1}$
$(\beta,k_1...[\alpha]_i...k_{n-1})\in I^{n} +_{j} I^{n}$ with
$\beta\in\{-,+\}$. The properties which are required for the
operations $+_j$ entail the following relations for the $\eta_{n,j}$ :
\bea
&&\eta_{n,j} \circ \delta_j^\pm = \phi_\pm \circ \delta_j^\pm \label{construction_+1}\\
&&\eta_{n,j} \circ \delta_i^\alpha =
(\delta_i^\alpha,\delta_i^\alpha)\circ \eta_{n-1,j-1}
\hbox{ if $i<j$}\label{construction_+2}\\
&&\eta_{n,j} \circ \delta_i^\alpha =
(\delta_i^\alpha,\delta_i^\alpha)\circ \eta_{n-1,j} \hbox{ if
  $i>j$}\label{construction_+3} \eea

\textit{Step 3}. The point is that it is difficult to find a formula for the
composition of all cells of $I^p +_j I^p$ (except in lower dimension).
It is simpler to find this formula in a free $\omega$-category
generated by a composable pasting scheme because composition means
union in such a context \cite{CPS}.  It turns out that $I^p +_j I^p$
is exactly the free globular $\omega$-category generated by the
composable pasting scheme defined as follows. Set $$(I^p +_j
I^p)_q=(\{-\}\p (I^p)_q)\cup (\{+\}\p (I^p)_q)/\equiv$$
where $\equiv$
is the equivalence relation induced by the binary relation
$$(-,k_1...[+]_j...k_{p-1})\equiv(+,k_1...[-]_j...k_{p-1})$$ for every
$k_1$, ..., $k_{p-1}$ in $\{-,0,+\}$ together with the binary
relations $E$ and $B$ defined by (with $x\in (I^p)_i$ and $y\in (I^p)_j$)
\beas
&&E^i_j=\{((a,x),(a,y))\in(I^p +_j I^p)\p (I^p +_j I^p)/x E^i_j y\hbox{ and }a\in\{-,+\}\}\\
&&B^i_j=\{((a,x),(a,y))\in(I^p +_j I^p)\p (I^p +_j I^p)/x B^i_j y\hbox{ and }a\in\{-,+\}\}
\eeas
Now we are in position to prove the following property $P(n)$ by
induction on $n$ : ``for any $j$ between $1$ and $n$, there exists one
and only one $\omega$-functor $\eta_{n,j}$ from $I^n$ to $I^n +_j I^n$
satisfying Condition~\ref{construction_+1},
Condition~\ref{construction_+2}, and Condition~\ref{construction_+3} ;
moreover $\eta_{n,j}(R(0_n))=R(\{(-,0_n),(+,0_n)\}$''. This latter
equality illustrates the interest of globular pasting schemes.

\textit{Step 4}. If $n=1$, we have to construct an $\omega$-functor from $I^1$ to $I^1
+_1 I^1$. The hypotheses lead us to set $\eta_{1,1}(R(-))=R((-,-))$
and $\eta_{1,1}(R(+))=R((+,+))$. There exists thus one and only one
suitable $\omega$-functor $\eta_{1,1}$ and this is the unique one
which satisfies $$\eta_{1,1}\left(R(0)\right)=R\left((-,0)\right)*_0
R\left((+,0)\right)=R\left((-,0),(+,0)\right).$$ So $P(1)$ is true.
Suppose we have proved $P(k)$ for $k<n$ where $n$ is a natural number
greater than $2$. We have to construct an $\omega$-functor $\eta_{n,j}$
for any $j$ between $1$ and $n$ from $I^n$ to $I^n +_j I^n$.  The
induction hypothesis and Condition~\ref{construction_+1}
Condition~\ref{construction_+2} and Condition~\ref{construction_+3}
entail the value of $\eta_{n,j}$ on all faces of $I^n$ of dimension at
most $n-1$.  It remains to prove that
$\eta_{n,j}(R(0_n))=R(\{(-,0_n),(+,0_n)\})$ is one and the only
solution. It suffices to verify that 
$s_{n-1}R(\{(-,0_n),(+,0_n)\})= \eta_{n,j}(s_{n-1} R(0_n))$
and that $t_{n-1}R(\{(-,0_n),(+,0_n)\})=\eta_{n,j}(t_{n-1}R(0_n))$. 
Let us verify the first equality. One has
\[s_{n-1}R(0_n)=R\left(\delta_1^{(-)^1}(0_{n-1}),\dots,\delta_n^{(-)^n}(0_{n-1})\right)\]
by the construction of $I^n$. By induction hypothesis, $\eta_{n,j}$ is $(n-1)$-extendable. Since composition
means union, and because of  Condition~\ref{construction_+1},
Condition~\ref{construction_+2}, and Condition~\ref{construction_+3}, one has 
\beas
&&\left[\bigcup_{h=1}^{h=j-1}(\delta_h^{(-)^h},\delta_h^{(-)^h})\circ \eta_{n-1,j-1}(0_{n-1})\right]
\cup \left( \phi_-\circ \delta_j^{(-)^j}(0_{n-1})\right)\\
&&\cup 
\left[\bigcup_{h=j+1}^{h=n}(\delta_h^{(-)^h},\delta_h^{(-)^h})\circ \eta_{n-1,j}(0_{n-1})\right]\\
&&=\left[\bigcup_{h=1}^{h=j-1}R\left((-,\delta_h^{(-)^h}(0_{n-1})),(+,\delta_h^{(-)^h}(0_{n-1}))\right)\right]\cup R\left(((-)^j,\delta_j^{(-)^j}(0_{n-1}))\right)\\
&&\cup \left[\bigcup_{h=j+1}^{h=n}R\left((-,\delta_h^{(-)^h}(0_{n-1})),(+,\delta_h^{(-)^h}(0_{n-1}))\right)\right] \\&&\subset\eta_{n,j}(s_{n-1} R(0_n))
\eeas

It suffices to verify that 
$$s_{n-1}R\left((- ,0_n),(+,0_n)\right)=R\left(\left\{((-)^j,\delta^{(-)^j}_j(0_{n-1})),(\pm ,\delta^{(-)^h}_h(0_{n-1}))/h\neq j\right\}\right)$$ in the pasting scheme 
$I^n +_j I^n$ to complete the proof.
\epf

Let us define a natural map $\square_n$ from $\tau_n\C$ (the set of
morphisms of $\C$ of dimension lower or equal than $n$) to $\omega
Cat(I^n,\C)$ by induction on $n$ as follows. One sets
$\square_0=\square_0^-$ and $\square_1=\square_1^-$.

\bp For any natural number $n$ greater or equal than $2$, there exists
a unique natural map $\square_n$ from $\C$ to $\omega Cat(I^n,\C)$ such
that
\begin{enumerate}
\item the equality $\square_n(x)(0_n)=x$ holds.
\item one has $\de_1^\alpha \square_n=\square_{n-1} d_{n-1}^{(-)^\alpha}$ for $\alpha=\pm$.
\item for $1<i\leqslant n$, one has $\de_i^\alpha \square_n = \epsilon_1
\de_{i-1}^\alpha \square_{n-1} s_{n-1}$.
\end{enumerate}
Moreover for $1\leqslant i \leqslant n$, we have $\de_i^\pm \square_n
s_n u= \de_i^\pm \square_n t_nu$ for any $(n+1)$-morphism $u$ and for
all $u\in\tau_n\C$, $\square_n(u)\in \gamma\mathcal{N}^\square(\C)_n$.
\ep

\bpf The induction equations define a fillable $(n-1)$-shell as defined
in Proposition~\ref{remplissage}. \epf

\bp For all $n\geq 0$, the evaluation map $ev_{0_n}:x\mapsto x(0_n)$
from $\omega Cat(I^n,\C)$ to $\C$ induces a bijection from $\gamma
\mathcal{N}^\square(\C)_n$ to $\tau_n\C$. \ep

\bpf Obvious for $n=0$ and $n=1$. Let us suppose that $n\geq 2$ and
let us proceed by induction on $n$. Since $ev_{0_n} \square_n (u)=u$
by the previous proposition, then the evaluation map $ev$ from
$\gamma\mathcal{N}^\square(\C)_n$ to $\tau_n\C$ is surjective. Now let us
prove that $x\in \gamma \mathcal{N}^\square(\C)_n$ and $y\in \gamma
\mathcal{N}^\square(\C)_n$ and $x(0_n)=y(0_n)=u$ imply $x=y$.  Since
$x$ and $y$ are in $\gamma \mathcal{N}^\square(\C)_n$, then one sees
immediately that the four elements $\de_1^\pm x$ and $\de_1^\pm y$ are
in $\gamma \mathcal{N}^\square(\C)_{n-1}$. Since all other
$\de_i^\alpha x$ and $\de_i^\alpha y$ are thin, then
$\de_1^-x(0_{n-1})=\de_1^-y(0_{n-1})=s_{n-1}u$ and
$\de_1^+x(0_{n-1})=\de_1^+y(0_{n-1})=t_{n-1}u$.  By induction
hypothesis, $\de_1^-x=\de_1^-y=\square_{n-1}(s_{n-1}u)$ and
$\de_1^+x=\de_1^+y=\square_{n-1}(t_{n-1}u)$. By hypothesis, one can
set $\de_j^\alpha x=\epsilon_1^{j-1} x^\alpha_j$ and $\de_j^\alpha
y=\epsilon_1^{j-1} y^\alpha_j$ for $2 \leqslant j\leqslant n$. And one
gets $x_j^\alpha=(\de_1^\alpha)^{j-1} \de_j^\alpha
x=(\de_1^\alpha)^{j}x= (\de_1^\alpha)^{j}y=y_j^\alpha$. Therefore
$\de_j^\alpha x=\de_j^\alpha y$ for all $\alpha\in\{-,+\}$ and all
$j\in[1,\dots,n]$. By Proposition~\ref{remplissage}, one gets $x=y$.
\epf

The above proof also shows that the map which associates to any cube
$x$ of the cubical singular nerve of $\C$ the cube
$\square_{dim(x)}(x(0_{dim(x)}))$ is exactly the usual folding
operator as exposed in \cite{phd-Al-Agl}.

Now let us remark that the free globular $\omega$-category generated
by a cubical set $K$ can be also obtained by considering the image by
the functor $\gamma$ of $\rho(K)$. Beware of the fact that in Al-Agl's
PhD, globular $\omega$-categories contain identity operators (his
$\omega$-categories are $\N$-graded). So the correct statement is
$(\gamma\rho(K))_n= \tau_nF(K)$ where $\tau_nF(K)$ is the $n$-category
obtained by keeping only the $p$-morphisms with $p\leqslant n$. It
suffices to prove the previous equality for $K=I^n$ since $\gamma$ is
a left adjoint functor
\cite{phd-Al-Agl}, therefore it commutes with all direct limits.

\begin{cor} Let $K$ be a cubical set. Then $\mathcal{N}^\square(F(K))$ is the
free cubical $\omega$-category $\rho(K)$ generated by $K$. \end{cor}

\bpf By \cite{phd-Al-Agl} Proposition 2.7.3, any $n$-cube $x$ of
$\rho(K)$ (resp. of $\mathcal{N}^\square(F(K))$) is determined by its
$(n-1)$-shell of $(n-1)$-faces $(\de_j^\pm x)_{1\leqslant j\leqslant
n+1}$ and by its image in $\gamma\rho(K)$ (resp. $\gamma
\mathcal{N}^\square(F(K))$).
\epf

Now take a cubical set $K$. As a consequence of the above remarks, it
is possible to construct $H_*^{gl}(K)$ and $H_*^\pm (K)$ and the two
morphisms $h_*^\pm $ by using the free cubical $\omega$-category
generated by $K$ instead of using the globular one. Let us still
denote by $\gamma \rho(K)$ the globular $\omega$-category obtained by
removing all identity elements.  It is exactly the free globular
$\omega$-category generated by $K$. We set
$H_*^{gl}(K):=H_*^{gl}(\gamma \rho(K))$ and since
$\mathcal{N}^\square(F(K))$ is the free cubical $\omega$-category
generated by $K$, we set $H_*^\pm (K)=H_*(\Z\rho(K)_*^\pm ,\de^\pm )$
where $$\rho(K)_n^\pm =\{x\in\rho(K)_n, \forall
i_1,\dots,i_{n-1},\de_{i_1}^\pm \dots \de_{i_{n-1}}^\pm x\in
\rho(K)_1\}.$$ The two morphisms $h_*^\pm $ from $H_*^{gl}(K)$ to
$H_*^\pm (K)$ are constructed like in
Proposition~\ref{construction_carre} : the only tool to be used is
again \cite{phd-Al-Agl} Proposition 2.7.3.

\section{Acknowledgments}

Thanks are due to Eric Goubault for an uncalculable number of
discussions about the geometry of automata, to Ronnie Brown and the
referees for their helpful comments, to Emmanuel Peyre for the picture
of the oriented $4$-cube, and to Marc Wambst for his remarks about the
redaction of the paper.

\vfill
%\newpage 

\section{The categories and functors of this work}

%\newpage
\[
\xymatrix{
{Sets^{\square^{op}}} \ar@{->}[dd]_F
&&  \ar@{->}[dd]_F\ar@{_{(}->}[ll] {Sets^{\square^{op}}_1}
&&
&&\\
&&&&&&\\
{\omega Cat}\ar@{->}[dd]_{\mathcal{N}^\square}
&& \ar@{->}[dd]_{\mathcal{N}^\square}\ar@{_{(}->}[ll] {\omega Cat_1}\ar@{-->}[rr]_{Gl(-)}\ar@/^30pt/[rrrr]^{H_*^{gl}}\ar[ddrr]_{C_*^{gl}}\ar@/_20pt/[rrrrdddddd]^{H_*^\alpha}
&&[Glob,Ab]\ar@{->}[rr]_{L_*(H)} 
&& Ab\\
&&&&&&\\
{Sets^{\Gamma^{op}}}
&& \ar@{_{(}->}[ll]{Sets^{\Gamma^{op}}_1}\ar[dddd]_{\mathcal{N}^\alpha}
&& {Comp(Ab)}\ar[uurr]_{H_*}
&&\\&&&&&&\\
&&&&&&\\&&&&&&\\
&&
{Sets^{\Delta^{op}}}\ar[rr]_{N}
&& {Comp(Ab)} \ar[rr]_{H_*}&& Ab\\}
\]

The above diagram is commutative in the sense that two different ways
between the same pair of points give the same transformation. All
these transformations are functors except $Gl(-)$. This diagram
summarizes all transformations or functors constructed in this paper.

\listoffigures

%\nocite{*}
%\bibliographystyle{alpha}
%\bibliography{article}

\end{document}